\documentclass[aos,preprint]{imsart}

\RequirePackage{amsthm,amsmath,amsfonts,amssymb}
\RequirePackage{mathtools}
\RequirePackage{mathrsfs}
\RequirePackage{bm}

\RequirePackage{graphicx}
\RequirePackage{array}
\RequirePackage{booktabs}
\RequirePackage{makecell}
\RequirePackage{multirow}
\RequirePackage{longtable}

\RequirePackage{enumitem}

\RequirePackage[authoryear]{natbib}

\RequirePackage[colorlinks,citecolor=blue,urlcolor=blue,linkcolor=blue]{hyperref}

\usepackage{amsmath}

\startlocaldefs

\numberwithin{equation}{section}

\theoremstyle{plain}
\newtheorem{theorem}{Theorem}
\newtheorem{lemma}{Lemma}
\newtheorem{proposition}{Proposition}
\newtheorem{corollary}{Corollary}

\theoremstyle{definition}
\newtheorem{assumption}{Assumption}
\newtheorem{remark}{Remark}

\makeatletter
\def\maxwidth{\ifdim\Gin@nat@width>\linewidth\linewidth\else\Gin@nat@width\fi}
\def\maxheight{\ifdim\Gin@nat@height>\textheight\textheight\else\Gin@nat@height\fi}
\makeatother
\setkeys{Gin}{width=\maxwidth,height=\maxheight,keepaspectratio}

\endlocaldefs

\begin{document}

\begin{frontmatter}

\title{Fairness-Regularized Principal Component Analysis}
\runtitle{Fair PCA}

\begin{aug}

\author[A]{\fnms{Chuang}~\snm{Xu}\ead[label=e1]{Chuang.Xu@anu.edu.au}}
\author[A]{\fnms{Yanrong}~\snm{Yang}\ead[label=e2]{Yanrong.Yang@anu.edu.au}}
\author[A]{\fnms{Andrew T. A.}~\snm{Wood}\ead[label=e3]{Andrew.Wood@anu.edu.au}}

\address[A]{Research School of Finance, Actuarial Studies and Statistics,
Australian National University\printead[presep={,\\ }]{e1,e2,e3}}

\end{aug}

\begin{abstract}
Principal component analysis (PCA) optimizes aggregate reconstruction accuracy but may represent heterogeneous groups unevenly. We study fairness-regularized PCA, which balances aggregate reconstruction risk against worst-group reconstruction regret, and establish its population and statistical properties. First, we characterize the reconstruction–fairness trade-off and show that, under a nondegeneracy condition, a small fairness penalty reduces worst-group regret to first order while increasing aggregate reconstruction risk only to second order. This reveals a vanishing initial marginal cost of improving fairness. We also identify conditions under which the regularization path remains locally constant. A minimax spectral characterization shows that, under an optimal effective eigengap condition, the convex relaxation is exact and Fair PCA is equivalent to ordinary PCA with endogenously adjusted group weights. Second, we develop a statistical theory for estimating fairness-regularized principal subspaces. A uniform eigengap over the set of possibly nonunique dual optimizers implies polyhedral identification and quadratic growth, yielding first-order perturbation bounds and convergence rates that accommodate growing dimensions, ranks, and numbers of groups. Third, we establish uniform recovery of the reconstruction–fairness path and theoretical guarantees for data-adaptive penalty selection. Simulations and an NHANES accelerometry application demonstrate substantial improvements in worst-group representation quality with minimal loss of aggregate reconstruction accuracy.


\end{abstract}



\end{frontmatter}

\section{Introduction}
\label{sec:fair-introduction}
Principal component analysis (PCA) provides a common
low-dimensional representation by maximizing aggregate
explained variance, or equivalently minimizing average
squared reconstruction error. When observations arise
from heterogeneous groups, however, a subspace that
is optimal on average will typically not represent each group
equally well. Differences in group covariance structures
may cause a shared principal subspace to retain important
directions of variation for some groups while poorly
approximating those of others. This raises two closely related questions:
how should group-specific losses in
representation quality be controlled,
and at what cost to aggregate reconstruction
accuracy can these losses be reduced?

Statistical notions of fairness depend on the learning task and the aspect of performance compared across groups \citep{HeLi2025,CherianCandes2024,BarocasEtAl2023}. Existing approaches examine, for example, demographic parity in regression, which requires predicted outcomes to have the same distribution across groups \citep{ChzhenSchreuder2022}; group-wise error control in classification, which controls specified error rates within each group \citep{RavaEtAl2026}, and equal-opportunity constraints, which limit disparities in false-negative rates across
groups \citep{FanEtAl2026}; counterfactual fairness through preprocessing, which seeks predictions that remain unchanged under a hypothetical change in a protected attribute, such as gender or race \citep{ChenEtAl2024}; and fairness--utility trade-offs in representation learning, which balance predictive usefulness against the encoding of protected information \citep{XuStrohmer2023}.

Our focus is the reconstruction fidelity of a common unsupervised linear representation. Even within PCA, fairness admits different criteria. One approach controls disparities in absolute reconstruction errors across groups. However, groups with different covariance structures may have different optimal rank-\(r\) reconstruction risks. Consequently, disparities in raw errors need not indicate unequal representation quality relative to what each group could achieve individually. We instead consider reconstruction regret: each group's excess reconstruction loss relative to its own optimal rank-\(r\) principal subspace. Controlling the largest regret across groups limits the worst additional loss incurred by adopting a common representation in place of a group-specific optimum. Ordinary PCA and minimax-regret PCA represent
two different priorities.
While ordinary PCA minimizes aggregate
reconstruction risk, it may leave some groups
with substantial excess reconstruction loss.
Conversely, minimizing worst-group regret
provides direct protection against the largest
group-specific loss, but may entail a substantial
increase in aggregate reconstruction risk.
In applications where both aggregate accuracy
and worst-group representation quality matter,
it is therefore natural to consider the
trade-off between these objectives rather
than optimizing either one in isolation.

The reconstruction-regret perspective originates in the fair PCA formulation of \citet{SamadiEtAl2018}. Subsequent studies extended this framework through multi-criteria formulations and computational algorithms for fairness-aware dimensionality reduction \citep{TantipongpipatEtAl2019,KamaniEtAl2022}. Related approaches have addressed disparities in groups' absolute reconstruction errors through penalized and eigendecomposition-based formulations \citep{VuEtAl2022,PelegrinaDuarte2024}. Focusing on minimax reconstruction regret, \citet{ShenEtAl2026} derived an eigenvalue-optimization characterization for the two-group setting, while \citet{MatakosEtAl2025} investigated sequential constructions of fair principal components that preserve nested subspaces. Other worst-case PCA formulations have also been studied \citep{ZalcbergWiesel2021}. More recently, \citet{WangEtAl2026StablePCA} developed a Mirror-Prox algorithm with statistical and optimization error bounds for worst-case explained-variance PCA and extended this computational approach to minimax-regret Fair PCA. From a statistical perspective, \citet{FriesEtAl2026} established consistency of empirical solutions to several worst-case PCA formulations, including reconstruction-regret-based PCA, under uniqueness of the population solution. A distinct line of research considers alternative fairness criteria based on sensitive-attribute information or discrepancies between projected group distributions \citep{OlfatAswani2019,LeeEtAl2022MMD,KleindessnerEtAl2023}.

Despite these advances, a comprehensive characterization of the population reconstruction–fairness trade-off and the statistical estimation of principal subspaces along this trade-off remain less understood, particularly in increasing-dimensional settings. 
In particular, we seek to quantify how
worst-group regret responds to departures
from ordinary PCA, the associated cost in
aggregate reconstruction accuracy, and how accurately the resulting principal
subspaces can be estimated from data.

Specifically, we develop a population and empirical framework for PCA regularized by worst-group reconstruction regret. For a common rank-$r$ orthogonal projector $P$, let $R(P)$ denote the weighted aggregate reconstruction risk, and let $R_a(P)$ and $R_a^\star$ denote the reconstruction risk under $P$ and the optimal rank-$r$ reconstruction risk for group $a$, respectively. We study the criterion
\[
R(P)+\lambda B(P),
\qquad
B(P)=\max_{1\leq a\leq G}
{R_a(P)-R_a^\star},
\]
where $\lambda\geq 0$ controls the emphasis on worst-group reconstruction regret. When $\lambda=0$, the criterion reduces to ordinary aggregate PCA, while increasing $\lambda$ progressively emphasizes fairness, approaching minimization of worst-group regret in the large-penalty limit. Building on a minimax spectral characterization of the regularized formulation, our main theoretical contributions are threefold. First, we characterize the local geometry of the reconstruction--fairness trade-off. Second, we establish first-order perturbation bounds for estimated principal subspaces that accommodate growing ambient dimensions, target ranks, and numbers of groups, and recover root-$n$ convergence in fixed-dimensional settings, even when the optimal dual weights are nonunique. Third, we establish uniform recovery of the regularization path and provide theoretical guarantees for data-driven penalty selection.

As a foundation for these contributions, we begin by studying a convex relaxation of the regularized Fair PCA problem, replacing the nonconvex set of rank-$r$ orthogonal projectors with its convex hull, the rank-$r$ Fantope \citep{VuChoLeiRohe2013,LeiVu2015}. The relaxed problem admits a minimax spectral representation in which worst-group reconstruction regret is expressed through nonnegative group weights summing to one. These weights, referred to as dual weights, are determined by an auxiliary optimization problem and, together with the fairness penalty, induce an adjustment to the usual weighted average of group covariance matrices. We refer to the resulting matrix as the effective covariance. We establish that if the effective covariance associated with at least one dual-optimal weight vector has a strict rank-$r$ eigengap, then the convex relaxation is exact and admits a unique solution given by the leading rank-$r$ spectral projector of that matrix. Consequently, regularized Fair PCA is equivalent to ordinary PCA applied to an endogenously reweighted covariance matrix. Importantly, the optimal principal subspace remains uniquely identified even when the dual-optimal weights are nonunique. This spectral characterization provides the foundation for our analysis of the population regularization path and statistical subspace recovery.

Our first main theoretical contribution is to characterize the population reconstruction--fairness trade-off along the regularization path. We show that worst-group reconstruction regret is nonincreasing, whereas aggregate reconstruction risk is nondecreasing, as the fairness penalty increases. We further establish a local expansion around ordinary PCA that quantifies the initial cost of improving fairness. Under a nondegeneracy condition, a small positive penalty reduces worst-group regret at first order while increasing aggregate reconstruction risk only at second order. Thus, the reconstruction cost is quadratic in the initial reduction in worst-group regret, implying a vanishing marginal reconstruction cost of fairness at ordinary PCA. We also identify a complementary regime in which the regularization path is locally constant: when the relevant nondegeneracy condition fails, the population-optimal principal subspace coincides with the ordinary PCA subspace for all sufficiently small penalties. Together, these results distinguish settings in which fairness can be improved at negligible initial reconstruction cost from those in which fairness regularization leaves the optimal representation unchanged locally.

Our second main theoretical contribution concerns statistical estimation of fairness-regularized principal subspaces. High-dimensional PCA and spectral projector estimation have been extensively studied, including conditions for consistency and quantitative perturbation bounds \citep{JohnstoneLu2009,KoltchinskiiLounici2017Spectral,YuWangSamworth2015}. Unlike ordinary PCA, however, the effective covariance in our framework depends on data-driven dual weights that may not be uniquely identified at the population level. We first establish uniform convergence of the empirical criterion and consistency of empirical Fantope optimizers towards the unique population principal projector, under a strict rank-\(r\) eigengap at at least one population dual optimum. To obtain sharper guarantees, we analyze the geometry of the entire dual-optimal set and establish a quadratic-growth property for the dual objective, ensuring that near-optimal weights lie close to the set of population-optimal weights rather than to any particular optimizer. Under a stronger uniform eigengap condition, this geometric characterization yields first-order projector perturbation bounds in terms of covariance estimation error, without requiring uniqueness of the dual weights. The resulting convergence guarantees explicitly capture the effects of spectral separation and dual-optimal geometry, accommodate growing ambient dimensions, target ranks, and numbers of groups, and recover root-$n$ rates in fixed-dimensional settings under standard covariance concentration conditions.

Our third main theoretical contribution concerns uniform statistical recovery of the population reconstruction--fairness path over compact intervals of positive penalties, covering both the principal subspaces and their aggregate reconstruction risks and worst-group regrets. These results provide a theoretical foundation for data-adaptive penalty selection according to a prescribed balance between reconstruction accuracy and fairness. Under suitable conditions, the resulting fitted representation achieves asymptotically optimal population performance relative to the regularization path. We further study penalty calibration based on a specified tolerance for worst-group reconstruction regret and derive a finite-sample bound on violations of the population fairness tolerance. Simulations under heterogeneous covariance models and an application to daily accelerometry profiles from the National Health and Nutrition Examination Survey illustrate the theoretical findings and the empirical reconstruction--fairness trade-off, including its performance on held-out observations.


The remainder of the paper is organized as follows.
Section~\ref{sec:fair-pca} introduces the formulation
and its spectral characterization.
Section~\ref{sec:fair-tradeoff} studies the population
reconstruction--fairness path, and
Section~\ref{sec:fair-sample-theory} develops the
statistical theory. Penalty selection is considered in
Section~\ref{sec:fair-selection}.
Sections~\ref{sec:fair-simulations} and
\ref{sec:fair-nhanes} present the numerical studies and
real-data application. Proofs and additional numerical
results are given in the Supplementary Material.

\section{Fair PCA: Formulation and Spectral Characterization}
\label{sec:fair-pca}

\subsection{Group-specific reconstruction regret}
\label{subsec:fair-pca-setup}

Consider a population comprising $G\geq 2$ groups. For each
$a\in\{1,\ldots,G\}$, let $X_a\in\mathbb R^p$ be a random vector
representing an observation from group $a$, with
\(
\mathbb E(X_a)=0,
\operatorname{Cov}(X_a)=\Sigma_a.
\)
Throughout, $1\leq r<p$ denotes the dimension of the principal
subspace of interest. Let $q_a\geq 0$, with
$\sum_{a=1}^G q_a=1$, be prespecified group weights, and define
the aggregate covariance matrix
\[
\overline\Sigma=\sum_{a=1}^G q_a\Sigma_a.
\]
The weights may represent population proportions or other
application-specific priorities and need not coincide with the
relative group sample sizes.
Let
\[
\mathcal P_r
=
\left\{
P\in\mathbb R^{p\times p}:
P=P^\top=P^2,\ \operatorname{tr}(P)=r
\right\}
\]
denote the set of rank-$r$ orthogonal projectors. For
$P\in\mathcal P_r$, the reconstruction risk in group $a$ is
\(
R_a(P)
=
\mathbb E\|X_a-PX_a\|^2
=
\operatorname{tr}\{(I-P)\Sigma_a\}.
\)
The corresponding aggregate reconstruction risk is
\begin{equation}
R(P)
=
\sum_{a=1}^G q_a R_a(P)
=
\operatorname{tr}\{(I-P)\overline\Sigma\}.
\label{eq:fair-aggregate-risk}
\end{equation}

For a symmetric matrix $A$, let
$\theta_1(A)\geq\cdots\geq\theta_p(A)$ denote its ordered
eigenvalues, and write
\(
\mathcal K_r(A)
=
\sum_{j=1}^r\theta_j(A)
\)
for its rank-$r$ Ky Fan functional. Ordinary population PCA
minimizes $R(P)$ over $\mathcal P_r$. In particular, if $P_r(A)$
denotes a projector onto a leading $r$-dimensional eigenspace
of $A$, an ordinary PCA solution is
\(
P_0^\star=P_r(\overline\Sigma).
\)
This projector is uniquely determined whenever
$\theta_r(\overline\Sigma)>
\theta_{r+1}(\overline\Sigma)$.

Ordinary PCA optimizes aggregate reconstruction accuracy without
explicitly controlling the reconstruction performance of individual
groups. To quantify the loss incurred by a group when using a
common representation, consider its group-specific optimal
rank-$r$ reconstruction risk,
\(
R_a^\star
=
\min_{P\in\mathcal P_r}R_a(P)
=
\operatorname{tr}(\Sigma_a)-\mathcal K_r(\Sigma_a).
\)
Writing
\(
c_a=\mathcal K_r(\Sigma_a),
\)
we define the reconstruction regret of group $a$ under $P$ by
\begin{equation}
\mathcal E_a(P)
=
R_a(P)-R_a^\star
=
c_a-\operatorname{tr}(P\Sigma_a).
\label{eq:fair-group-regret}
\end{equation}
The use of group-specific excess reconstruction
loss as a fairness criterion was introduced by
\citet{SamadiEtAl2018}. Here we investigate the
population and statistical properties of PCA
regularized by the largest such loss.
By the Ky Fan variational principle,
$\mathcal E_a(P)\geq 0$ for every $P\in\mathcal P_r$.

Unlike raw reconstruction risk, $\mathcal E_a(P)$ measures the
additional loss incurred by group $a$ relative to its own
optimal rank-$r$ representation. This distinction is relevant
when groups have different covariance spectra and hence
different irreducible reconstruction errors. In particular,
equal reconstruction risks need not imply that groups are
equally well represented relative to their respective
attainable benchmarks.

We measure the fairness of a common representation through
its worst-group reconstruction regret,
\begin{equation}
B(P)
=
\max_{1\leq a\leq G}\mathcal E_a(P).
\label{eq:fair-worst-regret}
\end{equation}
A smaller value of $B(P)$ indicates that the largest
group-specific excess reconstruction loss is better controlled.
The proposed fairness criterion therefore protects the
worst-represented group in terms of regret, rather than
requiring equality of group reconstruction errors.

\subsection{Regularized Fair PCA and its Fantope relaxation}
\label{subsec:fair-pca-relaxation}

To balance aggregate reconstruction accuracy and worst-group
reconstruction regret, we consider the regularized Fair PCA
criterion
\begin{equation}
P_\lambda^\star
\in
\arg\min_{P\in\mathcal P_r}
\left\{
R(P)+\lambda B(P)
\right\},
\label{eq:fair-penalized-projector}
\end{equation}
where $\lambda\geq 0$ controls the strength of the fairness
penalty. When $\lambda=0$, the criterion reduces to ordinary
aggregate PCA. As $\lambda$ increases, greater emphasis is
placed on reducing the largest group-specific reconstruction
regret. The resulting family of solutions describes a
trade-off between aggregate reconstruction accuracy and
worst-group protection, which will be studied in
Section~\ref{sec:fair-tradeoff}.

The optimization in \eqref{eq:fair-penalized-projector}
is nonconvex because the set $\mathcal P_r$ of rank-$r$
orthogonal projectors is nonconvex. To obtain a convex formulation, we consider
the convex hull of $\mathcal P_r$, the
rank-$r$ Fantope, which has also been used
in the estimation of sparse principal
subspaces \citep{VuChoLeiRohe2013,LeiVu2015}:
\[
\mathcal F_r
=
\left\{
H\in\mathbb R^{p\times p}:
H=H^\top,\ 
0\preceq H\preceq I,\ 
\operatorname{tr}(H)=r
\right\}.
\]
Here $A\preceq B$ means that $B-A$
is positive semidefinite. For $H\in\mathcal F_r$, extend the aggregate reconstruction
risk and group-specific regret by
\[
R(H)
=
\operatorname{tr}\{(I-H)\overline\Sigma\},
\qquad
\mathcal E_a(H)
=
c_a-\operatorname{tr}(H\Sigma_a),
\]
and define
\(
B(H)
=
\max_{1\leq a\leq G}\mathcal E_a(H).
\)
These extensions agree with the original criteria on
$\mathcal P_r$. Moreover, the Ky Fan variational principle
ensures that $\mathcal E_a(H)\geq 0$ for every
$H\in\mathcal F_r$.
We define the relaxed regularized Fair PCA oracle by
\begin{equation}
H_\lambda^\star
\in
\arg\min_{H\in\mathcal F_r}
\left\{
R(H)+\lambda B(H)
\right\}.
\label{eq:fair-fantope-oracle}
\end{equation}
Since $R$ is affine and $B$ is the pointwise maximum of
finitely many affine functions, the objective in
\eqref{eq:fair-fantope-oracle} is convex. The Fantope is
compact and convex, so a minimizer exists for every
$\lambda\geq 0$. A complementary formulation minimizes aggregate
reconstruction risk subject to a prescribed upper
bound on worst-group regret:
\[
\min_{H\in\mathcal F_r} R(H)
\quad\text{subject to}\quad
B(H)\leq\tau,
\]
where $\tau$ denotes the fairness tolerance.
Under suitable feasibility conditions, convex
duality relates this constrained formulation
to the regularized problem
\eqref{eq:fair-fantope-oracle}.
We focus on the regularized formulation and
return to fairness-tolerance calibration in
Section~\ref{subsec:fair-selection-tolerance}.

Although the Fantope is the convex hull of $\mathcal P_r$,
a minimizer of \eqref{eq:fair-fantope-oracle} need not be a
rank-$r$ projector. Thus the convex relaxation is not
automatically exact. In the following subsections, we
derive a minimax representation of the relaxed problem
and establish a spectral condition under which its
minimizer is unique and belongs to $\mathcal P_r$, thereby
recovering a solution of the original Fair PCA problem.

\subsection{Minimax duality and effective covariance}
\label{subsec:fair-pca-duality}
Let
\begin{equation}
\Delta_G
=
\left\{
\alpha\in\mathbb R^G:
\alpha_a\geq 0,\quad
\sum_{a=1}^G\alpha_a=1
\right\}
\label{eq:probability-simplex}
\end{equation}
denote the probability simplex. Since the maximum of
finitely many real numbers equals the maximum of their
convex combinations,
\[
B(H)
=
\max_{\alpha\in\Delta_G}
\sum_{a=1}^G\alpha_a\mathcal E_a(H).
\]
The relaxed Fair PCA criterion consequently admits the
minimax representation
\begin{equation}
\min_{H\in\mathcal F_r}
\max_{\alpha\in\Delta_G}
L_\lambda(H,\alpha),
\label{eq:fair-minimax}
\end{equation}
where
\(
L_\lambda(H,\alpha)
=
R(H)+
\lambda\sum_{a=1}^G
\alpha_a\mathcal E_a(H).
\)
Define the fairness-adjusted effective covariance matrix
\begin{equation}
M_\lambda(\alpha)
=
\overline\Sigma+
\lambda\sum_{a=1}^G\alpha_a\Sigma_a,
\qquad
\alpha\in\Delta_G.
\label{eq:fair-effective-covariance}
\end{equation}
Substitution into the minimax criterion gives
\begin{equation}
L_\lambda(H,\alpha)
=
\operatorname{tr}(\overline\Sigma)
+
\lambda\sum_{a=1}^G\alpha_a c_a
-
\operatorname{tr}\{H M_\lambda(\alpha)\}.
\label{eq:fair-lagrangian}
\end{equation}
Since $L_\lambda$ is continuous and affine in each argument
and both feasible sets are compact and convex, the minimax
theorem yields
\begin{equation}
\min_{H\in\mathcal F_r}
\max_{\alpha\in\Delta_G}L_\lambda(H,\alpha)
=
\max_{\alpha\in\Delta_G}
\min_{H\in\mathcal F_r}L_\lambda(H,\alpha).
\label{eq:fair-minimax-equality}
\end{equation}
For a fixed $\alpha$, the inner minimization is equivalent
to maximizing $\operatorname{tr}\{H M_\lambda(\alpha)\}$
over $\mathcal F_r$. By the Ky Fan variational principle,
\(
\max_{H\in\mathcal F_r}
\operatorname{tr}\{H M_\lambda(\alpha)\}
=
\mathcal K_r\{M_\lambda(\alpha)\}.
\)
We therefore define the dual objective
\begin{equation}
\phi_\lambda(\alpha)
=
\lambda\sum_{a=1}^G\alpha_a c_a
-
\mathcal K_r\{M_\lambda(\alpha)\}.
\label{eq:fair-dual-objective}
\end{equation}
Writing
\(
V(\lambda)
=
\min_{H\in\mathcal F_r}
\left\{
R(H)+\lambda B(H)
\right\}
\)
for the optimal value of the relaxed Fair PCA problem,
we obtain
\begin{equation}
V(\lambda)
=
\operatorname{tr}(\overline\Sigma)
+
\max_{\alpha\in\Delta_G}\phi_\lambda(\alpha).
\label{eq:fair-primal-dual-value}
\end{equation}
Let
\begin{equation}
\mathcal A_\lambda^\star
=
\arg\max_{\alpha\in\Delta_G}
\phi_\lambda(\alpha)
\label{eq:fair-dual-optimal-set}
\end{equation}
denote the set of optimal dual weights.

The Ky Fan functional is convex, whereas
$M_\lambda(\alpha)$ is affine in $\alpha$. Hence
$\phi_\lambda$ is concave on $\Delta_G$, and the dual
problem is a concave maximization problem over a
$(G-1)$-dimensional simplex. The dual-optimal set
$\mathcal A_\lambda^\star$ is nonempty, compact and
convex, but need not be a singleton.
\subsection{Spectral characterization and rank recovery}
\label{subsec:fair-pca-spectral}

We now establish conditions under which the convex
relaxation yields a unique rank-$r$ orthogonal projector.

\begin{assumption}[Optimal effective eigengap]
\label{ass:fair-optimal-eigengap}
For the prescribed $\lambda\geq 0$, there exists
$\alpha^\star\in\mathcal A_\lambda^\star$ such that
\begin{equation}
\theta_r\{M_\lambda(\alpha^\star)\}
>
\theta_{r+1}\{M_\lambda(\alpha^\star)\}.
\label{eq:fair-optimal-eigengap}
\end{equation}
\end{assumption}

\begin{theorem}[Spectral characterization of Fair PCA]
\label{thm:fair-spectral-characterization}
Suppose Assumption~\ref{ass:fair-optimal-eigengap} holds.
Then the relaxed problem
\eqref{eq:fair-fantope-oracle} has a unique minimizer,
which is a rank-$r$ orthogonal projector satisfying
\begin{equation}
H_\lambda^\star
=
P_\lambda^\star
=
P_r\{M_\lambda(\alpha^\star)\}.
\label{eq:fair-spectral-projector}
\end{equation}
Consequently, the Fantope relaxation is exact and
$P_\lambda^\star$ is also the unique minimizer of the
original projector-level problem
\eqref{eq:fair-penalized-projector}.

Furthermore, uniqueness of the primal solution does
not require uniqueness of the dual weights. For every
$\alpha\in\mathcal A_\lambda^\star$,
\begin{equation}
P_\lambda^\star
\in
\arg\max_{H\in\mathcal F_r}
\operatorname{tr}\{H M_\lambda(\alpha)\}.
\label{eq:fair-common-primal-projector}
\end{equation}
In particular, any dual-optimal weight satisfying a
strict rank-$r$ eigengap induces the same projector
$P_\lambda^\star$.
\end{theorem}

The theorem identifies Fair PCA with ordinary PCA
applied to a covariance matrix that is adaptively
adjusted to account for group-specific reconstruction
regret. Indeed,
\[
\begin{aligned}
M_\lambda(\alpha^\star)
=
\overline\Sigma+
\lambda\sum_{a=1}^G\alpha_a^\star\Sigma_a
=
(1+\lambda)
\sum_{a=1}^G
w_{\lambda,a}^\star\Sigma_a,
\end{aligned}
\]
where
\begin{equation}
w_{\lambda,a}^\star
=
\frac{q_a+\lambda\alpha_a^\star}{1+\lambda},
\qquad
a=1,\ldots,G.
\label{eq:fair-adaptive-weights}
\end{equation}
The weights are nonnegative and sum to one. Since
multiplication by a positive scalar does not change
eigenspaces, $P_\lambda^\star$ is the leading rank-$r$
spectral projector of
\(
\sum_{a=1}^G
w_{\lambda,a}^\star\Sigma_a.
\)
Thus Fair PCA can be viewed as an adaptive reweighting
of group covariance matrices, with the reweighting
determined endogenously by the worst-group regret
criterion.

For $\lambda>0$, saddle-point optimality also gives
\begin{equation}
\alpha_a^\star>0
\quad\Longrightarrow\quad
\mathcal E_a(P_\lambda^\star)
=
B(P_\lambda^\star).
\label{eq:fair-active-groups}
\end{equation}
Hence only groups attaining the worst reconstruction
regret can receive positive adversarial weight.
The converse need not hold: a group may attain
the maximum regret without receiving positive
weight in a particular dual optimum.

The dual representation also facilitates computation. Whenever
$M_\lambda(\alpha)$ has a strict rank-$r$ eigengap, the dual
objective is differentiable at $\alpha$, with
\begin{equation}
\frac{\partial\phi_\lambda(\alpha)}
{\partial\alpha_a}
=
\lambda
\left[
c_a-
\operatorname{tr}\{
P_r(M_\lambda(\alpha))\Sigma_a
\}
\right]
=
\lambda\,
\mathcal E_a\left(
P_r\{M_\lambda(\alpha)\}
\right),
\qquad a=1,\ldots,G.
\label{eq:fair-dual-gradient}
\end{equation}
Thus the dual gradient is determined by the group-specific
reconstruction regrets associated with the leading eigenspace
of the effective covariance.
More generally, differentiability is not required to obtain
first-order information. For any
\(
H_\alpha
\in
\arg\max_{H\in\mathcal F_r}
\operatorname{tr}\{H M_\lambda(\alpha)\},
\)
the vector $g(\alpha)\in\mathbb R^G$ with components
\begin{equation}
g_a(\alpha)
=
\lambda\mathcal E_a(H_\alpha),
\qquad a=1,\ldots,G,
\label{eq:fair-dual-supergradient}
\end{equation}
is a supergradient of the concave dual objective
$\phi_\lambda$ at $\alpha$.
Consequently, the dual problem can be solved using
first-order optimization methods over $\Delta_G$,
with gradient or supergradient information obtained
from leading-eigenspace computations.

When Assumption~\ref{ass:fair-optimal-eigengap}
does not hold, a relaxed optimizer need not be a rank-$r$ projector.
In particular, an arbitrary leading rank-$r$ projector
of a dual-optimal effective covariance need not solve
the original Fair PCA problem.

The representation
\eqref{eq:fair-spectral-projector} will be the basis
for our analysis of the population reconstruction--
fairness trade-off and the statistical properties
of the empirical Fair PCA estimator.

\section{Population Reconstruction--Fairness Trade-off}
\label{sec:fair-tradeoff}
We investigate the trade-off between aggregate reconstruction
risk and worst-group reconstruction regret induced by the
regularized Fair PCA criterion. We first establish global
properties of the regularization path and its marginal-price
interpretation. We then characterize the local behavior
around ordinary PCA, showing that, under a nondegeneracy
condition, worst-group regret decreases to first order
while aggregate reconstruction risk increases only to
second order.

The global results are stated for the Fantope formulation
\eqref{eq:fair-fantope-oracle} and do not require uniqueness
of the optimizer. Whenever the optimal effective eigengap
condition holds, they also describe the corresponding
rank-$r$ Fair PCA projector.
\subsection{Monotonicity and the price of fairness}
\label{subsec:fair-monotonicity}

The regularized Fair PCA solutions exhibit a monotone
reconstruction--fairness trade-off. For
$0\leq\lambda_1<\lambda_2$, let
$H_{\lambda_1}^\star$ and $H_{\lambda_2}^\star$
be any corresponding relaxed optimizers. Then
\[
B(H_{\lambda_2}^\star)
\leq B(H_{\lambda_1}^\star),
\qquad
R(H_{\lambda_2}^\star)
\geq R(H_{\lambda_1}^\star).
\]
Thus increasing the fairness penalty cannot increase
worst-group regret or decrease aggregate reconstruction
risk. 
These monotonicity properties follow directly from the
optimality of the regularized solutions at different
penalty levels. More precise finite-difference bounds
are provided in the Supplementary Material.

The optimal value function $V$, defined in
Section~\ref{subsec:fair-pca-duality}, provides a
complementary interpretation of this trade-off.
\begin{proposition}[Value function and marginal price]
\label{prop:fair-value-function}
The function $V$ is nondecreasing and concave on
$[0,\infty)$. For every $\lambda\geq0$, every optimizer
$H_\lambda^\star$ of \eqref{eq:fair-fantope-oracle},
and every $\lambda'\geq0$,
\begin{equation}
V(\lambda')
\leq
V(\lambda)
+
(\lambda'-\lambda)B(H_\lambda^\star).
\label{eq:fair-value-supergradient}
\end{equation}
Consequently, $B(H_\lambda^\star)$ is a supergradient
of $V$ at $\lambda$. At every differentiability point
$\lambda>0$,
\begin{equation}
V'(\lambda)=B(H_\lambda^\star).
\label{eq:fair-value-derivative}
\end{equation}
Furthermore, for any selection of optimizers along the
regularization path, the functions
\[
\lambda\longmapsto R(H_\lambda^\star),
\qquad
\lambda\longmapsto B(H_\lambda^\star)
\]
are respectively nondecreasing and nonincreasing.
At almost every $\lambda>0$,
\begin{equation}
\frac{d}{d\lambda}R(H_\lambda^\star)
=
-\lambda
\frac{d}{d\lambda}B(H_\lambda^\star).
\label{eq:fair-marginal-price}
\end{equation}
\end{proposition}

Equation~\eqref{eq:fair-marginal-price} identifies
$\lambda$ as the marginal reconstruction price of
reducing worst-group regret. At points where the regret
decreases at a strictly negative derivative, an
infinitesimal reduction in worst-group regret requires
an increase in aggregate reconstruction risk at a
rate equal to $\lambda$. The penalty parameter thus
quantifies the local exchange rate between the two
criteria, rather than directly specifying an attained
fairness level.
\subsection{Local geometry around ordinary PCA}
\label{subsec:fair-local-geometry}

We now examine the behavior of the regularized oracle
as $\lambda\downarrow0$. In contrast to the preceding
global results, the local analysis uses the spectral
structure of PCA to obtain precise expansions for
both reconstruction risk and worst-group regret.
Let
\[
\overline\Sigma
=
\sum_{i=1}^p\mu_i u_i u_i^\top,
\qquad
\mu_1\geq\cdots\geq\mu_p,
\]
be a spectral decomposition of the aggregate covariance
matrix, where $\mu_i=\theta_i(\overline\Sigma)$ and
$u_1,\ldots,u_p$ are corresponding orthonormal
eigenvectors.
We assume that ordinary PCA has a strict rank-$r$
eigengap,
\begin{equation}
\delta_0
=
\mu_r-\mu_{r+1}>0.
\label{eq:fair-ordinary-eigengap}
\end{equation}
Hence the ordinary PCA projector
\(
P_0^\star
=
\sum_{i=1}^r u_i u_i^\top
\)
is uniquely determined.
Define the set of groups attaining the largest
reconstruction regret under ordinary PCA by
\begin{equation}
\mathcal I_0
=
\left\{
a\in\{1,\ldots,G\}:
\mathcal E_a(P_0^\star)=B(P_0^\star)
\right\}.
\label{eq:fair-ordinary-active-set}
\end{equation}
The corresponding face of the probability simplex is
\begin{equation}
\Delta(\mathcal I_0)
=
\left\{
\alpha\in\Delta_G:
\alpha_a=0
\ \text{for }a\notin\mathcal I_0
\right\}.
\label{eq:fair-active-simplex}
\end{equation}
For $\alpha\in\Delta_G$, write
\[
S(\alpha)
=
\sum_{a=1}^G\alpha_a\Sigma_a,
\]
so that the effective covariance
\eqref{eq:fair-effective-covariance} becomes
\(
M_\lambda(\alpha)
=
\overline\Sigma+\lambda S(\alpha).
\)
The interaction between the fairness perturbation
and the ordinary principal subspace is captured by
\begin{equation}
\Gamma_0(\alpha)
=
\sum_{i=1}^r
\sum_{j=r+1}^p
\frac{
\{u_j^\top S(\alpha)u_i\}^2
}{
\mu_i-\mu_j
}.
\label{eq:fair-local-spectral-coupling}
\end{equation}
The denominators are strictly positive by
\eqref{eq:fair-ordinary-eigengap}, and therefore
$\Gamma_0(\alpha)\geq0$.
We define the local fairness-response coefficient
\begin{equation}
\zeta_0
=
\min_{\alpha\in\Delta(\mathcal I_0)}
\Gamma_0(\alpha).
\label{eq:fair-local-coefficient}
\end{equation}

The minimum is attained by compactness. The quantities $u_j^\top S(\alpha)u_i$ measure the interactions induced by $S(\alpha)$ between the principal directions retained and discarded by ordinary PCA, which we refer to as cross-subspace interactions. Minimization over the active-group weights allows these interactions to offset one another across groups, while the eigenvalue gaps determine their relative contributions to the coefficient. In particular, $\zeta_0=0$ if and only if there exists $\alpha\in\Delta(\mathcal I_0)$ for which all cross-subspace interactions vanish. As shown below, $\zeta_0$ governs the leading-order changes in both worst-group reconstruction regret and aggregate reconstruction risk as $\lambda\downarrow0$.

\begin{theorem}[Local reconstruction--fairness expansion]
\label{thm:fair-local-expansion}
Suppose \eqref{eq:fair-ordinary-eigengap} holds.
Then there exists $\lambda_{\mathrm{loc}}>0$
such that, for every
$0\leq\lambda<\lambda_{\mathrm{loc}}$,
the relaxed Fair PCA oracle is unique and is a
rank-$r$ orthogonal projector.
As $\lambda\downarrow0$, the optimal value function
satisfies
\begin{equation}
V(\lambda)
=
R(P_0^\star)
+
\lambda B(P_0^\star)
-
\zeta_0\lambda^2
+
o(\lambda^2).
\label{eq:fair-local-value-expansion}
\end{equation}
Moreover, the worst-group reconstruction regret
and aggregate reconstruction risk satisfy
\begin{equation}
B(P_0^\star)-B(P_\lambda^\star)
=
2\zeta_0\lambda+o(\lambda)
\label{eq:fair-local-regret-expansion}
\end{equation}
and
\begin{equation}
R(P_\lambda^\star)-R(P_0^\star)
=
\zeta_0\lambda^2+o(\lambda^2),
\label{eq:fair-local-risk-expansion}
\end{equation}
respectively.
\end{theorem}

Theorem~\ref{thm:fair-local-expansion} establishes
an asymmetry between fairness improvement and
reconstruction cost near ordinary PCA. If
$\zeta_0>0$, a small fairness penalty produces
a first-order decrease in worst-group regret,
whereas the accompanying increase in aggregate
reconstruction risk is only second order.
Thus the initial marginal reconstruction price
of fairness vanishes.

When $\zeta_0>0$, eliminating $\lambda$ from the
expansions in Theorem~\ref{thm:fair-local-expansion}
further yields
\begin{equation}
\begin{aligned}
R(P_\lambda^\star)-R(P_0^\star)
={}&
\frac{
\{B(P_0^\star)-B(P_\lambda^\star)\}^2
}{
4\zeta_0
}+
o\left(
\{B(P_0^\star)-B(P_\lambda^\star)\}^2
\right).
\end{aligned}
\label{eq:fair-local-quadratic-frontier}
\end{equation}
Thus, locally around ordinary PCA, the aggregate
reconstruction cost grows quadratically with the
reduction in worst-group regret, with leading
coefficient $1/(4\zeta_0)$.
A larger $\zeta_0$ therefore corresponds to a
smaller reconstruction cost for a given small
reduction in worst-group regret.
We next characterize when $\zeta_0$ vanishes
and establish the resulting local constancy
of the population regularization path.
\begin{proposition}[Local nondegeneracy and stability]
\label{prop:fair-local-nondegeneracy}
Suppose \eqref{eq:fair-ordinary-eigengap} holds.
Then $\zeta_0=0$ if and only if there exists
$\alpha\in\Delta(\mathcal I_0)$ such that
\begin{equation}
(I-P_0^\star)S(\alpha)P_0^\star=0.
\label{eq:fair-local-cross-block}
\end{equation}
Equivalently, $\zeta_0>0$ precisely when no convex
combination of the active group covariance matrices
has a vanishing cross-block relative to the ordinary
principal subspace.

Furthermore, if $\zeta_0=0$, there exists
$\lambda_{\mathrm{flat}}>0$ such that
\begin{equation}
P_\lambda^\star=P_0^\star,
\qquad
0\leq\lambda<\lambda_{\mathrm{flat}}.
\label{eq:fair-local-flat-path}
\end{equation}
\end{proposition}

The proposition distinguishes two possible local
behaviors. When $\zeta_0>0$, the fairness penalty
immediately changes the optimal reconstruction--fairness
trade-off. When $\zeta_0=0$, ordinary PCA remains
optimal for all sufficiently small penalty values.
In particular, a vanishing local coefficient
corresponds to an initial plateau of the
regularization path, rather than merely a
higher-order reduction in worst-group regret.

When the ordinary-PCA worst group is unique,
say $\mathcal I_0=\{a_0\}$, the coefficient simplifies to
\[
\zeta_0
=
\sum_{i=1}^r
\sum_{j=r+1}^p
\frac{
(u_j^\top\Sigma_{a_0}u_i)^2
}{
\mu_i-\mu_j
}.
\]
In this case, $\zeta_0>0$ if and only if
\(
(I-P_0^\star)\Sigma_{a_0}P_0^\star\neq0.
\)
For multiple active groups, the corresponding
cross-block terms may cancel under convex
combinations, so the positivity of $\zeta_0$
depends on the joint covariance geometry of
the active groups.
Two simple covariance examples illustrating
the cases $\zeta_0=0$ and $\zeta_0>0$
are provided in
Section~\ref{supp:subsec:fair-local-examples}
of the Supplement.

The preceding results characterize the
population behavior of Fair PCA. In the next
section, we investigate how accurately the
population oracle can be recovered from
empirical covariance matrices.
\section{Statistical Theory for Fair PCA}
\label{sec:fair-sample-theory}

We investigate estimation of the population Fair PCA
oracle from independent groupwise samples. We first
establish uniform convergence of the empirical criterion
and consistency of the estimated principal subspace
under the optimal effective eigengap condition of
Section~\ref{subsec:fair-pca-spectral}. We then obtain
a sharper projector rate by exploiting the geometry
of the population dual-optimal set. In particular,
we show that a uniform eigengap over this set
implies the quadratic-growth property needed for
first-order perturbation analysis, even when the
dual optimizer is not unique.

The results accommodate sequences of population
models in which the ambient dimension $p$, the
rank $r$, and the number of groups $G$ may increase
with sample size. Unless stated otherwise, these
quantities, the covariance matrices, and the penalty
parameter may depend on $n$, with the dependence
suppressed in the notation.

We use $\|\cdot\|_2$ for the Euclidean norm
of a vector, and $\|\cdot\|_{\mathrm{op}}$
and $\|\cdot\|_F$ for the matrix operator
and Frobenius norms, respectively. We also write
\(
\sigma
=
\max_{1\leq a\leq G}
\|\Sigma_a\|_{\mathrm{op}}
\)
for the largest population group covariance
operator norm.
The notation $O_p(\cdot)$ and $o_p(\cdot)$
has its usual stochastic meaning.

\subsection{Empirical Fair PCA and covariance concentration}
\label{subsec:fair-sample-setup}

For each group $a=1,\ldots,G$, suppose that
$X_{a,1},\ldots,X_{a,n_a}$ are independent
mean-zero random vectors in $\mathbb R^p$ with
covariance matrix $\Sigma_a$. Observations from
different groups are also assumed independent.
Write $n=\sum_{a=1}^G n_a$ and define
\(
\widehat\Sigma_a
=
\frac{1}{n_a}
\sum_{i=1}^{n_a}X_{a,i}X_{a,i}^\top,
\widehat{\overline\Sigma}
=
\sum_{a=1}^G q_a\widehat\Sigma_a,
\)
where the group weights $q_a$ are the
prespecified weights of Section~\ref{subsec:fair-pca-setup}.
The mean-zero assumption is made for notational
simplicity. When group means are unknown,
groupwise sample-centered covariance matrices
may be used instead. Under the same sub-Gaussian
condition on the centered observations, the
simultaneous covariance concentration rate
established below remains valid; see the
Supplementary Material.
The empirical counterparts of $c_a$ and
$\mathcal E_a(H)$ are
\(
\widehat c_a
=
\mathcal K_r(\widehat\Sigma_a),
\widehat{\mathcal E}_a(H)
=
\widehat c_a-\operatorname{tr}(H\widehat\Sigma_a).
\)
For $H\in\mathcal F_r$, define
\(
\widehat R(H)
=
\operatorname{tr}\{(I-H)\widehat{\overline\Sigma}\},
\widehat B(H)
=
\max_{1\leq a\leq G}
\widehat{\mathcal E}_a(H).
\)
The empirical Fair PCA estimator is any solution
\begin{equation}
\widehat H_\lambda
\in
\arg\min_{H\in\mathcal F_r}
\{\widehat R(H)+\lambda\widehat B(H)\}.
\label{eq:fair-sample-estimator}
\end{equation}
The empirical effective covariance and dual
objective are, respectively,
\[
\widehat M_\lambda(\alpha)
=
\widehat{\overline\Sigma}
+
\lambda\sum_{a=1}^G
\alpha_a\widehat\Sigma_a
\]
and
\[
\widehat\phi_\lambda(\alpha)
=
\lambda\sum_{a=1}^G\alpha_a\widehat c_a
-
\mathcal K_r\{\widehat M_\lambda(\alpha)\},
\qquad
\alpha\in\Delta_G.
\]
Thus $\widehat\phi_\lambda$ follows the same
constant-free convention as the population
dual objective $\phi_\lambda$.
Let
\[
\widehat{\mathcal A}_\lambda
=
\arg\max_{\alpha\in\Delta_G}
\widehat\phi_\lambda(\alpha).
\]
The population minimax characterization
carries over to the empirical problem.
In particular, whenever an empirical
dual optimizer $\widehat\alpha_\lambda$
satisfies a strict rank-$r$ eigengap,
the estimator is unique and obeys
\[
\widehat H_\lambda
=
P_r\{\widehat M_\lambda(\widehat\alpha_\lambda)\}.
\]
For the statistical analysis, write
\begin{equation}
e_n
=
\max_{1\leq a\leq G}
\|\widehat\Sigma_a-\Sigma_a\|_{\mathrm{op}}.
\label{eq:fair-covariance-error}
\end{equation}
We use a standard sub-Gaussian condition
to control this error.

\begin{assumption}[Sub-Gaussian observations]
\label{ass:fair-subgaussian}
There exists a constant $K<\infty$,
independent of $n$, such that
\[
\|u^\top X_{a,i}\|_{\psi_2}
\leq
K(u^\top\Sigma_a u)^{1/2}
\]
for every $u\in\mathbb R^p$, every group $a$,
and every observation $i$, where $\|\cdot\|_{\psi_2}$ denotes the
sub-Gaussian Orlicz norm (see, e.g.,
\citealp{vershynin2026highdimensional}).
\end{assumption}
For $\Sigma_a\neq0$, let
\[
\mathfrak r_a
=
\frac{\operatorname{tr}(\Sigma_a)}
{\|\Sigma_a\|_{\mathrm{op}}}
\]
denote its effective rank, with
$\mathfrak r_a=0$ when $\Sigma_a=0$.
Define the simultaneous covariance
concentration scale
\begin{equation}
\begin{aligned}
\eta_n
=
\max_{1\leq a\leq G}
\|\Sigma_a\|_{\mathrm{op}}
\Bigg\{
\sqrt{
\frac{\mathfrak r_a+1+\log(2G)}{n_a}
}
+
\frac{\mathfrak r_a+1+\log(2G)}{n_a}
\Bigg\}.
\end{aligned}
\label{eq:fair-covariance-rate}
\end{equation}
Under Assumption~\ref{ass:fair-subgaussian},
standard effective-rank covariance concentration
bounds (see, e.g.,
\citealp{vershynin2026highdimensional}),
combined with a union bound, yield
\begin{equation}
e_n=O_p(\eta_n).
\label{eq:fair-simultaneous-concentration}
\end{equation}
The logarithmic dependence on $G$ permits
the number of groups to increase with $n$.

For the statistical analysis, it is convenient
to remove the additive terms that do not depend
on $H$ from the population and empirical
reconstruction--fairness objectives. Define
\[
Q_\lambda(H)
=
R(H)+\lambda B(H)
-\operatorname{tr}(\overline\Sigma)
=
-\operatorname{tr}(H\overline\Sigma)
+\lambda B(H)
\]
and
\[
\widehat Q_\lambda(H)
=
\widehat R(H)+\lambda\widehat B(H)
-\operatorname{tr}(\widehat{\overline\Sigma})
=
-\operatorname{tr}(H\widehat{\overline\Sigma})
+\lambda\widehat B(H).
\]
These centered criteria have the same minimizers
as their respective original objectives.

\begin{proposition}[uniform convergence of the empirical objective]
\label{prop:fair-uniform-sample-criterion}
For every $\lambda\geq0$,
\begin{equation}
\sup_{H\in\mathcal F_r}
|\widehat Q_\lambda(H)-Q_\lambda(H)|
\leq
r(1+2\lambda)e_n.
\label{eq:fair-uniform-criterion-bound}
\end{equation}
Consequently, under
Assumption~\ref{ass:fair-subgaussian},
\[
\sup_{H\in\mathcal F_r}
|\widehat Q_\lambda(H)-Q_\lambda(H)|
=
O_p\{r(1+\lambda)\eta_n\}.
\]
The same deterministic bound
\[
\sup_{\alpha\in\Delta_G}
|\widehat\phi_\lambda(\alpha)
-\phi_\lambda(\alpha)|
\leq
r(1+2\lambda)e_n
\]
holds for the dual objectives.
\end{proposition}

The bounds depend explicitly on the rank $r$,
rather than on the ambient dimension $p$.
The dimension enters through the covariance
concentration scale $\eta_n$. In particular,
uniform criterion convergence holds whenever
$r(1+\lambda)\eta_n\to0$.

\subsection{Consistency under an optimal effective eigengap}
\label{subsec:fair-basic-consistency}

We first obtain a projector bound requiring
only the existence of one dual-optimal
effective covariance with a strict eigengap,
as in Theorem~\ref{thm:fair-spectral-characterization}.
Neither uniqueness of the dual weights nor
additional identification assumptions
are required.

\begin{theorem}[Basic projector estimation bound]
\label{thm:fair-basic-projector-rate}
Fix $\lambda\geq0$, possibly depending on $n$,
and suppose Assumption~\ref{ass:fair-optimal-eigengap}
holds. Choose
$\alpha_\lambda^\star\in\mathcal A_\lambda^\star$
satisfying \eqref{eq:fair-optimal-eigengap},
and write
\[
g_\lambda
=
\theta_r\{M_\lambda(\alpha_\lambda^\star)\}
-
\theta_{r+1}\{M_\lambda(\alpha_\lambda^\star)\}
>0.
\]
Then the population oracle
$P_\lambda^\star$ is uniquely determined.
Every empirical optimizer
$\widehat H_\lambda$ satisfies
\begin{equation}
\|\widehat H_\lambda-P_\lambda^\star\|_F
\leq
2
\left\{
\frac{r(1+2\lambda)e_n}{g_\lambda}
\right\}^{1/2}.
\label{eq:fair-basic-projector-bound}
\end{equation}
Under Assumption~\ref{ass:fair-subgaussian},
\begin{equation}
\|\widehat H_\lambda-P_\lambda^\star\|_F
=
O_p\left\{
\sqrt{
\frac{r(1+\lambda)\eta_n}{g_\lambda}
}
\right\}.
\label{eq:fair-basic-projector-rate}
\end{equation}
Consequently, if
\(
\frac{r(1+\lambda)\eta_n}{g_\lambda}
\longrightarrow0,
\)
then
\(
\|\widehat H_\lambda-P_\lambda^\star\|_F
\xrightarrow{p}0.
\)
\end{theorem}
The proof combines uniform criterion control
with the quadratic spectral margin induced
by the optimal effective eigengap. Specifically,
the population objective satisfies
\[
Q_\lambda(H)-Q_\lambda(P_\lambda^\star)
\geq
\frac{g_\lambda}{2}
\|H-P_\lambda^\star\|_F^2,
\qquad H\in\mathcal F_r.
\]
This margin holds even if the population
dual-optimal set contains multiple weight
vectors.

The resulting rate is an objective-based
$M$-estimation bound. For example, when
$r$, $\lambda$, and the eigengap are bounded
appropriately and $\eta_n=O(n^{-1/2})$,
it gives an $O_p(n^{-1/4})$ projector rate.
This is generally conservative relative
to the first-order spectral perturbation
behavior of ordinary PCA
\citep{YuWangSamworth2015,
KoltchinskiiLounici2017Spectral}.
We next show how the dual structure of
Fair PCA can recover a first-order rate.

\subsection{Sharper projector recovery from dual-set geometry}
\label{subsec:fair-sharp-rate}

The difficulty in sharpening
\eqref{eq:fair-basic-projector-rate} is that
the empirical effective covariance depends
on the estimated dual weights. Moreover,
the population dual optimizer need not be
unique even when the population principal
subspace is uniquely determined.

We therefore study identification of the
entire dual-optimal set
$\mathcal A_\lambda^\star$, rather than
the error relative to one arbitrarily
chosen dual optimizer. The central structural
result is that, under a uniform optimal-set
eigengap, the Fair PCA dual objective
automatically satisfies quadratic growth
away from its optimal set.

\subsubsection{Structural identification of the dual-optimal set}
Throughout this subsection, suppose $\lambda>0$.
Define the rank-$r$ eigengap of a symmetric
matrix $A$ by
\[
\operatorname{gap}_r(A)
=
\theta_r(A)-\theta_{r+1}(A).
\]

\begin{assumption}[Uniform optimal-set eigengap]
\label{ass:fair-uniform-optimal-gap}
There exists $\delta_\lambda>0$ such that
\begin{equation}
\inf_{\alpha\in\mathcal A_\lambda^\star}
\operatorname{gap}_r\{M_\lambda(\alpha)\}
\geq
2\delta_\lambda.
\label{eq:fair-uniform-optimal-gap}
\end{equation}
\end{assumption}

Unlike the basic eigengap condition,
Assumption~\ref{ass:fair-uniform-optimal-gap}
requires spectral separation at every
population dual optimizer. It does not
require the dual optimizer to be unique.
By Theorem~\ref{thm:fair-spectral-characterization},
all population dual-optimal effective
covariances share the same unique
rank-$r$ projector $P_\lambda^\star$.
\begin{remark}[Scope of the uniform eigengap condition]
Assumption~\ref{ass:fair-uniform-optimal-gap}
is stronger than
Assumption~\ref{ass:fair-optimal-eigengap},
which requires a strict eigengap at
only one population dual optimizer.
Under the latter assumption,
the population Fair PCA projector
$P_\lambda^\star$ is unique and is a
leading rank-$r$ spectral projector
of $M_\lambda(\alpha)$ for every
$\alpha\in\mathcal A_\lambda^\star$.
Consequently, the effective covariances
are block diagonal relative to
$\operatorname{ran}(P_\lambda^\star)$
and its orthogonal complement.
On $\mathcal A_\lambda^\star$, the
rank-$r$ eigengap equals the difference
between the smallest eigenvalue of
the leading block and the largest
eigenvalue of the complementary block.
Since both blocks depend affinely on
$\alpha$, this difference is a
nonnegative concave function of
$\alpha$.
It follows that
\[
\operatorname{gap}_r\{M_\lambda(\alpha)\}>0,
\qquad
\alpha\in
\operatorname{ri}(\mathcal A_\lambda^\star),
\]
where $\operatorname{ri}$ denotes
relative interior. Indeed, any
relative-interior point distinct
from a positive-gap dual optimizer
can be expressed as a convex
combination of that optimizer and
another dual optimizer, with both
coefficients strictly positive.
The conclusion follows from
concavity and nonnegativity of the
eigengap on $\mathcal A_\lambda^\star$.

This relative-interior property does
not, however, suffice for the
first-order perturbation results below.
Empirical dual optimizers may approach
boundary points at which the eigengap
vanishes, and the dual-set
characterization and quadratic-growth
argument used in the sharp-rate analysis
rely on uniform spectral separation.
A relaxation would therefore require
additional conditions controlling the
localization and stability of the
relevant empirical dual optimizers.
We retain
Assumption~\ref{ass:fair-uniform-optimal-gap}
to obtain perturbation guarantees
without imposing such selection
conditions.
A concrete example in Section~\ref{supp:subsec:fair-eigengap-separation} of the Supplementary Material illustrates that Assumption~\ref{ass:fair-optimal-eigengap} may hold even when the eigengap vanishes at a boundary point of the dual-optimal set, thereby violating Assumption~\ref{ass:fair-uniform-optimal-gap}.
\end{remark}
Define the groups attaining maximal regret
at this projector by
\[
\mathcal I_\lambda^\star
=
\{a:
\mathcal E_a(P_\lambda^\star)
=
B(P_\lambda^\star)\}.
\]
For $\alpha\in\Delta_G$, write
\[
T_\lambda(\alpha)
=
(I-P_\lambda^\star)
M_\lambda(\alpha)P_\lambda^\star.
\]
Thus $T_\lambda(\alpha)$ measures the
cross-subspace interaction of the effective
covariance relative to the population
Fair PCA principal subspace.

\begin{theorem}[Dual-set geometry and quadratic growth]
\label{thm:fair-derived-dual-growth}
Suppose $\lambda>0$ and
Assumption~\ref{ass:fair-uniform-optimal-gap}
holds. Then the population dual-optimal
set admits the characterization
\begin{equation}
\mathcal A_\lambda^\star
=
\left\{
\alpha\in\Delta_G:
\alpha_a=0\ \text{ for all }a\notin\mathcal I_\lambda^\star, T_\lambda(\alpha)=0
\right\}.
\label{eq:fair-polyhedral-dual-set}
\end{equation}
In particular, $\mathcal A_\lambda^\star$
is a nonempty compact convex polyhedron.
Moreover, there exists
$\underline\mu_\lambda>0$ such that
\begin{equation}
\phi_\lambda^\star-\phi_\lambda(\alpha)
\geq
\frac{\underline\mu_\lambda}{2}
\operatorname{dist}_2^2
(\alpha,\mathcal A_\lambda^\star),
\qquad
\alpha\in\Delta_G,
\label{eq:fair-derived-dual-growth}
\end{equation}
where
\(
\phi_\lambda^\star
=
\max_{\alpha\in\Delta_G}\phi_\lambda(\alpha)
\)
and $\operatorname{dist}_2
(\alpha,\mathcal A_\lambda^\star)
=
\inf_{\beta\in\mathcal A_\lambda^\star}
\|\alpha-\beta\|_2$.

The constant $\underline\mu_\lambda$
can be chosen explicitly in terms of
the group covariance matrices and a
Hoffman error-bound constant for
\eqref{eq:fair-polyhedral-dual-set}.
\end{theorem}

The result follows from an exact decomposition
of the population dual suboptimality:
\begin{equation}
\begin{aligned}
\phi_\lambda^\star-\phi_\lambda(\alpha)
=
\lambda\sum_{a=1}^G
\alpha_a
\{B(P_\lambda^\star)
-\mathcal E_a(P_\lambda^\star)\}
+
\mathcal K_r\{M_\lambda(\alpha)\}
-
\operatorname{tr}
\{P_\lambda^\star M_\lambda(\alpha)\}.
\end{aligned}
\label{eq:fair-dual-gap-decomposition}
\end{equation}
Both terms on the right-hand side are
nonnegative. The first measures weight
assigned to groups that are not active
at the population optimum. The second
measures the spectral optimality loss
of using $P_\lambda^\star$ for the
effective covariance $M_\lambda(\alpha)$.

The latter loss controls the squared
cross-block norm:
\begin{equation}
\mathcal K_r\{M_\lambda(\alpha)\}
-
\operatorname{tr}
\{P_\lambda^\star M_\lambda(\alpha)\}
\geq
\frac{
\|T_\lambda(\alpha)\|_F^2
}{
4(1+\lambda)\sigma
},
\label{eq:fair-spectral-crossblock-bound}
\end{equation}
where
\(
\sigma
=
\max_{1\leq a\leq G}
\|\Sigma_a\|_{\mathrm{op}},
\)
as defined at the beginning of this section.
Together with the polyhedral characterization
\eqref{eq:fair-polyhedral-dual-set},
this gives the quadratic-growth inequality
through a Hoffman error bound. The detailed
argument is provided in the Supplementary
Material.

For completeness, one admissible
quadratic-growth constant can be specified
as follows. Let $\kappa_\lambda>0$ satisfy
\begin{equation}
\begin{aligned}
\operatorname{dist}_2
(\alpha,\mathcal A_\lambda^\star)
\leq
\kappa_\lambda
\Bigg\{
\sum_{a\notin\mathcal I_\lambda^\star}
\alpha_a
+
\|T_\lambda(\alpha)\|_F
\Bigg\},
\end{aligned}
\label{eq:fair-hoffman-bound}
\end{equation}
for every $\alpha\in\Delta_G$.
Such a finite constant exists by
the polyhedral structure above.
If inactive groups exist, put
\[
d_\lambda
=
\min_{a\notin\mathcal I_\lambda^\star}
\{B(P_\lambda^\star)
-\mathcal E_a(P_\lambda^\star)\}>0.
\]
Then an admissible choice is
\begin{equation}
\underline\mu_\lambda
=
\frac{1}{\kappa_\lambda^2}
\begin{cases}
\displaystyle
\min\left\{
\lambda d_\lambda,
\frac{1}{4(1+\lambda)\sigma}
\right\},
&
\mathcal I_\lambda^\star
\neq\{1,\ldots,G\},
\\[3mm]
\displaystyle
\frac{1}{4(1+\lambda)\sigma},
&
\mathcal I_\lambda^\star
=\{1,\ldots,G\}.
\end{cases}
\label{eq:fair-derived-growth-modulus}
\end{equation}

The quadratic-growth condition is therefore
a consequence of the spectral and polyhedral
structure of Fair PCA.
Nevertheless, the derived constant may
deteriorate along a sequence of
high-dimensional models. The rate theorem
below explicitly accounts for this
possibility.

\subsubsection{First-order estimation of the principal subspace}

The population dual-optimal set can contain
multiple weight vectors, but
Assumption~\ref{ass:fair-uniform-optimal-gap}
ensures that they identify a common
principal subspace. The empirical dual
optimizer can therefore be compared with
its closest point in
$\mathcal A_\lambda^\star$.

The global quadratic-growth bound
\eqref{eq:fair-derived-dual-growth},
combined with uniform dual-objective
convergence, first provides quantitative
localization of the empirical weights.
Within the resulting eigengap neighborhood,
the population and empirical dual objectives
are differentiable. Controlling their
gradient difference then yields a
first-order dual-weight bound, from which
the projector rate follows by spectral
perturbation.
Define
\begin{equation}
s_{\lambda,n}
=
\frac{\lambda r\sqrt G}
{\underline\mu_\lambda}
\left\{
1+
\frac{(1+\lambda)\sigma}
{\delta_\lambda}
\right\}
\eta_n
\label{eq:fair-sharp-dual-scale}
\end{equation}
and
\begin{equation}
b_{\lambda,n}
=
(1+\lambda)\eta_n
+
\lambda\sigma\sqrt G\,s_{\lambda,n}.
\label{eq:fair-sharp-effective-scale}
\end{equation}
These quantities measure, respectively,
the dual-weight estimation scale and
the resulting effective-covariance
perturbation scale.
For the deterministic perturbation analysis, define
\begin{equation}
\mathfrak c_\lambda
=
(1+\lambda)
+
\frac{\lambda^2 rG\sigma}
{\underline\mu_\lambda}
\left\{
1+
\frac{(1+\lambda)\sigma}{\delta_\lambda}
\right\}.
\label{eq:fair-sharp-amplification}
\end{equation}
Thus $b_{\lambda,n}=\mathfrak c_\lambda\eta_n$.

\begin{theorem}[First-order projector perturbation bound]
\label{thm:fair-sharp-projector-rate}
Suppose $\lambda>0$ and
Assumption~\ref{ass:fair-uniform-optimal-gap}
holds. Let $\underline\mu_\lambda$ be the
derived quadratic-growth constant in
\eqref{eq:fair-derived-growth-modulus}.
There exist universal constants $c,C>0$
such that, whenever
\begin{equation}
\mathfrak c_\lambda e_n
\leq c\delta_\lambda,
\label{eq:fair-deterministic-sharp-condition}
\end{equation}
every empirical dual optimizer
$\widehat\alpha_\lambda$ satisfies
\begin{equation}
\begin{aligned}
\operatorname{dist}_2
(\widehat\alpha_\lambda,
\mathcal A_\lambda^\star)
\leq
C\frac{\lambda r\sqrt G}
{\underline\mu_\lambda}
\left\{
1+\frac{(1+\lambda)\sigma}
{\delta_\lambda}
\right\}e_n,
\end{aligned}
\label{eq:fair-deterministic-dual-bound}
\end{equation}
where \(e_n\) is defined in \eqref{eq:fair-covariance-error}.
Moreover, the empirical Fair PCA optimizer
is unique and is a rank-$r$ orthogonal
projector satisfying
\(
\widehat H_\lambda
=
P_r\{
\widehat M_\lambda(\widehat\alpha_\lambda)
\}.
\)
Furthermore,
\begin{equation}
\|\widehat H_\lambda-P_\lambda^\star\|_F
\leq
C\frac{\sqrt r}{\delta_\lambda}
\mathfrak c_\lambda e_n.
\label{eq:fair-deterministic-projector-bound}
\end{equation}
\end{theorem}

Theorem~\ref{thm:fair-sharp-projector-rate}
is a deterministic stability result.
It does not require a distributional
assumption or an asymptotic limit.
The small-perturbation condition \eqref{eq:fair-deterministic-sharp-condition} ensures
that the empirical dual optimizers remain
in a neighborhood where the relevant
spectral gaps are preserved.

Combining this result with simultaneous
covariance concentration yields the
following high-dimensional rate.

\begin{corollary}[High-dimensional projector rate]
\label{cor:fair-sharp-projector-rate}
Suppose $\lambda>0$ and
Assumptions~\ref{ass:fair-subgaussian}
and~\ref{ass:fair-uniform-optimal-gap}
hold. Suppose further that
\begin{equation}
\frac{b_{\lambda,n}}{\delta_\lambda}
\longrightarrow0.
\label{eq:fair-sharp-rate-condition}
\end{equation}
Then every empirical dual optimizer
$\widehat\alpha_\lambda$ satisfies
\begin{equation}
\operatorname{dist}_2
(\widehat\alpha_\lambda,
\mathcal A_\lambda^\star)
=
O_p(s_{\lambda,n}).
\label{eq:fair-sharp-dual-rate}
\end{equation}
With probability tending to one,
the empirical Fair PCA optimizer is unique
and is a rank-$r$ orthogonal projector, i.e.
\(
\widehat H_\lambda
=
P_r\{
\widehat M_\lambda(\widehat\alpha_\lambda)
\}.
\)
Moreover,
\begin{equation}
\|\widehat H_\lambda-P_\lambda^\star\|_F
=
O_p\left(
\frac{\sqrt r}{\delta_\lambda}
b_{\lambda,n}
\right).
\label{eq:fair-sharp-projector-rate}
\end{equation}
In particular, if
\begin{equation}
\frac{\sqrt r\,b_{\lambda,n}}
{\delta_\lambda}
\longrightarrow0,
\label{eq:fair-sharp-consistency-condition}
\end{equation}
then
\(
\|\widehat H_\lambda-P_\lambda^\star\|_F
\xrightarrow{p}0.
\)
\end{corollary}
Theorem~\ref{thm:fair-sharp-projector-rate}
establishes first-order stability of the empirical
Fair PCA projector under deterministic covariance
perturbations. Corollary
\ref{cor:fair-sharp-projector-rate}
translates this bound into a stochastic convergence
rate. The condition
\eqref{eq:fair-sharp-rate-condition}
ensures that the deterministic perturbation
requirement holds with probability tending to
one. Neither result requires uniqueness of the
population dual optimizer or a uniform positive
lower bound on $\underline\mu_\lambda$.

To interpret the rate, suppose that
$r$, $G$, $\lambda$, $\sigma$,
$\delta_\lambda^{-1}$, and
$\underline\mu_\lambda^{-1}$
remain bounded along the sequence of
population models. Then
\(
s_{\lambda,n}=O(\eta_n),
b_{\lambda,n}=O(\eta_n),
\)
and Corollary~\ref{cor:fair-sharp-projector-rate}
gives
\(
\|\widehat H_\lambda-P_\lambda^\star\|_F
=
O_p(\eta_n).
\)
Thus the first-order covariance-estimation
rate is recovered despite possible
nonuniqueness of the population dual weights.

More generally, the theorem allows
$p$, $r$, and $G$ to increase and the
spectral and identification constants
to vary with $n$. The resulting bound
makes explicit how covariance estimation,
spectral separation, and identification
of the dual-optimal set jointly determine
the accuracy of the estimated Fair PCA
subspace.

The requirement of dual identification
is specific to the sharper analysis:
the basic consistency result in
Theorem~\ref{thm:fair-basic-projector-rate}
continues to apply under the weaker
existence-type eigengap condition.
\section{Uniform Path Recovery and Penalty Selection}
\label{sec:fair-selection}
The fairness penalty controls the population
trade-off between aggregate reconstruction risk
and worst-group regret, and there is no
universally preferred penalty level. In
applications, its choice may instead reflect
a prespecified reconstruction--fairness
preference or an acceptable level of
worst-group regret.

We first establish uniform recovery of
the empirical reconstruction--fairness path
using the perturbation theory developed
in Section~\ref{sec:fair-sample-theory}.
We then derive statistical guarantees for
two approaches to penalty selection:
one based on a prescribed preference
over reconstruction and fairness, and
the other based on calibrating the penalty
to a specified fairness tolerance.
\subsection{Uniform recovery of the
reconstruction--fairness path}
\label{subsec:fair-selection-path}
Let
\[
\Lambda=[\lambda_-,\lambda_+]
\subset(0,\infty)
\]
be a compact interval of candidate penalties.
For each $\lambda\in\Lambda$, let
$P_\lambda^\star$ denote the population
Fair PCA projector, whenever uniquely
defined, and let $\widehat H_\lambda$
be an empirical optimizer of
\eqref{eq:fair-sample-estimator}.
To measure reconstruction loss relative
to ordinary PCA, define
\(
R_0
=
\min_{H\in\mathcal F_r}R(H),
\widehat R_0
=
\min_{H\in\mathcal F_r}\widehat R(H),
\)
and
\(
C(H)=R(H)-R_0,
\widehat C(H)=\widehat R(H)-\widehat R_0.
\)
Both quantities are nonnegative.
The subtraction of the ordinary-PCA
benchmark also removes the total-covariance
trace from their empirical difference.
Recall that $B$ and $\widehat B$ denote
the population and empirical worst-group
reconstruction regret criteria, defined
in Sections~\ref{sec:fair-pca}
and~\ref{subsec:fair-sample-setup},
respectively.
We define the population and empirical
performance paths by
\begin{equation}
\mathbf f(\lambda)
=
\begin{pmatrix}
C(P_\lambda^\star)\\
B(P_\lambda^\star)
\end{pmatrix},
\qquad
\widehat{\mathbf f}(\lambda)
=
\begin{pmatrix}
\widehat C(\widehat H_\lambda)\\
\widehat B(\widehat H_\lambda)
\end{pmatrix}.
\label{eq:fair-selection-performance-path}
\end{equation}
We derive a uniform performance bound
from the deterministic projector
perturbation theorem.
Suppose Assumption~\ref{ass:fair-uniform-optimal-gap}
holds at every $\lambda\in\Lambda$.
Recall the eigengap constants
$\delta_\lambda$ and the perturbation
amplification factors
$\mathfrak c_\lambda$ from Section
\ref{subsec:fair-sharp-rate}.
Define
\begin{equation}
\delta_\Lambda
=
\inf_{\lambda\in\Lambda}\delta_\lambda,
\qquad
\mathfrak c_\Lambda
=
\sup_{\lambda\in\Lambda}
\mathfrak c_\lambda.
\label{eq:fair-selection-uniform-constants}
\end{equation}
These quantities may depend on $n$.
The following result requires
$\delta_\Lambda>0$ and
$\mathfrak c_\Lambda<\infty$.

\begin{theorem}[Uniform reconstruction--fairness
path recovery]
\label{thm:fair-selection-uniform-path}
Suppose Assumption \ref{ass:fair-uniform-optimal-gap}
holds for every $\lambda\in\Lambda$,
with
$0<\delta_\Lambda$ and
$\mathfrak c_\Lambda<\infty$.
There exist universal constants
$c,C>0$ such that, whenever
\begin{equation}
\mathfrak c_\Lambda e_n
\leq c\delta_\Lambda,
\label{eq:fair-selection-uniform-condition}
\end{equation}
the empirical Fair PCA optimizer is
unique and is a rank-$r$ orthogonal
projector for every $\lambda\in\Lambda$.
Moreover,
\begin{equation}
\sup_{\lambda\in\Lambda}
\|\widehat H_\lambda-P_\lambda^\star\|_F
\leq
C\frac{\sqrt r\,\mathfrak c_\Lambda e_n}
{\delta_\Lambda},
\label{eq:fair-selection-uniform-projector}
\end{equation}
and
\begin{equation}
\begin{aligned}
\Delta_n
:=
\sup_{\lambda\in\Lambda}
\|\widehat{\mathbf f}(\lambda)
-\mathbf f(\lambda)\|_\infty\leq
2re_n+
C\frac{r\sigma\mathfrak c_\Lambda e_n}
{\delta_\Lambda},
\end{aligned}
\label{eq:fair-selection-uniform-performance}
\end{equation}
where
$\sigma=\max_a\|\Sigma_a\|_{\mathrm{op}}$.
Under Assumption~\ref{ass:fair-subgaussian},
if
\begin{equation}
\frac{\mathfrak c_\Lambda\eta_n}
{\delta_\Lambda}
\longrightarrow0,
\label{eq:fair-selection-uniform-localization-rate}
\end{equation}
then
\begin{equation}
\Delta_n
=
O_p\left\{
r\eta_n
\left(
1+
\frac{\sigma\mathfrak c_\Lambda}
{\delta_\Lambda}
\right)
\right\}.
\label{eq:fair-selection-path-rate}
\end{equation}
In particular, the empirical performance
path is uniformly consistent whenever,
in addition,
\begin{equation}
r\eta_n
\left(
1+
\frac{\sigma\mathfrak c_\Lambda}
{\delta_\Lambda}
\right)
\longrightarrow0.
\label{eq:fair-selection-path-consistency}
\end{equation}
\end{theorem}

Theorem~\ref{thm:fair-selection-uniform-path}
extends the fixed-penalty recovery results
of Section~\ref{sec:fair-sample-theory}
to an entire penalty interval.
The bound is uniform over both the
penalty and the corresponding fitted
representation. In particular,
\eqref{eq:fair-selection-uniform-performance}
controls the empirical reconstruction
and fairness criteria simultaneously.

The rate depends explicitly on the
effective spectral eigengaps and the
strength of the dual quadratic-growth
bound through $\delta_\Lambda$ and
$\mathfrak c_\Lambda$.
No uniqueness of the population dual
optimizer is required.
The performance of the fitted empirical
projector under the population criteria
also satisfies the same order of
uniform approximation:
\begin{equation}
\begin{aligned}
\sup_{\lambda\in\Lambda}
\left\|
\begin{pmatrix}
C(\widehat H_\lambda)\\
B(\widehat H_\lambda)
\end{pmatrix}
-
\mathbf f(\lambda)
\right\|_\infty
\leq
C\frac{
r\sigma\mathfrak c_\Lambda e_n
}{
\delta_\Lambda
}
\end{aligned}
\label{eq:fair-selection-fitted-performance}
\end{equation}
on the event
\eqref{eq:fair-selection-uniform-condition}.
Thus uniform recovery concerns not only
the observable empirical path, but also
the population performance of the fitted
representations.
\subsection{Preference-based penalty selection}
\label{subsec:fair-selection-preference}

A practitioner may express a preference
over the two performance criteria through
a function
$\Phi:\mathbb R_+^2\to\mathbb R$.
We do not impose a universally preferred
trade-off: the choice of $\Phi$ reflects
the relative importance assigned to
reconstruction accuracy and fairness.
Let $\Lambda_n\subseteq\Lambda$ be a
nonempty finite set of candidate penalties,
and define its covering radius by
\[
h_n
=
\sup_{\lambda\in\Lambda}
\min_{\mu\in\Lambda_n}
|\lambda-\mu|.
\]
Define the population and empirical
preference criteria
\(
D(\lambda)
=
\Phi\{\mathbf f(\lambda)\},
\widehat D(\lambda)
=
\Phi\{\widehat{\mathbf f}(\lambda)\},
\)
and select
\begin{equation}
\widehat\lambda
\in
\arg\min_{\lambda\in\Lambda_n}
\widehat D(\lambda).
\label{eq:fair-selection-empirical-penalty}
\end{equation}
The finite candidate set ensures that
the empirical minimum is attained and
allows the regularization path to be
computed on a prescribed grid.

The following result shows that the
fitted representation selected by
\eqref{eq:fair-selection-empirical-penalty}
achieves nearly optimal population
preference performance over the
regularization path.

\begin{corollary}[Oracle performance of
adaptively selected Fair PCA]
\label{cor:fair-selection-oracle}
Under the assumptions of
Theorem~\ref{thm:fair-selection-uniform-path},
let $\Phi$ be Lipschitz continuous
with constant $L_\Phi$ under the
$\ell_\infty$ norm.
On the event
\eqref{eq:fair-selection-uniform-condition},
\begin{equation}
\begin{aligned}
\left|
\Phi\{C(\widehat H_{\widehat\lambda}),
B(\widehat H_{\widehat\lambda})\}
-
\min_{\lambda\in\Lambda}D(\lambda)
\right|\leq
L_\Phi
\left\{
\Delta_n+2re_n+
\frac{2r\sigma^2}{\delta_\Lambda}h_n
\right\}.
\end{aligned}
\label{eq:fair-selection-oracle-performance}
\end{equation}
Under Assumption~\ref{ass:fair-subgaussian}
and conditions
\eqref{eq:fair-selection-uniform-localization-rate}
and
\eqref{eq:fair-selection-path-consistency},
if, in addition,
\begin{equation}
\frac{r\sigma^2}{\delta_\Lambda}h_n
\longrightarrow0,
\label{eq:fair-selection-grid-rate}
\end{equation}
then
\begin{equation}
\begin{aligned}
&
\Phi\{C(\widehat H_{\widehat\lambda}),
B(\widehat H_{\widehat\lambda})\}-
\min_{\lambda\in\Lambda}D(\lambda)
\xrightarrow{p}0.
\end{aligned}
\label{eq:fair-selection-oracle-consistency}
\end{equation}
\end{corollary}

The oracle inequality separates
statistical estimation error from
grid approximation error. Uniform
recovery of the reconstruction--fairness
path controls the former, while the
latter can be made negligible by refining
the candidate grid. Consequently, the
fitted empirical representation attains
asymptotically optimal population
preference performance relative to the
regularization path, without requiring
consistency of the selected penalty.
For example, the Lipschitz preference
$\Phi(c,b)=\max\{c/s_C,b/s_B\}$,
with prespecified $s_C,s_B>0$,
balances the two normalized performance
losses.

Further results on oracle performance,
consistency of the selected penalty,
and recovery of the corresponding
principal subspace are given in
Section~\ref{supp:subsec:fair-selection-preference}
of the Supplement.

\subsection{Penalty selection from a prescribed
fairness tolerance}
\label{subsec:fair-selection-tolerance}

An alternative to specifying a preference
function is to prescribe an acceptable
level $\tau$ of worst-group reconstruction
regret. We show that such a tolerance can
be used to calibrate the penalty in our
regularized Fair PCA formulation.
Let
\[
B_{\min}
=
\min_{H\in\mathcal F_r}B(H).
\]
For $\tau>B_{\min}$, convex duality
gives an associated population penalty
\begin{equation}
\lambda_\tau^\star
\in
\arg\max_{\lambda\geq0}
\{V(\lambda)-\lambda\tau\},
\label{eq:fair-selection-population-tolerance}
\end{equation}
where $V$ is the regularized value
function defined in Section
\ref{subsec:fair-pca-duality}.
The maximum is attained, and there
exists a corresponding penalized
optimizer satisfying
\(
B(H_{\lambda_\tau^\star}^\star)
\leq\tau,
\lambda_\tau^\star
\{B(H_{\lambda_\tau^\star}^\star)-\tau\}
=0.
\)
Thus, when $\lambda_\tau^\star>0$,
the selected optimizer attains the
prescribed tolerance exactly.
If a population dual-optimal
effective covariance at
$\lambda_\tau^\star$ has a strict
rank-$r$ eigengap, Theorem
\ref{thm:fair-spectral-characterization}
ensures that the penalized optimizer
is unique and equals the Fair PCA
projector $P_{\lambda_\tau^\star}^\star$.
The same calibration can be performed
using the empirical regularized
value function
\[
\widehat V(\lambda)
=
\min_{H\in\mathcal F_r}
\{\widehat R(H)
+\lambda\widehat B(H)\}.
\]
We select
\begin{equation}
\widehat\lambda(\tau)
\in
\arg\max_{\lambda\geq0}
\{\widehat V(\lambda)-\lambda\tau\},
\label{eq:fair-selection-tolerance-penalty}
\end{equation}
and solve the empirical regularized
Fair PCA problem at this penalty.

\begin{proposition}[Fairness-tolerance
calibration]
\label{prop:fair-selection-tolerance-calibration}
Suppose $\tau>B_{\min}$ and
\[
2re_n<\tau-B_{\min}.
\]
Then the maximization problem
\eqref{eq:fair-selection-tolerance-penalty}
admits a finite solution.
For any optimal penalty
$\widehat\lambda(\tau)$,
an empirical penalized optimizer
$\widehat H_{\widehat\lambda(\tau)}$
can be selected such that
\begin{equation}
\begin{aligned}
\widehat B
(\widehat H_{\widehat\lambda(\tau)})
\leq\tau,\quad
\widehat\lambda(\tau)
\{\widehat B
(\widehat H_{\widehat\lambda(\tau)})
-\tau\}
=0.
\end{aligned}
\label{eq:fair-selection-tolerance-slackness}
\end{equation}
Consequently,
\begin{equation}
\left[
B(\widehat H_{\widehat\lambda(\tau)})
-\tau
\right]_+
\leq2re_n.
\label{eq:fair-selection-tolerance-violation}
\end{equation}
If $\widehat\lambda(\tau)>0$,
then
\begin{equation}
\left|
B(\widehat H_{\widehat\lambda(\tau)})
-\tau
\right|
\leq2re_n.
\label{eq:fair-selection-tolerance-accuracy}
\end{equation}
\end{proposition}

The proposition provides a direct
statistical guarantee for penalty
selection from a prescribed fairness
tolerance. Under
Assumption~\ref{ass:fair-subgaussian},
the population fairness violation is
$O_p(r\eta_n)$ whenever the stated
feasibility condition holds with
probability tending to one.
If the calibrated penalty is positive
with probability tending to one,
the same rate holds for the absolute
deviation from the prescribed tolerance.

When the empirical regularized optimizer
at $\widehat\lambda(\tau)$ is unique,
it necessarily satisfies
\eqref{eq:fair-selection-tolerance-slackness}.
If the optimizer is nonunique,
strong duality guarantees the existence
of a solution satisfying these conditions,
but an arbitrary regularized optimizer
need not do so. The calibration guarantee applies to
a Fantope optimizer and does not,
without an additional eigengap
condition, ensure that the selected
solution is a rank-$r$ projector.

In practice, the calibrated penalty
may be determined by a one-dimensional
search over $\lambda$, with each evaluation
requiring a regularized Fair PCA fit.
The optimal worst-group regret is
nonincreasing in $\lambda$, allowing
bisection-based calibration when the
regret path is continuous and the
target tolerance lies between the
regrets at two trial penalties.
Alternatively, the fairness-constrained
problem can be solved directly as a
semidefinite program.
\subsection{Computational implementation}
\label{subsec:fair-selection-computation}

We summarize the computational procedure
for obtaining the empirical Fair PCA
representation. The same regularized
fitting procedure applies whether the
penalty is prescribed or selected
adaptively.

\begin{enumerate}
\item \textit{Empirical covariance estimation.}
Given the group observations, target rank
$r$, and prespecified group weights
$q_1,\ldots,q_G$, compute the group
covariance matrices $\widehat\Sigma_a$,
the average covariance
$\widehat{\overline\Sigma}$, and
$\widehat c_a=\mathcal K_r(\widehat\Sigma_a)$,
as defined in Section~
\ref{subsec:fair-sample-setup}.

\item \textit{Regularized Fair PCA fitting.}
For any given penalty $\lambda>0$,
solve the empirical dual problem
\(
\widehat\alpha_\lambda
\in
\arg\max_{\alpha\in\Delta_G}
\widehat\phi_\lambda(\alpha),
\)
where $\widehat\phi_\lambda$ is defined
in Section~\ref{subsec:fair-sample-setup}.
If the effective covariance
$\widehat M_\lambda(\widehat\alpha_\lambda)$
has a strict rank-$r$ eigengap, the
empirical Fair PCA estimator is
\(
\widehat H_\lambda
=
P_r\{
\widehat M_\lambda(\widehat\alpha_\lambda)
\}.
\)
Its range gives the estimated principal
subspace. In the absence of the required
eigengap, a primal optimizer of
\eqref{eq:fair-sample-estimator} can
instead be computed directly over
$\mathcal F_r$. For $\lambda=0$,
the procedure reduces to ordinary PCA.

\item \textit{Penalty selection and final fit.}
For a practitioner-specified penalty,
the preceding fitting step is performed
directly. For preference-based selection,
it is applied to each candidate
$\lambda\in\Lambda_n$, and the final
penalty is selected according to
\eqref{eq:fair-selection-empirical-penalty}.
Alternatively, a prescribed fairness
tolerance can be used to determine
$\widehat\lambda(\tau)$ through
\eqref{eq:fair-selection-tolerance-penalty},
with regularized fits evaluated as
needed during calibration. In the latter
case, a penalized optimizer satisfying
\eqref{eq:fair-selection-tolerance-slackness}
is selected. The resulting fitted
representation is then returned.
\end{enumerate}

The empirical dual objective is concave
over the probability simplex, so standard
convex optimization methods can be used
to compute the dual weights. Numerical
optimality can be assessed through the
primal--dual gap, while the rank-$r$
eigengap determines whether the leading
spectral projector is uniquely defined.
\section{Simulation Studies}
\label{sec:fair-simulations}
We investigate the finite-sample
reconstruction--fairness trade-off under
heterogeneous group covariances, with
particular attention to the effects of
dimension and group sample sizes.
Additional numerical illustrations in
Section~\ref{supp:sec:fair-theory-numerics}
of the Supplementary Material examine
the local reconstruction--fairness expansions
and the first-order statistical recovery
of the Fair PCA subspace, including
a setting with nonunique population
dual-optimal weights.
Further simulations in
Section~\ref{supp:sec:numerical}
compare finite-penalty Fair PCA with
ordinary sample PCA and empirical
minimax-regret PCA, highlighting the
potential benefits of moderate
regularization.

We consider a heterogeneous spiked
covariance model with group-specific
signal subspaces.
Consider $G=3$ groups with
\(
\Sigma_a
=
\nu_a I_p+U_aD_aU_a^\top,
a=1,2,3,
\)
where
\(
(\nu_1,\nu_2,\nu_3)=(1,1.3,2),
\)
and
\(
D_1=\operatorname{diag}(8,5),
D_2=\operatorname{diag}(9,4.5),
D_3=\operatorname{diag}(10,4).
\)
Writing $e_j$ for the $j$th canonical basis
vector, the group-specific signal directions
are given by
$U_1
=
(e_1,e_2),$ $U_2
=
(\cos(0.68)e_1+\sin(0.68)e_3,\ e_2),$ and $U_3
=
(\cos(1.02)e_1+\sin(1.02)e_4,\
 \cos(0.45)e_2+\sin(0.45)e_5).$
We set $r=2$ and prescribe the group weights
$(q_1,q_2,q_3)=(0.60,0.25,0.15)$.
Independent observations are generated from
$N_p(0,\Sigma_a)$ within each group.
We consider $p\in\{30,100,200\}$ and
two configurations of group sample sizes,
\[
(n_1,n_2,n_3)
\in
\{(300,150,100),(900,450,300)\}.
\]
For each configuration, 1,000 Monte Carlo
replications are performed.
The empirical Fair PCA estimator
$\widehat H_\lambda$ is computed using
the dual spectral formulation of
Section~\ref{sec:fair-sample-theory}.
We report results for
$\lambda\in\{0,0.2,0.5,1,2\}$.
All reconstruction risks and regrets are
evaluated using the population covariance
matrices.

\begin{figure}[t]
\centering
\includegraphics[width=\textwidth]
{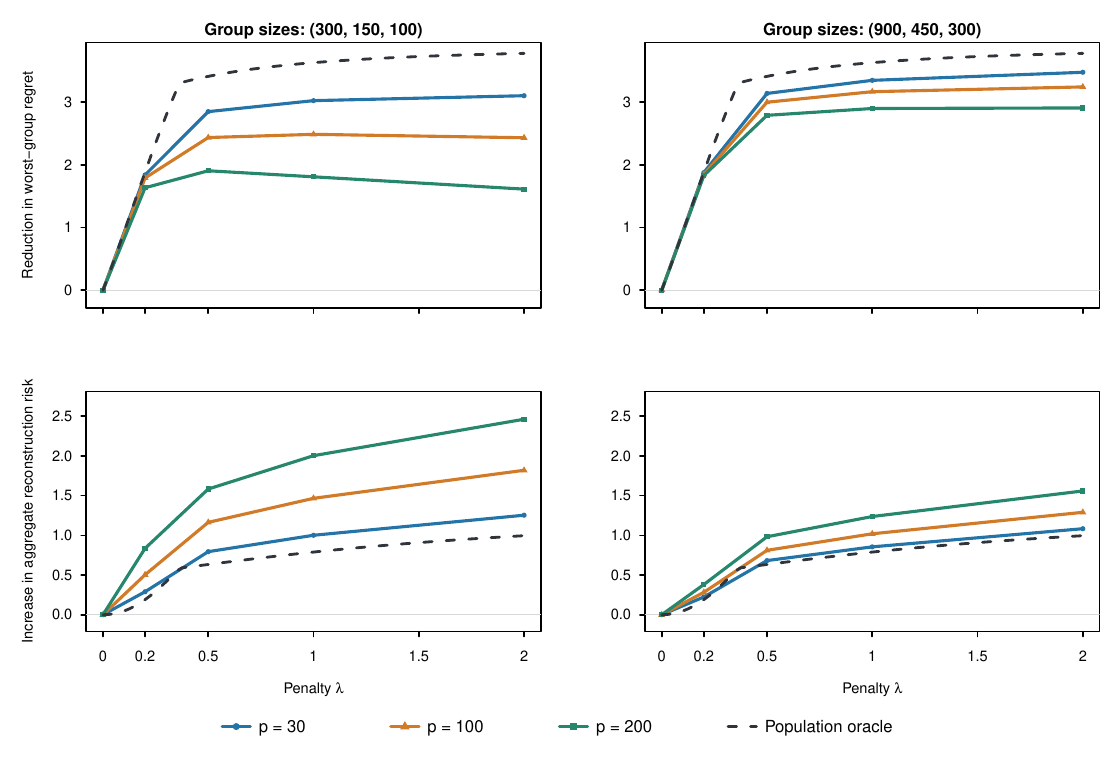}
\caption{
Reconstruction--fairness trade-off under
heterogeneous group covariances.
The upper panels show the Monte Carlo mean
reduction in worst-group regret,
$B(\widehat H_0)-B(\widehat H_\lambda)$,
and the lower panels show the corresponding
increase in aggregate reconstruction risk,
$R(\widehat H_\lambda)-R(\widehat H_0)$.
Both are evaluated using the population
covariance matrices.
Solid curves represent empirical Fair PCA,
while dashed curves show the population-oracle
differences relative to $P_0^\star$.
Results are based on 1,000 Monte Carlo
replications per configuration.
}
\label{fig:fair-sim-tradeoff}
\end{figure}

Figure~\ref{fig:fair-sim-tradeoff}
displays the reduction in worst-group
reconstruction regret and the corresponding
change in aggregate reconstruction risk
relative to ordinary PCA.
The population-oracle trade-off is unchanged
across dimensions in this design, since the
additional coordinates contribute only
isotropic noise. The empirical curves
therefore illustrate the effect of dimension
and group sample sizes on the attainable
finite-sample trade-off.
For reference,
Table~\ref{tab:fair-sim-ordinary-pca}
reports the population performance of
ordinary sample PCA, corresponding to
$\lambda=0$.
\begin{table}[t]
\centering
\caption{
Ordinary sample PCA performance at
$\lambda=0$ and population PCA-captured
variance. Reconstruction risks and regrets
are Monte Carlo means evaluated using the
population covariance matrices over
1,000 replications.
}
\label{tab:fair-sim-ordinary-pca}
\small
\begin{tabular}{cc cc cc}
\hline
&
&
\multicolumn{2}{c}{$(300,150,100)$}
&
\multicolumn{2}{c}{$(900,450,300)$}
\\
\cline{3-4}\cline{5-6}
$p$
&
$\mathcal K_2(\overline\Sigma)$
&
$R(\widehat H_0)$
&
$B(\widehat H_0)$
&
$R(\widehat H_0)$
&
$B(\widehat H_0)$
\\
\hline
30  & 13.926 & 36.260  & 6.977 & 36.152  & 6.931 \\
100 & 13.926 & 122.364 & 7.188 & 122.025 & 7.013 \\
200 & 13.926 & 245.346 & 7.479 & 244.697 & 7.102 \\
\hline
\end{tabular}
\end{table}
The quantity $\mathcal K_2(\overline\Sigma)$
is the maximum aggregate variance captured
by a rank-two projector. Since
\(
R(P_0^\star)
=
\operatorname{tr}(\overline\Sigma)
-
\mathcal K_2(\overline\Sigma),
\)
it also determines the minimum population
reconstruction risk.
In this design,
$\mathcal K_2(\overline\Sigma)$ remains
constant across dimensions, whereas the
baseline reconstruction risk increases
with $p$ because of the additional
isotropic noise coordinates.
Moderate fairness penalties yield substantial
reductions in worst-group regret at relatively
modest reconstruction costs.
For example, at $\lambda=0.5$ and $p=30$,
under the smaller-sample configuration,
mean worst-group regret decreases from
$6.977$ to $4.127$, a reduction of
$40.9\%$, while aggregate reconstruction
risk increases by only $0.794$. The reconstruction cost of $0.794$
corresponds to approximately $5.7\%$
of the variance captured by population
ordinary PCA.
At $p=200$, the corresponding regret
reduction is $25.5\%$, with a reconstruction
cost of $1.583$.
Increasing the group sample sizes to
$(900,450,300)$ improves the trade-off
at $p=200$, yielding a $39.3\%$ reduction
in worst-group regret at a reconstruction
cost of $0.982$.

These results illustrate the influence of
covariance-estimation accuracy on the
empirical reconstruction--fairness trade-off.
In this design, increasing $p$ adds
isotropic-noise coordinates without changing
the underlying signal structure or the
population trade-off. The resulting increase
in the effective rank of the group covariance
matrices makes their estimation more difficult
at fixed sample sizes, whereas larger samples
bring the empirical trade-off closer to its
population counterpart.
Moreover, the population-evaluated regret
$B(\widehat H_\lambda)$ need not be
nonincreasing in $\lambda$, in contrast
to the population-oracle regret
$B(P_\lambda^\star)$.
For $p=200$ under the smaller-sample
configuration, the Monte Carlo mean of
$B(\widehat H_\lambda)$ decreases up to
$\lambda=0.5$ but increases slightly
at larger penalties.

\section{Application to NHANES accelerometry}
\label{sec:fair-nhanes}
Daily accelerometry profiles provide a natural
setting for studying dimension reduction across
heterogeneous groups. A common low-dimensional
representation facilitates comparisons across
participants, but differences in within-group
covariance structures may cause ordinary PCA
to preserve the dominant modes of variation
of some groups better than others. We investigate
this issue using accelerometry data from the
2011--2012 cycle of the National Health and
Nutrition Examination Survey (NHANES), examining
whether Fair PCA can reduce the largest
group-specific reconstruction regret while
maintaining aggregate reconstruction accuracy.

\subsection{Data and preprocessing}
\label{subsec:fair-nhanes-data}

The NHANES physical activity monitor component collected
wrist-worn triaxial accelerometer measurements over seven
consecutive monitoring days. We use the publicly available
minute-level summaries, expressed in Monitor-Independent
Movement Summary (MIMS) units, together with the monitoring-day
and demographic records.\footnote{Data and documentation
are available at
\url{https://wwwn.cdc.gov/Nchs/Data/Nhanes/Public/2011/DataFiles/PAXMIN_G.htm}.}

For each participant, we construct a 24-hour activity
profile by dividing the day into 288 consecutive
five-minute intervals. Each coordinate is the mean
of all valid minute-level MIMS measurements within
the corresponding clock-time interval, pooled across
eligible monitoring days. Thus, each participant
contributes a single vector in $\mathbb R^{288}$,
rather than separate observations for individual
monitoring days.

We retain minutes with nonnegative triaxial MIMS
measurements (\texttt{PAXMTSM}), no quality-control
flags (\texttt{PAXQFM}=0), at least 55 seconds
of recorded data (\texttt{PAXTSM}), and predicted
wake or sleep wear status (\texttt{PAXPREDM}).
A monitoring day is eligible if it contains
at least 960 valid minutes. We retain participants
with at least three eligible days and at least
50\% valid-minute coverage in every five-minute
clock-time interval across their eligible days.
No missing activity values are imputed.

These criteria yield 6,055 participants, partitioned into four
age groups: 6--17 years ($n_1=1739$), 18--39 years
($n_2=1537$), 40--64 years ($n_3=1787$), and 65 years or older
($n_4=992$). We adopt equal reconstruction weights,
$q_a=1/4$, so that each age group contributes equally to the
aggregate reconstruction criterion. For covariance estimation, activity profiles
are centered separately within each age group.
The weights $q_a$ specify the aggregate
reconstruction criterion and are distinct
from NHANES survey sampling weights.
We do not incorporate the complex survey
design; the analysis therefore evaluates
the reconstruction--fairness trade-off
within the retained cohort rather than
providing nationally representative inference.

\begin{figure}[t]
\centering
\includegraphics[width=\linewidth]
{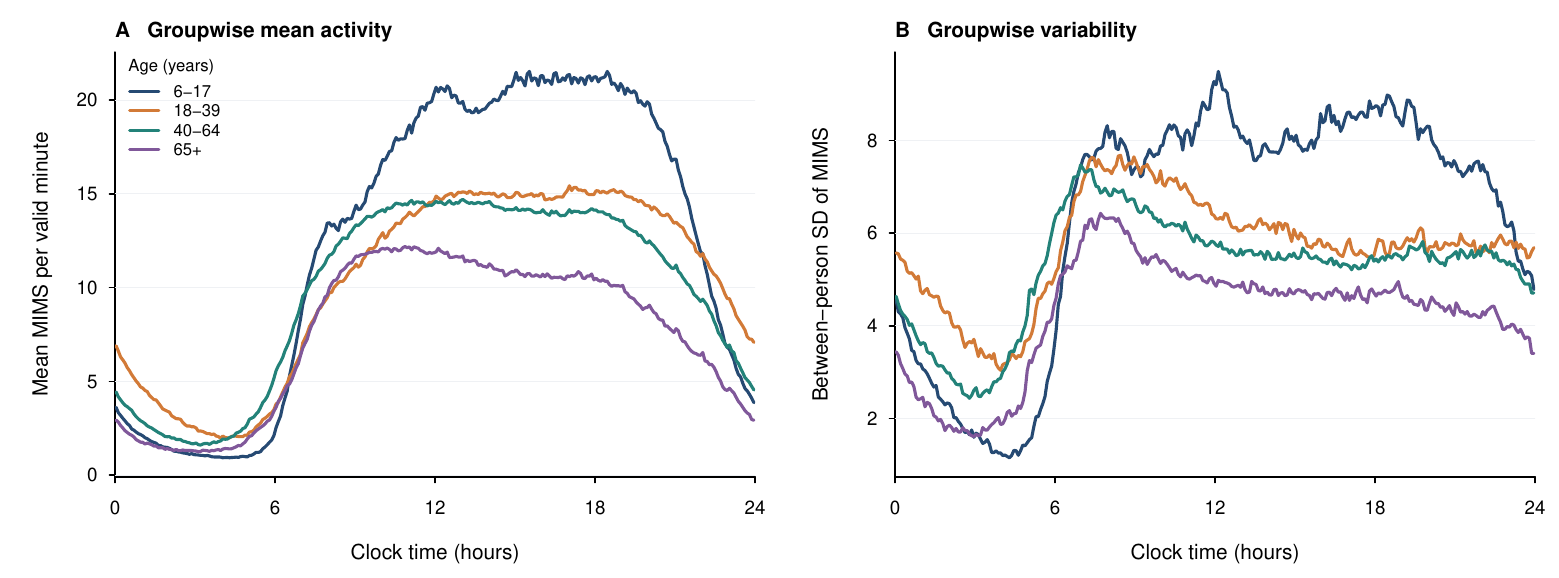}
\caption{
NHANES daily activity profiles by age group.
Panel A shows the mean MIMS activity per valid minute,
and Panel B shows the corresponding between-participant
standard deviation. Both are evaluated at five-minute
clock-time intervals using the 6,055 retained participants.
}
\label{fig:fair-nhanes-activity}
\end{figure}

Figure~\ref{fig:fair-nhanes-activity}
displays the groupwise mean activity profiles
and their pointwise between-participant
standard deviations, computed from the
original profiles before within-group centering. Activity is generally
concentrated in daytime hours, with differences
across age groups. Participants aged 6--17
exhibit higher daytime mean activity and
greater between-participant variability over
substantial portions of the day. Although the mean profiles
provide descriptive context, Fair PCA is determined by the
within-group covariance structures rather than differences
between group means.

\subsection{Experimental design and evaluation}
\label{subsec:fair-nhanes-design}

We consider rank $r=2$, for which ordinary PCA explains
approximately 44.9\% of the weighted sample variance,
and examine the penalty grid
\(
\lambda\in
\{0,\,0.025,\,0.05,\,0.1,\,0.2,\,0.4\}.
\)
Participants are randomly divided into training
and test samples in an 80:20 ratio, with
stratification by age group and sex to maintain
comparable sample compositions.
This procedure is repeated for 30 random
splits to assess the stability of the
empirical reconstruction--fairness trade-off.
Sex is used only for stratification; the
fairness criterion is defined over the
four age groups.
For each split, we compute the empirical Fair PCA estimator
$\widehat H_\lambda$ from the training sample, as defined
in \eqref{eq:fair-sample-estimator}. The same fitted projector
is then evaluated on both training and test observations.
In particular, the test sample is not used to construct
$\widehat H_\lambda$.

For clarity, let $D\in\{\mathrm{tr},\mathrm{te}\}$ denote
the training or test partition, and let
$\widehat\Sigma_{D,a}$ be the corresponding groupwise-centered
sample covariance matrix. The evaluation criteria are
\(
\widehat R_D(H)
=
\sum_{a=1}^4
q_a\operatorname{tr}
\{(I-H)\widehat\Sigma_{D,a}\},
\)
and
\(
\widehat B_D(H)
=
\max_{1\leq a\leq4}
\left\{
\mathcal K_r(\widehat\Sigma_{D,a})
-
\operatorname{tr}(H\widehat\Sigma_{D,a})
\right\}.
\)
Thus, group-specific PCA oracle risks are evaluated
separately within each partition.

We report performance relative to the ordinary PCA
projector $\widehat H_0$ fitted on the same training sample.
Specifically, the percentage reduction in worst-group regret
and the percentage reconstruction cost on partition $D$ are
\(
100
\frac{
\widehat B_D(\widehat H_0)
-
\widehat B_D(\widehat H_\lambda)
}{
\widehat B_D(\widehat H_0)
}
\)
and
\(
100
\frac{
\widehat R_D(\widehat H_\lambda)
-
\widehat R_D(\widehat H_0)
}{
\widehat R_D(\widehat H_0)
},
\)
respectively. Positive values indicate an improvement
in worst-group regret for the first measure and a
deterioration in reconstruction accuracy for the second.
Importantly, test reconstruction cost is measured relative
to the training-fitted ordinary PCA projector, not an
ordinary PCA projector refitted to the test sample.

\subsection{Empirical reconstruction--fairness trade-off}
\label{subsec:fair-nhanes-results}

\begin{figure}[t]
\centering
\includegraphics[width=\linewidth]
{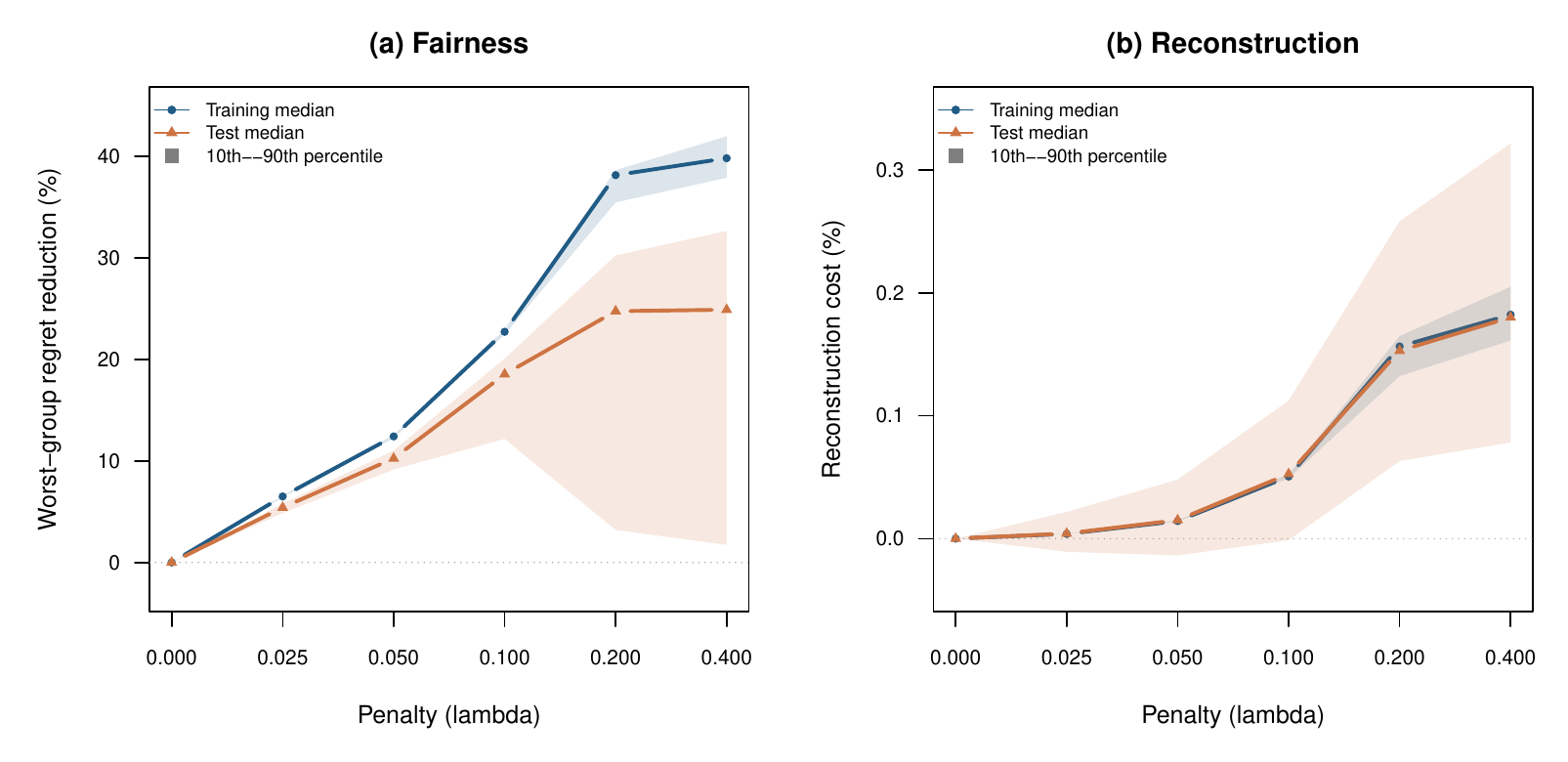}
\caption{
Reconstruction--fairness paths for NHANES
accelerometry ($p=288$, $r=2$).
Panel A shows the percentage reduction in
worst-group regret, and Panel B the percentage
reconstruction cost, both relative to ordinary
PCA. Curves are pointwise medians across
30 train/test splits, with shaded 10th--90th
percentile bands. All projectors are fitted
on training data; the bands describe variation
across overlapping partitions, not confidence
intervals.
}
\label{fig:fair-nhanes-tradeoff}
\end{figure}

Figure~\ref{fig:fair-nhanes-tradeoff} displays the
training and test reconstruction--fairness paths
over 30 random splits. At $\lambda=0.1$,
the median worst-group regret reduction is
22.72\% on the training samples and 18.54\%
on the test samples, with a median test
reconstruction cost of only 0.053\%.
Test regret improves in 29 of the 30 splits. Under the alternative normalization by
the maximum rank-two variance captured
on the test sample, the median reconstruction
cost remains small at $0.065\%$; see Section \ref{supp:subsec:fair-nhanes-captured} of the
Supplementary Material.
Increasing $\lambda$ to $0.2$ yields a larger
median test regret reduction of 24.75\%, but
also increases the reconstruction cost to
0.153\%. At $\lambda=0.4$, the corresponding
figures are 24.89\% and 0.180\%.
Thus, further penalization produces diminishing
fairness gains at increasing reconstruction cost.
The wider test bands at larger penalties also
indicate greater sensitivity to the random
training--test partition.

Table~\ref{tab:fair-nhanes-group-regret}
shows how the improvement is distributed
across age groups. Under ordinary PCA,
participants aged 6--17 have the largest
test regret in 29 of the 30 splits.
At $\lambda=0.1$, Fair PCA reduces this
group's regret in every split, while the
other groups incur additional regret.
This redistribution reduces the worst-group
regret in most splits, without requiring
improvements for every group.

\begin{table}[t]
\centering
\caption{
Group-specific test reconstruction regrets
under ordinary PCA and Fair PCA at
$\lambda=0.1$. Entries are medians across
30 splits; changes are computed within
each split before taking medians.
}
\label{tab:fair-nhanes-group-regret}
\begin{tabular}{lrrr}
\toprule
Age group
& Ordinary PCA
& Fair PCA
& Change \\
\midrule
6--17  & 291.30 & 237.10 & $-53.78$ \\
18--39 & 150.51 & 179.03 & $+23.59$ \\
40--64 &  84.68 & 108.10 & $+21.60$ \\
65+    & 132.32 & 149.94 & $+18.01$ \\
\bottomrule
\end{tabular}
\end{table}

Overall, Fair PCA achieves substantial reductions
in worst-group reconstruction regret at small
aggregate reconstruction costs. The differences
between training and test performance, including
one unfavorable test split at $\lambda=0.1$,
highlight the role of finite-sample variability.
The penalty $\lambda=0.1$ serves as an
illustrative compromise rather than an
independently validated choice. The full penalty paths in
Figure~\ref{fig:fair-nhanes-tradeoff} provide the
primary assessment of the reconstruction--fairness
relationship.

\subsection{Sensitivity and interpretation}
\label{subsec:fair-nhanes-sensitivity}

Additional analyses examine sensitivity to
the projection rank, reconstruction weights,
and temporal resolution. The favorable
reconstruction--fairness trade-off persists
across the examined ranks and temporal
resolutions, while its magnitude depends
on the reconstruction weights. We also
consider an alternative normalization of
reconstruction cost and the minimax-regret
endpoint $\lambda=\infty$. The latter
achieves the greatest training fairness
improvement but does not yield the largest
median test improvement in this application.
Detailed results are reported in
Section~\ref{supp:sec:fair-nhanes}
of the Supplementary Material. 
\section{Discussion}
\label{sec:fair-discussion}

Fairness-regularized PCA provides a way to
control group-specific reconstruction regret
while retaining a spectral interpretation
through fairness-adjusted covariance weights.
Our population analysis identifies conditions
under which small improvements in worst-group
regret incur only second-order reconstruction
costs, as well as conditions under which the
regularization path is locally constant.
The statistical theory quantifies the recovery
of the resulting principal subspaces from
empirical covariance matrices, including when
the optimal fairness-adjusted weights are
nonunique.

The simulations and NHANES application illustrate
the potential for meaningful reductions in
worst-group regret with limited reconstruction
loss. Nevertheless, the favorable local
trade-off does not imply negligible costs
at arbitrary penalty levels, and the
appropriate balance between reconstruction
accuracy and fairness remains application
dependent. Moreover, the fairness criterion
concerns reconstruction regret rather than
broader notions of fairness in downstream
statistical decisions.

\begin{supplement}
\stitle{Supplementary material for
Fairness-Regularized Principal Component Analysis}
\sdescription{
Contains proofs and auxiliary theoretical results,
additional finite-sample simulations, and
supplementary analyses of NHANES accelerometry.
}
\end{supplement}

\bibliographystyle{imsart-nameyear}
\bibliography{references}

\clearpage

\section*{Supplementary material}

\setcounter{theorem}{0}
\setcounter{lemma}{0}
\setcounter{proposition}{0}
\setcounter{corollary}{0}
\setcounter{table}{0}
\setcounter{figure}{0}
\setcounter{equation}{0}

\renewcommand{\thetheorem}{S\arabic{theorem}}
\renewcommand{\thelemma}{S\arabic{lemma}}
\renewcommand{\theproposition}{S\arabic{proposition}}
\renewcommand{\thecorollary}{S\arabic{corollary}}
\renewcommand{\thesection}{S.\arabic{section}}
\renewcommand{\thefigure}{S.\arabic{figure}}
\renewcommand{\thetable}{S.\arabic{table}}
\section{Overview and notation}
\label{supp:sec:overview}

This supplementary document accompanies
``Fairness-Regularized PCA.'' It contains
proofs and auxiliary results supporting the
theoretical developments in the main paper,
together with additional numerical investigations.

Section~\ref{supp:sec:spectral} establishes
the spectral characterization and exactness
of the Fantope relaxation.
Section~\ref{supp:sec:fair-tradeoff} provides
the global reconstruction--fairness properties,
the technical arguments underlying the local
expansion around ordinary PCA, and illustrative
examples of local degeneracy and nondegeneracy.
Section~\ref{supp:sec:fair-sample-theory}
develops the covariance concentration,
spectral perturbation, and dual-set geometry
results used to establish statistical
consistency and projector recovery.
Section~\ref{supp:sec:fair-selection} contains
the proofs for uniform path recovery and
data-adaptive penalty selection.
The remaining sections provide numerical
illustrations of the theoretical results
(Section~\ref{supp:sec:fair-theory-numerics}),
additional simulations comparing
finite-penalty and minimax-regret Fair PCA
(Section~\ref{supp:sec:numerical}),
and further analyses of NHANES accelerometry
(Section~\ref{supp:sec:fair-nhanes}).

We follow the notation of the main paper
throughout. For convenience, we recall the
principal objects used in the proofs.
There are $G$ groups with covariance matrices
$\Sigma_1,\ldots,\Sigma_G$ and prespecified
weights $q_1,\ldots,q_G$, where
$\sum_{a=1}^G q_a=1$.
For a target rank $1\leq r<p$, let
\[
\overline\Sigma
=
\sum_{a=1}^G q_a\Sigma_a,
\qquad
\mathcal K_r(A)
=
\sum_{j=1}^r\theta_j(A),
\]
where the eigenvalues $\theta_1(A)\geq\cdots\geq\theta_p(A)$ of a symmetric matrix
$A$ are arranged in nonincreasing order.
We write $P_r(A)$ for a projector onto
a leading rank-$r$ eigenspace of $A$.

Let $\mathcal P_r$ denote the set of
rank-$r$ orthogonal projectors and
$\mathcal F_r$ its convex hull, the
rank-$r$ Fantope. For $H\in\mathcal F_r$,
the aggregate reconstruction risk and
worst-group reconstruction regret are
\[
\begin{aligned}
R(H)
=
\operatorname{tr}
\{(I-H)\overline\Sigma\}, \qquad
B(H)
=
\max_{1\leq a\leq G}
\left\{
\mathcal K_r(\Sigma_a)
-
\operatorname{tr}(H\Sigma_a)
\right\}.
\end{aligned}
\]
The population regularized optimizer is
denoted by $H_\lambda^\star$.
When the relevant eigengap condition holds,
it is the unique rank-$r$ projector
$P_\lambda^\star$.

The minimax formulation uses the probability
simplex $\Delta_G$, the effective covariance
\[
M_\lambda(\alpha)
=
\overline\Sigma
+
\lambda\sum_{a=1}^G
\alpha_a\Sigma_a,
\qquad
\alpha\in\Delta_G,
\]
and the population dual-optimal set
$\mathcal A_\lambda^\star$.
Empirical quantities are indicated by hats;
in particular, $\widehat\Sigma_a$ and
$\widehat H_\lambda$ denote the sample
covariance matrices and empirical Fair PCA
optimizer, respectively.

We use $\|\cdot\|_2$ for the Euclidean norm,
and $\|\cdot\|_{\mathrm{op}}$, $\|\cdot\|_F$,
and $\|\cdot\|_*$ for the matrix operator,
Frobenius, and nuclear norms, respectively.
For symmetric matrices, $A\preceq B$ means
that $B-A$ is positive semidefinite.
The notation $O_p(\cdot)$ and $o_p(\cdot)$
denotes the usual stochastic orders.
For a symmetric matrix $A$, we write
\[
\operatorname{gap}_r(A)
=
\theta_r(A)-\theta_{r+1}(A)
\]
for its rank-$r$ eigengap. For
$\alpha\in\mathbb R^G$ and a nonempty
set $\mathcal A\subseteq\mathbb R^G$,
let
\[
\operatorname{dist}_2(\alpha,\mathcal A)
=
\inf_{\beta\in\mathcal A}
\|\alpha-\beta\|_2
\]
denote the Euclidean distance to
$\mathcal A$.

The statistical arguments draw on spectral
perturbation and covariance concentration
techniques. Relevant background is provided by
\citet{YuWangSamworth2015},
\citet{KoltchinskiiLounici2017Spectral},
and \citet{vershynin2026highdimensional}.
References for specific auxiliary inequalities
and error bounds are given at their points
of use.

All assumptions and regularity conditions
are those specified in the main paper or
explicitly stated in the corresponding
supplementary results. When parameters
such as $p$, $r$, $G$, and the covariance
matrices vary with sample size, this
dependence is suppressed unless needed
for clarity. Additional notation used
only in particular arguments is introduced
locally.

\section{Proofs for Section 2: Spectral Characterization}
\label{supp:sec:spectral}

We prove the spectral characterization theorem
by combining the minimax identity with the
Ky Fan variational principle. The optimal
effective eigengap ensures uniqueness of
the relaxed solution and exact recovery
of a rank-$r$ projector.
\subsection{Proof of the spectral characterization theorem}
\label{supp:subsec:spectral-proof}
\begin{proof}[Proof of Theorem~\ref{thm:fair-spectral-characterization}]
Let $H^\star$ be any minimizer of
\eqref{eq:fair-fantope-oracle}, and let
$\alpha^\star\in\mathcal A_\lambda^\star$ satisfy
\eqref{eq:fair-optimal-eigengap}. By minimax equality,
\[
\min_{H\in\mathcal F_r}
L_\lambda(H,\alpha^\star)
=
V(\lambda)
=
\max_{\alpha\in\Delta_G}
L_\lambda(H^\star,\alpha).
\]
Since
\[
\min_{H\in\mathcal F_r}
L_\lambda(H,\alpha^\star)
\leq
L_\lambda(H^\star,\alpha^\star)
\leq
\max_{\alpha\in\Delta_G}
L_\lambda(H^\star,\alpha),
\]
equality must hold throughout. Thus $H^\star$ minimizes
$L_\lambda(\,\cdot\,,\alpha^\star)$ over $\mathcal F_r$.
By \eqref{eq:fair-lagrangian}, this is equivalent to
\[
H^\star
\in
\arg\max_{H\in\mathcal F_r}
\operatorname{tr}\{H M_\lambda(\alpha^\star)\}.
\]
The strict eigengap in
\eqref{eq:fair-optimal-eigengap} implies, by the Ky Fan
variational principle, that the maximizer is unique and
equals $P_r\{M_\lambda(\alpha^\star)\}$. Since the same
argument applies to every relaxed minimizer,
the relaxed problem has a unique solution, and
\[
H_\lambda^\star
=
P_r\{M_\lambda(\alpha^\star)\}
\in\mathcal P_r.
\]

Because $\mathcal P_r\subset\mathcal F_r$, the relaxed
optimal value is no greater than the optimal value of
the projector-level problem. The relaxed optimizer
belongs to $\mathcal P_r$ and attains the relaxed
optimal value, so the two optimal values coincide.
Uniqueness over $\mathcal F_r$ also gives uniqueness
over $\mathcal P_r$.

Finally, the minimax equality and the preceding
sandwich argument hold for any
$\alpha\in\mathcal A_\lambda^\star$, which establishes
\eqref{eq:fair-common-primal-projector}. If such an
$\alpha$ has a strict rank-$r$ eigengap, its leading
spectral projector is unique and must equal
$P_\lambda^\star$.
\end{proof}

\section{Proofs for Section 3: Population Reconstruction--Fairness Trade-off}
\label{supp:sec:fair-tradeoff}
We first establish global properties of the
regularization path, including monotonicity
and the marginal reconstruction price of
fairness. We then develop the spectral
expansions and active-set arguments underlying
the local reconstruction--fairness analysis,
and characterize the conditions for a
nontrivial local fairness response.
\subsection{Global reconstruction--fairness properties}
\label{supp:subsec:fair-global}
\begin{proposition}[Monotonicity and finite-difference bounds]
\label{prop:fair-monotone-path}
Let $0\leq\lambda_1<\lambda_2$, and choose
\[
H_i\in
\arg\min_{H\in\mathcal F_r}
\{R(H)+\lambda_i B(H)\},
\qquad i=1,2.
\]
Then
\begin{equation}
B(H_2)\leq B(H_1),
\qquad
R(H_2)\geq R(H_1).
\label{eq:fair-path-monotonicity}
\end{equation}
Moreover,
\begin{equation}
\begin{aligned}
\lambda_1\{B(H_1)-B(H_2)\}
\leq R(H_2)-R(H_1)\leq
\lambda_2\{B(H_1)-B(H_2)\}.
\end{aligned}
\label{eq:fair-finite-difference-bound}
\end{equation}
\end{proposition}
\begin{proof}[Proof of Proposition~\ref{prop:fair-monotone-path}]
Let
\[
H_i\in
\arg\min_{H\in\mathcal F_r}
\{R(H)+\lambda_iB(H)\},
\qquad i=1,2,
\]
where $0\leq\lambda_1<\lambda_2$.
Optimality of $H_1$ at $\lambda_1$ gives
\[
R(H_1)+\lambda_1B(H_1)
\leq
R(H_2)+\lambda_1B(H_2),
\]
and hence
\begin{equation}
R(H_2)-R(H_1)
\geq
\lambda_1\{B(H_1)-B(H_2)\}.
\label{eq:supp:fair-monotone-lower}
\end{equation}
Similarly, optimality of $H_2$ at $\lambda_2$ yields
\begin{equation}
R(H_2)-R(H_1)
\leq
\lambda_2\{B(H_1)-B(H_2)\}.
\label{eq:supp:fair-monotone-upper}
\end{equation}
Combining the two inequalities,
\[
(\lambda_2-\lambda_1)
\{B(H_1)-B(H_2)\}\geq0.
\]
Since $\lambda_2>\lambda_1$, we obtain
\[
B(H_2)\leq B(H_1).
\]
Substitution into
\eqref{eq:supp:fair-monotone-lower} gives
$R(H_2)\geq R(H_1)$. The two-sided bound
\eqref{eq:fair-finite-difference-bound} follows
directly from
\eqref{eq:supp:fair-monotone-lower} and
\eqref{eq:supp:fair-monotone-upper}.
\end{proof}

\begin{proof}[Proof of Proposition~\ref{prop:fair-value-function}]
For every fixed $H\in\mathcal F_r$, the function
\[
\lambda\longmapsto R(H)+\lambda B(H)
\]
is affine. Since $V$ is the pointwise infimum
of these functions, it is concave on
$[0,\infty)$.

Moreover, $B(H)\geq0$ for every
$H\in\mathcal F_r$. Consequently, each
affine function is nondecreasing in $\lambda$,
and so is their pointwise infimum $V$.

Let $H_\lambda^\star$ be an optimizer at
$\lambda$. For every $\lambda'\geq0$,
\[
\begin{aligned}
V(\lambda')
&\leq
R(H_\lambda^\star)
+\lambda'B(H_\lambda^\star)=
V(\lambda)
+(\lambda'-\lambda)B(H_\lambda^\star).
\end{aligned}
\]
This proves the supergradient inequality
\eqref{eq:fair-value-supergradient}.

Because $V$ is finite and concave, it is
differentiable almost everywhere on
$(0,\infty)$. At a differentiability point,
its supergradient is unique, and therefore
\[
V'(\lambda)=B(H_\lambda^\star).
\]
By Proposition~\ref{prop:fair-monotone-path},
for any selection of optimizers the functions
\[
\lambda\longmapsto B(H_\lambda^\star),
\qquad
\lambda\longmapsto R(H_\lambda^\star)
\]
are monotone. Hence both are differentiable
almost everywhere.
At every point where all three functions
are differentiable, the identity
\[
V(\lambda)
=
R(H_\lambda^\star)
+
\lambda B(H_\lambda^\star)
\]
can be differentiated to give
\[
V'(\lambda)
=
\frac{d}{d\lambda}R(H_\lambda^\star)
+
B(H_\lambda^\star)
+
\lambda
\frac{d}{d\lambda}B(H_\lambda^\star).
\]
Using
$V'(\lambda)=B(H_\lambda^\star)$,
we obtain
\[
\frac{d}{d\lambda}R(H_\lambda^\star)
=
-\lambda
\frac{d}{d\lambda}B(H_\lambda^\star).
\]
This proves \eqref{eq:fair-marginal-price}
almost everywhere on $(0,\infty)$.
\end{proof}
\subsection{Auxiliary results for the local expansion}
\label{supp:subsec:fair-local-lemmas}

Throughout the local analysis, the dimension $p$,
the rank $r$, the number of groups $G$, and the
population covariance matrices are fixed as
$\lambda\downarrow0$.
Recall the spectral decomposition
\[
\overline\Sigma
=
\sum_{i=1}^p\mu_i u_i u_i^\top,
\qquad
\mu_1\geq\cdots\geq\mu_p,
\]
where $u_1,\ldots,u_p$ form an orthonormal basis
of eigenvectors. Under the ordinary-PCA eigengap
condition \eqref{eq:fair-ordinary-eigengap},
the leading rank-$r$ projector is uniquely given by
\(
P_0^\star
=
\sum_{i=1}^r u_i u_i^\top.
\)

We first establish a second-order expansion of the
Ky Fan functional that is uniform over the simplex.

\begin{lemma}[Uniform second-order Ky Fan expansion]
\label{lem:supp:fair-ky-fan-expansion}
Suppose that the ordinary-PCA eigengap satisfies
$\delta_0=\mu_r-\mu_{r+1}>0$.
Then, as $\lambda\downarrow0$,
\begin{equation}
\begin{aligned}
\mathcal K_r\{
\overline\Sigma+\lambda S(\alpha)
\}
=
\mathcal K_r(\overline\Sigma)
+
\lambda\operatorname{tr}
\{P_0^\star S(\alpha)\}+
\lambda^2\Gamma_0(\alpha)
+
o(\lambda^2),
\end{aligned}
\label{eq:supp:fair-ky-fan-expansion}
\end{equation}
uniformly over $\alpha\in\Delta_G$, where
$\Gamma_0(\alpha)$ is defined in
\eqref{eq:fair-local-spectral-coupling}.
\end{lemma}
\begin{proof}
Let
\[
m=\max_{1\leq a\leq G}
\|\Sigma_a\|_{\mathrm{op}}.
\]
Since $\alpha\in\Delta_G$,
\[
\|S(\alpha)\|_{\mathrm{op}}
\leq m.
\]
For
\[
A_\alpha(t)=\overline\Sigma+tS(\alpha),
\]
Weyl's inequality gives
\[
\begin{aligned}
\theta_r\{A_\alpha(t)\}
-
\theta_{r+1}\{A_\alpha(t)\}\geq
\delta_0-2|t|m.
\end{aligned}
\]
Consequently, there exists $t_0>0$ such that
\[
\theta_r\{A_\alpha(t)\}
-
\theta_{r+1}\{A_\alpha(t)\}
\geq\frac{\delta_0}{2}
\]
uniformly over $\alpha\in\Delta_G$ and
$|t|\leq t_0$.
Write
\[
P_\alpha(t)=P_r\{A_\alpha(t)\},
\qquad
f_\alpha(t)=\mathcal K_r\{A_\alpha(t)\}.
\]
The uniform separation of the leading spectral
cluster implies that $P_\alpha(t)$ and
$f_\alpha(t)$ are smooth in $t$ on this
neighborhood. In particular,
\[
f_\alpha'(t)
=
\operatorname{tr}\{P_\alpha(t)S(\alpha)\},
\]
so that
\begin{equation}
f_\alpha'(0)
=
\operatorname{tr}\{P_0^\star S(\alpha)\}.
\label{eq:supp:fair-ky-fan-first}
\end{equation}
The derivative of the leading spectral projector
at $t=0$ is
\begin{equation}
\begin{aligned}
P_\alpha'(0)
=
\sum_{i=1}^r\sum_{j=r+1}^p
\frac{u_j^\top S(\alpha)u_i}
{\mu_i-\mu_j}\times
(u_ju_i^\top+u_iu_j^\top).
\end{aligned}
\label{eq:supp:fair-projector-derivative}
\end{equation}
This expression remains valid when eigenvalues
within either spectral cluster have multiplicities,
because the two clusters are separated by
$\delta_0>0$.
Differentiating $f_\alpha'$ once more gives
\[
f_\alpha''(0)
=
\operatorname{tr}\{P_\alpha'(0)S(\alpha)\}.
\]
Substituting
\eqref{eq:supp:fair-projector-derivative}
and using the symmetry of $S(\alpha)$ yields
\[
\begin{aligned}
f_\alpha''(0)
&=
2\sum_{i=1}^r\sum_{j=r+1}^p
\frac{
\{u_j^\top S(\alpha)u_i\}^2
}{
\mu_i-\mu_j
}=2\Gamma_0(\alpha).
\end{aligned}
\]
Thus the second-order Taylor expansion
at $t=0$ has the asserted coefficients.
It remains to establish uniformity. The set
\(
\{S(\alpha):\alpha\in\Delta_G\}
\)
is compact, and the leading spectral cluster
of $A_\alpha(t)$ remains uniformly separated
for $|t|\leq t_0$. By the resolvent representation
of spectral projectors, the second derivative
$f_\alpha''(t)$ is jointly continuous in
$(t,\alpha)$ on this compact parameter set.
Consequently,
\[
\omega(h)
:=
\sup_{\substack{\alpha\in\Delta_G\\|t|\leq h}}
|f_\alpha''(t)-f_\alpha''(0)|
\longrightarrow0
\qquad\text{as }h\downarrow0.
\]
Taylor's theorem with integral remainder gives
\[
\begin{aligned}
f_\alpha(\lambda)-f_\alpha(0)
-\lambda f_\alpha'(0)
-\frac{\lambda^2}{2}f_\alpha''(0)=
\int_0^\lambda
(\lambda-t)
\{f_\alpha''(t)-f_\alpha''(0)\}\,dt.
\end{aligned}
\]
Hence, uniformly over $\alpha\in\Delta_G$,
\[
\begin{aligned}
\Big|
f_\alpha(\lambda)-f_\alpha(0)
-\lambda f_\alpha'(0)
-\frac{\lambda^2}{2}f_\alpha''(0)
\Big|
\leq
\frac{\lambda^2}{2}\omega(\lambda)
=o(\lambda^2).
\end{aligned}
\]
Substituting the expressions for
$f_\alpha'(0)$ and $f_\alpha''(0)$
establishes the asserted uniform expansion.
\end{proof}

The next lemma identifies the second-order
contribution of the active groups to the
optimized dual objective.

\begin{lemma}[Second-order active-set reduction]
\label{lem:supp:fair-active-reduction}
Suppose that the ordinary-PCA eigengap satisfies
$\delta_0=\mu_r-\mu_{r+1}>0$.
Define
\[
f(\alpha)
=
\sum_{a=1}^G
\alpha_a\mathcal E_a(P_0^\star),
\qquad
\alpha\in\Delta_G.
\]
Suppose that
$\varepsilon_\lambda:\Delta_G\to\mathbb R$
is continuous and satisfies
\[
\sup_{\alpha\in\Delta_G}
|\varepsilon_\lambda(\alpha)|
=o(\lambda^2).
\]
Then
\begin{equation}
\begin{aligned}
\max_{\alpha\in\Delta_G}
\big\{
\lambda f(\alpha)
-\lambda^2\Gamma_0(\alpha)
+\varepsilon_\lambda(\alpha)
\big\}
=
\lambda B(P_0^\star)-\lambda^2\zeta_0
+o(\lambda^2),
\end{aligned}
\label{eq:supp:fair-active-reduction}
\end{equation}
where
$\Gamma_0(\alpha)$ is defined in
\eqref{eq:fair-local-spectral-coupling}.
\end{lemma}

\begin{proof}
Write
\[
b_0=B(P_0^\star).
\]
By definition of the ordinary-PCA active set,
\[
\max_{\alpha\in\Delta_G}f(\alpha)=b_0,
\]
and the set of maximizers is precisely
$\Delta(\mathcal I_0)$.
Since $\Gamma_0$ is continuous and
$\Delta(\mathcal I_0)$ is compact,
there exists
\[
\bar\alpha\in
\arg\min_{\alpha\in\Delta(\mathcal I_0)}
\Gamma_0(\alpha).
\]
Thus
\[
f(\bar\alpha)=b_0,
\qquad
\Gamma_0(\bar\alpha)=\zeta_0.
\]
Evaluating the objective at $\bar\alpha$
gives the lower bound
\begin{equation}
\begin{aligned}
\max_{\alpha\in\Delta_G}
\big\{
\lambda f(\alpha)
-\lambda^2\Gamma_0(\alpha)
+\varepsilon_\lambda(\alpha)
\big\}
\geq
\lambda b_0-\lambda^2\zeta_0
+o(\lambda^2).
\end{aligned}
\label{eq:supp:fair-active-lower}
\end{equation}

For the upper bound, let $\alpha_\lambda$
be a maximizer. Existence follows from
continuity and compactness. Comparison
with $\bar\alpha$ gives
\[
\begin{aligned}
\lambda f(\alpha_\lambda)
-\lambda^2\Gamma_0(\alpha_\lambda)
+\varepsilon_\lambda(\alpha_\lambda)\geq
\lambda b_0-\lambda^2\zeta_0
+\varepsilon_\lambda(\bar\alpha).
\end{aligned}
\]
Rearranging,
\[
\begin{aligned}
\lambda\{b_0-f(\alpha_\lambda)\}
\leq{}&
\lambda^2\{
\zeta_0-\Gamma_0(\alpha_\lambda)
\}+
\varepsilon_\lambda(\alpha_\lambda)
-\varepsilon_\lambda(\bar\alpha).
\end{aligned}
\]
Since $\Gamma_0$ is bounded on the
compact simplex and the remainder is
uniformly $o(\lambda^2)$,
\[
0\leq b_0-f(\alpha_\lambda)=O(\lambda).
\]
In particular,
\[
f(\alpha_\lambda)\longrightarrow b_0.
\]
By compactness of $\Delta_G$ and continuity of $f$,
every accumulation point of $\alpha_\lambda$
belongs to
\(
\arg\max_{\alpha\in\Delta_G}f(\alpha)
=
\Delta(\mathcal I_0).
\)
Since $\Gamma_0$ is continuous and
$\zeta_0$ is its minimum over
$\Delta(\mathcal I_0)$, it follows that
\[
\liminf_{\lambda\downarrow0}
\Gamma_0(\alpha_\lambda)\geq\zeta_0.
\]
Equivalently,
\[
\Gamma_0(\alpha_\lambda)
\geq\zeta_0-o(1).
\]
Using $f(\alpha_\lambda)\leq b_0$,
we obtain
\[
\begin{aligned}
\lambda f(\alpha_\lambda)
-\lambda^2\Gamma_0(\alpha_\lambda)
+\varepsilon_\lambda(\alpha_\lambda)\leq
\lambda b_0
-\lambda^2\zeta_0
+o(\lambda^2).
\end{aligned}
\]
Combining this inequality with
\eqref{eq:supp:fair-active-lower}
establishes the result.
\end{proof}

\subsection{Proof of the local reconstruction--fairness expansion}
\label{supp:subsec:fair-local-theorem}

\begin{proof}[Proof of Theorem~\ref{thm:fair-local-expansion}]
Recall that
\(
M_\lambda(\alpha)
=
\overline\Sigma+\lambda S(\alpha).
\)
By the same Weyl inequality argument used in the
proof of Lemma~\ref{lem:supp:fair-ky-fan-expansion},
there exists $\lambda_{\mathrm{loc}}>0$ such that
\[
\theta_r\{M_\lambda(\alpha)\}
-
\theta_{r+1}\{M_\lambda(\alpha)\}
\geq\frac{\delta_0}{2}
\]
for every $\alpha\in\Delta_G$ and
$0\leq\lambda<\lambda_{\mathrm{loc}}$.

In particular, every dual-optimal effective
covariance has a strict rank-$r$ eigengap.
Theorem~\ref{thm:fair-spectral-characterization}
therefore implies that the relaxed optimizer
is unique and is a rank-$r$ orthogonal projector
throughout this interval.
We next derive the expansion of $V$.
Recall from \eqref{eq:fair-primal-dual-value}
that
\[
\begin{aligned}
V(\lambda)
=
\operatorname{tr}(\overline\Sigma)
+
\max_{\alpha\in\Delta_G}
\Big[
\lambda\sum_{a=1}^G\alpha_a c_a
-
\mathcal K_r\{
\overline\Sigma+\lambda S(\alpha)
\}
\Big].
\end{aligned}
\]
By Lemma~\ref{lem:supp:fair-ky-fan-expansion},
uniformly over $\alpha\in\Delta_G$,
\[
\begin{aligned}
\mathcal K_r\{
\overline\Sigma+\lambda S(\alpha)
\}
=
\mathcal K_r(\overline\Sigma)
+
\lambda\operatorname{tr}
\{P_0^\star S(\alpha)\}&+
\lambda^2\Gamma_0(\alpha)
+
o(\lambda^2).
\end{aligned}
\]
Observe that
\[
\begin{aligned}
\sum_{a=1}^G\alpha_a c_a
-
\operatorname{tr}\{P_0^\star S(\alpha)\}=
\sum_{a=1}^G
\alpha_a
\{c_a-\operatorname{tr}(P_0^\star\Sigma_a)\}=
\sum_{a=1}^G
\alpha_a\mathcal E_a(P_0^\star)
=f(\alpha),
\end{aligned}
\]
where $f$ is defined in
Lemma~\ref{lem:supp:fair-active-reduction}.
Consequently,
\[
\begin{aligned}
V(\lambda)
=
\operatorname{tr}(\overline\Sigma)
-\mathcal K_r(\overline\Sigma)+
\max_{\alpha\in\Delta_G}
\{
\lambda f(\alpha)
-\lambda^2\Gamma_0(\alpha)
+\varepsilon_\lambda(\alpha)
\},
\end{aligned}
\]
where the uniform Ky Fan expansion gives
\[
\sup_{\alpha\in\Delta_G}
|\varepsilon_\lambda(\alpha)|
=o(\lambda^2).
\]
Since
\[
\operatorname{tr}(\overline\Sigma)
-\mathcal K_r(\overline\Sigma)
=
R(P_0^\star),
\]
Lemma~\ref{lem:supp:fair-active-reduction}
yields
\[
V(\lambda)
=
R(P_0^\star)
+\lambda B(P_0^\star)
-\zeta_0\lambda^2
+o(\lambda^2).
\]
This proves
\eqref{eq:fair-local-value-expansion}.

We now derive the regret expansion.
For every sufficiently small $\lambda>0$,
the relaxed optimizer is unique. By Danskin's
theorem, $V$ is differentiable at such $\lambda$,
with
\begin{equation}
V'(\lambda)
=
B(H_\lambda^\star)
=
B(P_\lambda^\star).
\label{eq:supp:fair-local-envelope}
\end{equation}

To extract the derivative expansion from
\eqref{eq:fair-local-value-expansion},
we use concavity of $V$. Fix
$0<c<1<h$. For sufficiently small
$\lambda>0$,
\begin{equation}
\frac{
V(h\lambda)-V(\lambda)
}{
(h-1)\lambda
}
\leq
V'(\lambda)
\leq
\frac{
V(\lambda)-V(c\lambda)
}{
(1-c)\lambda
}.
\label{eq:supp:fair-derivative-sandwich}
\end{equation}
Write $b_0=B(P_0^\star)$. The value
expansion gives
\[
\begin{aligned}
\frac{
V(h\lambda)-V(\lambda)
}{
(h-1)\lambda
}
&=
b_0-\zeta_0(h+1)\lambda+o(\lambda),
\\
\frac{
V(\lambda)-V(c\lambda)
}{
(1-c)\lambda
}
&=
b_0-\zeta_0(1+c)\lambda+o(\lambda).
\end{aligned}
\]
It follows that
\[
\begin{aligned}
-\zeta_0(h+1)
&\leq
\liminf_{\lambda\downarrow0}
\frac{V'(\lambda)-b_0}{\lambda}\leq
\limsup_{\lambda\downarrow0}
\frac{V'(\lambda)-b_0}{\lambda}
\leq
-\zeta_0(1+c).
\end{aligned}
\]
Letting $h\downarrow1$ and
$c\uparrow1$, we obtain
\[
V'(\lambda)
=
b_0-2\zeta_0\lambda+o(\lambda).
\]
Together with
\eqref{eq:supp:fair-local-envelope},
this proves
\[
B(P_0^\star)-B(P_\lambda^\star)
=
2\zeta_0\lambda+o(\lambda),
\]
which is
\eqref{eq:fair-local-regret-expansion}.
Finally, by optimality,
\[
V(\lambda)
=
R(P_\lambda^\star)
+\lambda B(P_\lambda^\star).
\]
Therefore,
\[
\begin{aligned}
R(P_\lambda^\star)
&=
V(\lambda)-\lambda B(P_\lambda^\star)=
R(P_0^\star)
+\lambda b_0
-\zeta_0\lambda^2
+o(\lambda^2)-
\lambda\{b_0-2\zeta_0\lambda+o(\lambda)\}\\
&=
R(P_0^\star)
+\zeta_0\lambda^2
+o(\lambda^2).
\end{aligned}
\]
This establishes
\eqref{eq:fair-local-risk-expansion}
and completes the proof.
\end{proof}

\subsection{Local nondegeneracy and stability}
\label{supp:subsec:fair-nondegeneracy}
\begin{proof}[Proof of Proposition~\ref{prop:fair-local-nondegeneracy}]
By definition,
\[
\Gamma_0(\alpha)
=
\sum_{i=1}^r
\sum_{j=r+1}^p
\frac{
\{u_j^\top S(\alpha)u_i\}^2
}{
\mu_i-\mu_j
}.
\]
Under \eqref{eq:fair-ordinary-eigengap},
all denominators are strictly positive.
Hence
\[
\Gamma_0(\alpha)=0
\]
if and only if
\[
u_j^\top S(\alpha)u_i=0
\qquad
\text{for all }i\leq r<j.
\]
Since
\[
P_0^\star=\sum_{i=1}^r u_i u_i^\top,
\]
these conditions are equivalent to
\[
(I-P_0^\star)S(\alpha)P_0^\star=0.
\]
Since $\Gamma_0$ is continuous and nonnegative
and $\Delta(\mathcal I_0)$ is compact,
the minimum defining $\zeta_0$ is attained.
Consequently, $\zeta_0=0$ if and only if
the preceding matrix identity holds for some
$\alpha\in\Delta(\mathcal I_0)$.
This proves the first assertion.

We next establish the stability assertion.
Suppose $\zeta_0=0$, and choose
$\bar\alpha\in\Delta(\mathcal I_0)$
such that
\begin{equation}
(I-P_0^\star)S(\bar\alpha)P_0^\star=0.
\label{eq:supp:fair-vanishing-cross-block}
\end{equation}
Since $P_0^\star$ is a spectral projector
of $\overline\Sigma$, the latter is block
diagonal relative to the orthogonal
decomposition
\[
\mathbb R^p
=
\operatorname{ran}(P_0^\star)
\oplus
\operatorname{ran}(I-P_0^\star).
\]
Moreover, \eqref{eq:supp:fair-vanishing-cross-block}
and the symmetry of $S(\bar\alpha)$ imply
that $S(\bar\alpha)$ is block diagonal
relative to the same decomposition.
Consequently, so is
\[
M_\lambda(\bar\alpha)
=
\overline\Sigma+\lambda S(\bar\alpha).
\]

At $\lambda=0$, the leading and trailing
blocks are separated by the strict gap
$\delta_0>0$. By Weyl's inequality,
there exists $\lambda_{\mathrm{flat}}>0$
such that this separation persists for
every $0\leq\lambda<\lambda_{\mathrm{flat}}$.
It follows that
\begin{equation}
P_r\{M_\lambda(\bar\alpha)\}
=
P_0^\star,
\qquad
0\leq\lambda<\lambda_{\mathrm{flat}}.
\label{eq:supp:fair-stable-projector}
\end{equation}
By the definition of the active set,
\[
\mathcal E_a(P_0^\star)
=
B(P_0^\star),
\qquad a\in\mathcal I_0.
\]
Since $\bar\alpha$ is supported on
$\mathcal I_0$,
\[
\sum_{a=1}^G
\bar\alpha_a\mathcal E_a(P_0^\star)
=
B(P_0^\star).
\]
On the other hand, for any
$\alpha\in\Delta_G$,
\[
\sum_{a=1}^G
\alpha_a\mathcal E_a(P_0^\star)
\leq B(P_0^\star).
\]
Recalling that
\[
L_\lambda(H,\alpha)
=
R(H)+\lambda\sum_{a=1}^G
\alpha_a\mathcal E_a(H),
\]
we therefore obtain
\begin{equation}
L_\lambda(P_0^\star,\alpha)
\leq
L_\lambda(P_0^\star,\bar\alpha),
\qquad \alpha\in\Delta_G.
\label{eq:supp:fair-saddle-first}
\end{equation}
On the other hand,
\eqref{eq:supp:fair-stable-projector}
and the Ky Fan variational principle
give
\[
P_0^\star
\in
\arg\max_{H\in\mathcal F_r}
\operatorname{tr}\{
H M_\lambda(\bar\alpha)
\}.
\]
Using the representation
\eqref{eq:fair-lagrangian},
we obtain
\begin{equation}
L_\lambda(P_0^\star,\bar\alpha)
\leq
L_\lambda(H,\bar\alpha),
\qquad
H\in\mathcal F_r.
\label{eq:supp:fair-saddle-second}
\end{equation}
Combining
\eqref{eq:supp:fair-saddle-first}
and
\eqref{eq:supp:fair-saddle-second},
the pair $(P_0^\star,\bar\alpha)$
is a saddle point of $L_\lambda$.
Consequently, $P_0^\star$ is a
minimizer of the relaxed Fair PCA
objective, and $\bar\alpha$ is
dual-optimal.

Since $M_\lambda(\bar\alpha)$ has a
strict rank-$r$ eigengap throughout
$[0,\lambda_{\mathrm{flat}})$,
Theorem~\ref{thm:fair-spectral-characterization}
guarantees uniqueness of the relaxed
minimizer. Thus
\[
H_\lambda^\star
=
P_\lambda^\star
=
P_0^\star,
\qquad
0\leq\lambda<\lambda_{\mathrm{flat}}.
\]
This proves \eqref{eq:fair-local-flat-path}.
\end{proof}

\subsection{Examples of local degeneracy and nondegeneracy}
\label{supp:subsec:fair-local-examples}

We illustrate the distinction between
$\zeta_0=0$ and $\zeta_0>0$ through two
simple covariance configurations.
In both examples, let $p=2$, $r=1$,
$G=2$, and $q_1=q_2=1/2$.
Write $e_1=(1,0)^\top$ and
$e_2=(0,1)^\top$.

\medskip
\noindent\textbf{Example 1: Cancellation among active groups.}
Consider
\[
\Sigma_1=
\begin{pmatrix}
3&1\\
1&1
\end{pmatrix},
\qquad
\Sigma_2=
\begin{pmatrix}
3&-1\\
-1&1
\end{pmatrix}.
\]
The aggregate covariance is
\[
\overline\Sigma
=
\frac{\Sigma_1+\Sigma_2}{2}
=
\begin{pmatrix}
3&0\\
0&1
\end{pmatrix}.
\]
Hence $P_0^\star=e_1e_1^\top$,
with eigengap $\delta_0=2$.
Both groups have the same largest eigenvalue
$2+\sqrt{2}$, giving
\[
\mathcal E_1(P_0^\star)
=
\mathcal E_2(P_0^\star)
=
\sqrt{2}-1.
\]
Thus $\mathcal I_0=\{1,2\}$.
For $\alpha=(t,1-t)$, $0\leq t\leq1$,
the cross-subspace interaction is
\[
e_2^\top S(\alpha)e_1
=
t-(1-t)=2t-1.
\]
Consequently,
\[
\Gamma_0(\alpha)
=
\frac{(2t-1)^2}{2},
\qquad
\zeta_0
=
\min_{0\leq t\leq1}
\frac{(2t-1)^2}{2}
=0.
\]
Although the two individual covariance
matrices have nonzero cross-subspace
interactions, these interactions cancel
when $t=1/2$.

In fact, this example exhibits a stronger
property than local constancy.
For any $H\in\mathcal F_1$,
\[
\begin{aligned}
B(H)
&\geq
\frac{\mathcal E_1(H)+\mathcal E_2(H)}{2}=
2+\sqrt{2}
-\operatorname{tr}(H\overline\Sigma)\geq \sqrt{2}-1
=B(P_0^\star).
\end{aligned}
\]
Since $P_0^\star$ also uniquely minimizes
aggregate reconstruction risk, it is the
unique minimizer of
$R(H)+\lambda B(H)$ for every finite
$\lambda\geq0$.
Thus the population Fair PCA projector
remains equal to ordinary PCA throughout
the regularization path in this symmetric
example.

\medskip
\noindent\textbf{Example 2: Nondegeneracy with a unique active group.}
Now consider
\[
\Sigma_1=
\begin{pmatrix}
2&1\\
1&2
\end{pmatrix},
\qquad
\Sigma_2=
\begin{pmatrix}
4&-1\\
-1&1
\end{pmatrix}.
\]
Their aggregate covariance is
\[
\overline\Sigma
=
\begin{pmatrix}
3&0\\
0&3/2
\end{pmatrix},
\]
so again $P_0^\star=e_1e_1^\top$,
now with eigengap $\delta_0=3/2$.
The reconstruction regrets at ordinary
PCA are
\[
\mathcal E_1(P_0^\star)=1,
\qquad
\mathcal E_2(P_0^\star)
=
\frac{\sqrt{13}-3}{2}<1.
\]
Therefore $\mathcal I_0=\{1\}$,
and $\Delta(\mathcal I_0)$ contains only
$\alpha=(1,0)$.
Since
\[
e_2^\top\Sigma_1e_1=1,
\]
we obtain
\[
\zeta_0
=
\Gamma_0(1,0)
=
\frac{1}{3-3/2}
=
\frac23>0.
\]
For all sufficiently small positive
$\lambda$, group 1 remains the unique
worst-regret group at the population
Fair PCA solution.
The corresponding effective covariance is
\[
M_\lambda(1,0)
=
\begin{pmatrix}
3+2\lambda&\lambda\\
\lambda&3/2+2\lambda
\end{pmatrix}.
\]
Its leading unit eigenvector can be written as
\[
v_\lambda=
\begin{pmatrix}
\cos\theta_\lambda\\
\sin\theta_\lambda
\end{pmatrix},
\qquad
\theta_\lambda
=
\frac12\arctan\left(\frac{4\lambda}{3}\right).
\]
Hence $P_\lambda^\star
=v_\lambda v_\lambda^\top$, and
\[
\theta_\lambda
=
\frac23\lambda+o(\lambda).
\]
Thus the population-optimal subspace
changes for every sufficiently small
positive penalty.
For sufficiently small $\lambda>0$, group 1
remains the unique worst-regret group.
Since
\[
\mathcal E_1(P_\lambda^\star)
=
1-\sin(2\theta_\lambda),
\]
we obtain
\[
\begin{aligned}
B(P_0^\star)-B(P_\lambda^\star)
=
\sin(2\theta_\lambda)=
\frac{4\lambda}{\sqrt{9+16\lambda^2}}=
\frac43\lambda+O(\lambda^3).
\end{aligned}
\]
Similarly, direct calculation gives
\[
\begin{aligned}
R(P_\lambda^\star)-R(P_0^\star)=
\frac32\sin^2\theta_\lambda=
\frac34
\left(
1-\frac{3}{\sqrt{9+16\lambda^2}}
\right)&=
\frac23\lambda^2+O(\lambda^4).
\end{aligned}
\]
Thus, the exact calculations recover the
leading coefficients $2\zeta_0=4/3$
and $\zeta_0=2/3$ predicted by
Theorem~\ref{thm:fair-local-expansion}.

These examples highlight the importance of restricting the minimization defining $\zeta_0$ to active groups. In both cases, the cross-subspace interactions of the two group covariance matrices cancel under the aggregate weights $q_1=q_2=1/2$. In Example 1, both groups are active, so these weights are admissible in the minimization defining $\zeta_0$, yielding $\zeta_0=0$. In contrast, only group 1 is active in Example 2, preventing its cross-subspace interaction from being offset by that of group 2 and hence yielding $\zeta_0>0$. Thus, cancellation across all groups does not necessarily imply local degeneracy; the decisive condition is whether such cancellation can occur among the groups attaining the largest reconstruction regret at ordinary PCA.

\section{Proofs for Section 4: Statistical Theory for Fair PCA}
\label{supp:sec:fair-sample-theory}

This section contains proofs of the statistical
results in Section~\ref{sec:fair-sample-theory}.
Throughout, the dimension $p$, the target rank $r$,
the number of groups $G$, the group sample sizes,
and the population covariance matrices may depend
on $n$. As in the main text, this dependence is
suppressed unless needed for clarity.

\subsection{Covariance concentration and uniform
criterion convergence}
\label{supp:subsec:fair-sample-concentration}

We first recall the simultaneous covariance
concentration bound used in Section~
\ref{subsec:fair-sample-setup}.

\begin{lemma}[Simultaneous covariance concentration]
\label{lem:supp:fair-simultaneous-concentration}
Suppose Assumption~\ref{ass:fair-subgaussian}
holds. For $x\geq0$, define
\[
\begin{aligned}
\eta_n(x)
=
\max_{1\leq a\leq G}
\|\Sigma_a\|_{\mathrm{op}}
\Bigg\{
\sqrt{
\frac{
\mathfrak r_a+1+x+\log(2G)
}{
n_a
}
}
+
\frac{
\mathfrak r_a+1+x+\log(2G)
}{
n_a
}
\Bigg\}.
\end{aligned}
\]
Then there exists a constant $C_K>0$,
depending only on the sub-Gaussian constant
$K$, such that
\[
\mathbb P\{e_n>C_K\eta_n(x)\}
\leq e^{-x},
\]
where $e_n
=
\max_{1\leq a\leq G}
\|\widehat\Sigma_a-\Sigma_a\|_{\mathrm{op}}$ is defined in \eqref{eq:fair-covariance-error}.
In particular,
\[
e_n=O_p(\eta_n),
\]
where $\eta_n=\eta_n(0)$ is defined in
\eqref{eq:fair-covariance-rate}.
\end{lemma}

\begin{proof}
Under Assumption~\ref{ass:fair-subgaussian},
the standard effective-rank concentration
inequality for sample covariance matrices from \citet{vershynin2026highdimensional}
gives, for every $t\geq1$,
\[
\begin{aligned}
\mathbb P\Bigg[
\|\widehat\Sigma_a-\Sigma_a\|_{\mathrm{op}}
>
C_K\|\Sigma_a\|_{\mathrm{op}}
\left\{
\sqrt{\frac{\mathfrak r_a+t}{n_a}}
+
\frac{\mathfrak r_a+t}{n_a}
\right\}
\Bigg]
\leq2e^{-t}.
\end{aligned}
\]
The conclusion also holds when $\Sigma_a=0$,
in which case $X_{a,i}=0$ almost surely.
Taking
\[
t=1+x+\log(2G)
\]
and applying the union bound over
$a=1,\ldots,G$ yields
\[
\begin{aligned}
\mathbb P\{e_n>C_K\eta_n(x)\}
\leq
2G\exp\{-1-x-\log(2G)\}
=
e^{-1-x}
\leq e^{-x}.
\end{aligned}
\]
It remains to establish the stochastic-order
statement. For each fixed $x\geq0$,
\[
\mathfrak r_a+1+x+\log(2G)
\leq
(1+x)\{\mathfrak r_a+1+\log(2G)\}.
\]
Consequently,
\[
\eta_n(x)\leq(1+x)\eta_n.
\]
The asserted bound
$e_n=O_p(\eta_n)$ follows by choosing
$x$ sufficiently large but fixed.
\end{proof}
\begin{lemma}[Covariance concentration with estimated group means]
\label{lem:supp:fair-centered-covariance}
For each group $a$, suppose that
$X_{a,1},\ldots,X_{a,n_a}$ are independent
observations with population mean
$m_a=\mathbb E(X_{a,i})$ and covariance
matrix $\Sigma_a$. Assume that the
population-centered observations
$X_{a,i}-m_a$ satisfy the sub-Gaussian
condition in
Assumption~\ref{ass:fair-subgaussian},
with a common constant $K$.
Define
\[
\overline X_a
=
\frac{1}{n_a}
\sum_{i=1}^{n_a}X_{a,i},
\qquad
\widetilde\Sigma_a
=
\frac{1}{n_a}
\sum_{i=1}^{n_a}
(X_{a,i}-\overline X_a)
(X_{a,i}-\overline X_a)^\top.
\]
Let
\[
\widetilde e_n
=
\max_{1\leq a\leq G}
\|\widetilde\Sigma_a-\Sigma_a\|_{\mathrm{op}}.
\]
Then there exists a constant $C_K>0$,
depending only on the sub-Gaussian constant
$K$, such that, for every
$x\geq0$,
\[
\mathbb P\{
\widetilde e_n>C_K\eta_n(x)
\}
\leq 2e^{-x},
\]
where $\eta_n(x)$ is defined in
Lemma~\ref{lem:supp:fair-simultaneous-concentration}.
Consequently,
\[
\widetilde e_n=O_p(\eta_n).
\]
\end{lemma}

\begin{proof}
Let $Z_{a,i}=X_{a,i}-m_a$.
Since
\[
X_{a,i}-\overline X_a
=
Z_{a,i}-\overline Z_a,
\]
the sample-centered covariance matrix is
unchanged by subtracting the population
mean. Hence, without loss of generality,
we may assume $m_a=0$ throughout the proof.
Since
\[
\widetilde\Sigma_a
=
\widehat\Sigma_a
-
\overline X_a\overline X_a^\top,
\]
the triangle inequality gives
\[
\widetilde e_n
\leq
e_n+
\max_{1\leq a\leq G}
\|\overline X_a\|_2^2.
\]
It therefore suffices to control the
groupwise sample means simultaneously.

Under Assumption~\ref{ass:fair-subgaussian},
independence and the sub-Gaussian moment
generating function bound imply that
\[
\mathbb E\exp(u^\top\overline X_a)
\leq
\exp\left\{
\frac{C_K}{2n_a}
u^\top\Sigma_a u
\right\},
\qquad u\in\mathbb R^p.
\]
We use $C_K$ to denote a positive constant
depending only on $K$, whose value may
change from line to line.

For completeness, let $Z$ be an independent
standard Gaussian vector in $\mathbb R^p$.
For sufficiently small $s>0$, the Gaussian
moment generating function identity yields
\[
\begin{aligned}
\mathbb E
\exp\{s\|\overline X_a\|_2^2\}
&=
\mathbb E_Z
\mathbb E
\exp\{\sqrt{2s}\,Z^\top\overline X_a\}\leq
\mathbb E_Z
\exp\left\{
\frac{C_Ks}{n_a}
Z^\top\Sigma_a Z
\right\}.
\end{aligned}
\]
Applying the exponential Markov inequality
and the standard Gaussian quadratic-form
bound therefore gives, for every $t\geq0$,
\[
\begin{aligned}
\mathbb P\Bigg[
\|\overline X_a\|_2^2
>
\frac{C_K}{n_a}
\Big\{
\operatorname{tr}(\Sigma_a)
+
2\sqrt{
\operatorname{tr}(\Sigma_a^2)t
}
+
2\|\Sigma_a\|_{\mathrm{op}}t
\Big\}
\Bigg]
\leq e^{-t}.
\end{aligned}
\]
Using
\[
\operatorname{tr}(\Sigma_a)
=
\mathfrak r_a\|\Sigma_a\|_{\mathrm{op}},
\qquad
\operatorname{tr}(\Sigma_a^2)
\leq
\mathfrak r_a
\|\Sigma_a\|_{\mathrm{op}}^2,
\]
we obtain
\[
\mathbb P\left[
\|\overline X_a\|_2^2
>
C_K\|\Sigma_a\|_{\mathrm{op}}
\frac{\mathfrak r_a+t}{n_a}
\right]
\leq e^{-t}.
\]
The case $\Sigma_a=0$ is immediate, since
$X_{a,i}=0$ almost surely.
Taking $t=1+x+\log(2G)$ and applying
a union bound over the groups shows that
\[
\begin{aligned}
\mathbb P\Bigg[
\max_{1\leq a\leq G}
\|\overline X_a\|_2^2
>
C_K
\max_{1\leq a\leq G}
\|\Sigma_a\|_{\mathrm{op}}
\frac{
\mathfrak r_a+1+x+\log(2G)
}{n_a}
\Bigg]
\leq e^{-x}.
\end{aligned}
\]
The threshold on the left is bounded by
$C_K\eta_n(x)$.
Combining this bound with
Lemma~\ref{lem:supp:fair-simultaneous-concentration}
and adjusting $C_K$ gives
\[
\mathbb P\{
\widetilde e_n>C_K\eta_n(x)
\}
\leq2e^{-x}.
\]
As in the proof of the preceding lemma,
$\eta_n(x)\leq(1+x)\eta_n$ for fixed
$x\geq0$. Hence
$\widetilde e_n=O_p(\eta_n)$.
\end{proof}
We will repeatedly use the Ky Fan variational
principle
\[
\mathcal K_r(A)
=
\sum_{j=1}^r\theta_j(A)
=
\max_{H\in\mathcal F_r}\operatorname{tr}(HA)
\]
for symmetric matrices $A$. In particular,
\begin{equation}
|\mathcal K_r(A)-\mathcal K_r(B)|
\leq
r\|A-B\|_{\mathrm{op}}.
\label{eq:supp:fair-ky-fan-lipschitz}
\end{equation}
Indeed, the inequality follows by comparing
the two variational problems and using
$\|H\|_*=r$ for $H\in\mathcal F_r$.

\begin{proof}[Proof of Proposition~
\ref{prop:fair-uniform-sample-criterion}]
Recall that
\[
Q_\lambda(H)
=
-\operatorname{tr}(H\overline\Sigma)
+\lambda B(H)
\]
and
\[
\widehat Q_\lambda(H)
=
-\operatorname{tr}(H\widehat{\overline\Sigma})
+\lambda\widehat B(H).
\]
For any $H\in\mathcal F_r$,
\[
\begin{aligned}
|\widehat Q_\lambda(H)-Q_\lambda(H)|
\leq{}&
\left|
\operatorname{tr}
\{H(\widehat{\overline\Sigma}
-\overline\Sigma)\}
\right|
+
\lambda|\widehat B(H)-B(H)|.
\end{aligned}
\]
Since $H\in\mathcal F_r$,
\[
\|H\|_*=\operatorname{tr}(H)=r.
\]
Moreover,
\[
\widehat{\overline\Sigma}
-\overline\Sigma
=
\sum_{a=1}^G
q_a(\widehat\Sigma_a-\Sigma_a).
\]
The nonnegativity and normalization
of the group weights therefore give
\[
\|\widehat{\overline\Sigma}
-\overline\Sigma\|_{\mathrm{op}}
\leq e_n.
\]
Hence
\[
\left|
\operatorname{tr}
\{H(\widehat{\overline\Sigma}
-\overline\Sigma)\}
\right|
\leq re_n.
\]
For every group $a$,
\[
\begin{aligned}
|\widehat{\mathcal E}_a(H)
-\mathcal E_a(H)|
\leq
|\widehat c_a-c_a|
+
\left|
\operatorname{tr}
\{H(\widehat\Sigma_a-\Sigma_a)\}
\right|\leq
2r\|\widehat\Sigma_a-\Sigma_a\|_{\mathrm{op}}
\leq2re_n,
\end{aligned}
\]
where we used
\eqref{eq:supp:fair-ky-fan-lipschitz}.
The elementary inequality
\[
|\max_a x_a-\max_a y_a|
\leq\max_a|x_a-y_a|
\]
then implies
\[
|\widehat B(H)-B(H)|
\leq2re_n.
\]
Combining the preceding bounds yields
\[
|\widehat Q_\lambda(H)-Q_\lambda(H)|
\leq r(1+2\lambda)e_n.
\]
The bound is independent of $H$, proving
\eqref{eq:fair-uniform-criterion-bound}.
For the dual objectives, recall that
\[
\widehat\phi_\lambda(\alpha)
=
\lambda\sum_{a=1}^G
\alpha_a\widehat c_a
-
\mathcal K_r\{\widehat M_\lambda(\alpha)\}.
\]
Since $\alpha\in\Delta_G$,
\[
\begin{aligned}
\|\widehat M_\lambda(\alpha)
-M_\lambda(\alpha)\|_{\mathrm{op}}
\leq
\|\widehat{\overline\Sigma}
-\overline\Sigma\|_{\mathrm{op}}
+
\lambda\sum_{a=1}^G
\alpha_a
\|\widehat\Sigma_a-\Sigma_a\|_{\mathrm{op}}
\leq(1+\lambda)e_n.
\end{aligned}
\]
Therefore,
\[
\begin{aligned}
|\widehat\phi_\lambda(\alpha)
-\phi_\lambda(\alpha)|
&\leq
\lambda
\sum_{a=1}^G
\alpha_a|\widehat c_a-c_a|+
|\mathcal K_r\{\widehat M_\lambda(\alpha)\}
-\mathcal K_r\{M_\lambda(\alpha)\}|
\leq
\lambda re_n+r(1+\lambda)e_n\\
&=
r(1+2\lambda)e_n.
\end{aligned}
\]
Taking the supremum over $\Delta_G$
establishes the dual bound.
Finally,
Lemma~\ref{lem:supp:fair-simultaneous-concentration}
gives $e_n=O_p(\eta_n)$, completing
the proof of the stochastic statements.
\end{proof}

\subsection{Spectral margin and basic projector consistency}
\label{supp:subsec:fair-basic-consistency}

The basic projector rate follows from the
quadratic margin of a linear spectral objective
over the Fantope.

\begin{lemma}[Spectral quadratic margin]
\label{lem:supp:fair-spectral-margin}
Let $A$ be a real symmetric $p\times p$ matrix,
with ordered eigenvalues
$\nu_1\geq\cdots\geq\nu_p$, and suppose
\[
g=\nu_r-\nu_{r+1}>0.
\]
Let $P=P_r(A)$. Then, for every
$H\in\mathcal F_r$,
\begin{equation}
\operatorname{tr}\{A(P-H)\}
\geq
\frac g2\|P-H\|_F^2.
\label{eq:supp:fair-spectral-margin}
\end{equation}
\end{lemma}

\begin{proof}
Let $v_1,\ldots,v_p$ be an orthonormal
eigenbasis of $A$, ordered according to
$\nu_1\geq\cdots\geq\nu_p$. Write
\[
P=\sum_{i=1}^r v_iv_i^\top
\]
and define
\[
s(H)
=
\sum_{i=1}^r
\{1-v_i^\top Hv_i\}.
\]
Since $\operatorname{tr}(H)=r$,
\[
s(H)
=
\sum_{j=r+1}^p v_j^\top Hv_j.
\]
The quantities $v_i^\top H v_i$ are the
diagonal entries of $H$ expressed in the
eigenbasis of $A$.
Furthermore, because $0\preceq H\preceq I$,
all the diagonal quantities
$v_i^\top Hv_i$ belong to $[0,1]$.
Consequently,
\[
\begin{aligned}
\operatorname{tr}\{A(P-H)\}
=
\sum_{i=1}^r
\nu_i\{1-v_i^\top Hv_i\}
-
\sum_{j=r+1}^p
\nu_jv_j^\top Hv_j
\geq
(\nu_r-\nu_{r+1})s(H)
=gs(H).
\end{aligned}
\]
Since $H^2\preceq H$ and
$\operatorname{tr}(H)=r$,
\[
\begin{aligned}
\|P-H\|_F^2
=
\operatorname{tr}(P^2)
+\operatorname{tr}(H^2)
-2\operatorname{tr}(PH)
\leq
2r-2\operatorname{tr}(PH)
=2s(H).
\end{aligned}
\]
Combining the last two inequalities
proves \eqref{eq:supp:fair-spectral-margin}.
\end{proof}

\begin{proof}[Proof of Theorem~
\ref{thm:fair-basic-projector-rate}]
Let
$\alpha_\lambda^\star\in
\mathcal A_\lambda^\star$
satisfy the eigengap condition in the theorem.
By Theorem~
\ref{thm:fair-spectral-characterization},
the population relaxed oracle is uniquely
given by
\[
P_\lambda^\star
=
P_r\{M_\lambda(\alpha_\lambda^\star)\}.
\]
Recall the Lagrangian
$L_\lambda(H,\alpha)$ from
\eqref{eq:fair-lagrangian}. Since
\[
Q_\lambda(H)
=
\max_{\alpha\in\Delta_G}
L_\lambda(H,\alpha)
-
\operatorname{tr}(\overline\Sigma),
\]
the additive constants cancel when taking
differences. Since $P_\lambda^\star$ and
$\alpha_\lambda^\star$ are primal and dual
optimal, respectively, the minimax equality
\eqref{eq:fair-minimax-equality}
implies that they form a saddle point.
In particular,
\[
L_\lambda(P_\lambda^\star,
\alpha_\lambda^\star)
=
\max_{\alpha\in\Delta_G}
L_\lambda(P_\lambda^\star,\alpha).
\]
Consequently,
for every $H\in\mathcal F_r$,
\[
\begin{aligned}
Q_\lambda(H)-Q_\lambda(P_\lambda^\star)
=
\max_{\alpha\in\Delta_G}
L_\lambda(H,\alpha)
-
\max_{\alpha\in\Delta_G}
L_\lambda(P_\lambda^\star,\alpha)
&\geq
L_\lambda(H,\alpha_\lambda^\star)
-
L_\lambda(P_\lambda^\star,
\alpha_\lambda^\star)\\
&=
\operatorname{tr}
\left[
M_\lambda(\alpha_\lambda^\star)
(P_\lambda^\star-H)
\right].
\end{aligned}
\]
Applying
Lemma~\ref{lem:supp:fair-spectral-margin}
gives
\begin{equation}
Q_\lambda(H)-Q_\lambda(P_\lambda^\star)
\geq
\frac{g_\lambda}{2}
\|H-P_\lambda^\star\|_F^2.
\label{eq:supp:fair-population-margin}
\end{equation}
Now let $\widehat H_\lambda$ be any
empirical optimizer. By its optimality,
\[
\widehat Q_\lambda(\widehat H_\lambda)
\leq
\widehat Q_\lambda(P_\lambda^\star).
\]
Therefore,
\[
\begin{aligned}
Q_\lambda(\widehat H_\lambda)
-Q_\lambda(P_\lambda^\star)
\leq
2\sup_{H\in\mathcal F_r}
|\widehat Q_\lambda(H)-Q_\lambda(H)|
\leq
2r(1+2\lambda)e_n,
\end{aligned}
\]
where the last inequality follows from
Proposition~\ref{prop:fair-uniform-sample-criterion}.
Combining this with
\eqref{eq:supp:fair-population-margin}
yields
\[
\frac{g_\lambda}{2}
\|\widehat H_\lambda-P_\lambda^\star\|_F^2
\leq
2r(1+2\lambda)e_n.
\]
Rearranging proves
\eqref{eq:fair-basic-projector-bound}.
Under Assumption~\ref{ass:fair-subgaussian},
we have $e_n=O_p(\eta_n)$.
Hence
\[
\|\widehat H_\lambda-P_\lambda^\star\|_F
=
O_p\left\{
\sqrt{
\frac{r(1+\lambda)\eta_n}{g_\lambda}
}
\right\}.
\]
The consistency assertion follows
under the stated rate condition.
\end{proof}

\subsection{Geometry of the dual-optimal set}
\label{supp:subsec:fair-dual-geometry}

We next establish the structural result underlying
the sharper projector rate.
Recall the dual objective $\phi_\lambda$
and the dual-optimal set $\mathcal A_\lambda^\star$
defined in \eqref{eq:fair-dual-objective}
and \eqref{eq:fair-dual-optimal-set},
respectively, and write
\[
\phi_\lambda^\star
=
\max_{\alpha\in\Delta_G}\phi_\lambda(\alpha).
\]
We first establish an exact decomposition
of the dual optimality gap.

\begin{lemma}[Exact dual-gap decomposition]
\label{lem:supp:fair-dual-gap-decomposition}
Suppose Assumption~\ref{ass:fair-optimal-eigengap}
holds.
Then, for every $\alpha\in\Delta_G$,
\[
\begin{aligned}
\phi_\lambda^\star-\phi_\lambda(\alpha)
=
\lambda\sum_{a=1}^G
\alpha_a
\{B(P_\lambda^\star)-\mathcal E_a(P_\lambda^\star)\}
+
\mathcal K_r\{M_\lambda(\alpha)\}
-
\operatorname{tr}\{P_\lambda^\star M_\lambda(\alpha)\}.
\end{aligned}
\]
Both terms on the right-hand side
are nonnegative.
\end{lemma}

\begin{proof}
By population primal--dual equality,
\[
\begin{aligned}
\phi_\lambda^\star
=
V(\lambda)-\operatorname{tr}(\overline\Sigma)
=
-\operatorname{tr}(P_\lambda^\star\overline\Sigma)
+\lambda B(P_\lambda^\star).
\end{aligned}
\]
On the other hand,
\[
\begin{aligned}
\phi_\lambda(\alpha)
=
\lambda\sum_{a=1}^G\alpha_a c_a
-
\mathcal K_r\{M_\lambda(\alpha)\}
=
-\operatorname{tr}(P_\lambda^\star\overline\Sigma)
+
\lambda\sum_{a=1}^G
\alpha_a\mathcal E_a(P_\lambda^\star)
-
\left[
\mathcal K_r\{M_\lambda(\alpha)\}
-
\operatorname{tr}\{P_\lambda^\star M_\lambda(\alpha)\}
\right].
\end{aligned}
\]
Subtracting the last identity from
the expression for $\phi_\lambda^\star$
gives the decomposition.

The first term is nonnegative because
$\mathcal E_a(P_\lambda^\star)\leq B(P_\lambda^\star)$.
The second is nonnegative by the
Ky Fan variational principle.
\end{proof}

The next lemma quantifies the spectral
optimality loss appearing in
\eqref{eq:fair-dual-gap-decomposition}.

\begin{lemma}[Spectral loss and cross-block mixing]
\label{lem:supp:fair-crossblock-mixing}
Let $A$ be a real symmetric $p\times p$ matrix
with $\|A\|_{\mathrm{op}}\leq L$
for some $L>0$, and let $P\in\mathcal P_r$.
Then
\begin{equation}
\mathcal K_r(A)-\operatorname{tr}(PA)
\geq
\frac{
\|(I-P)AP\|_F^2
}{
4L
}.
\label{eq:supp:fair-crossblock-mixing}
\end{equation}
The numerical constant is independent
of $p$ and $r$.
\end{lemma}

\begin{proof}
Choose an orthonormal basis whose first $r$
vectors span the range of $P$ and whose
remaining vectors span its orthogonal
complement. In this basis, write
\[
P=
\begin{pmatrix}
I_r&0\\
0&0
\end{pmatrix},
\qquad
A=
\begin{pmatrix}
A_{11}&X^\top\\
X&A_{22}
\end{pmatrix}.
\]
Then
\[
\|X\|_F=\|(I-P)AP\|_F.
\]
Define the skew-symmetric matrix
\[
K=
\begin{pmatrix}
0&-X^\top\\
X&0
\end{pmatrix}
\]
and the family of rank-$r$ projectors
\[
P_t=e^{tK}Pe^{-tK},
\qquad t\geq0.
\]
Let
\[
g(t)=\operatorname{tr}(P_tA).
\]
Differentiating $P_t=e^{tK}Pe^{-tK}$
with respect to $t$ gives
\[
P_t'=KP_t-P_tK
\]
and
\[
P_t''
=
K^2P_t-2KP_tK+P_tK^2.
\]
Consequently,
\[
\begin{aligned}
g'(0)
=
\operatorname{tr}\{(KP-PK)A\}
=
2\|X\|_F^2.
\end{aligned}
\]
Because $P_t$ is an orthogonal projector,
the Schatten-norm inequalities give
\[
\|K^2P_t\|_*
\leq\|K\|_F^2,
\qquad
\|P_tK^2\|_*
\leq\|K\|_F^2,
\]
and
\[
\|KP_tK\|_*
\leq\|K\|_F^2.
\]
Consequently,
\[
\|P_t''\|_*
\leq4\|K\|_F^2
=
8\|X\|_F^2.
\]
It follows that
\[
|g''(t)|
\leq
\|A\|_{\mathrm{op}}\|P_t''\|_*
\leq
8L\|X\|_F^2.
\]
Taylor's theorem therefore yields
\[
\begin{aligned}
g(t)-g(0)
\geq
tg'(0)-4Lt^2\|X\|_F^2
=
(2t-4Lt^2)\|X\|_F^2.
\end{aligned}
\]
Taking $t=(4L)^{-1}$ gives
\[
g(t)-g(0)
\geq
\frac{\|X\|_F^2}{4L}.
\]
Finally, $P_t\in\mathcal P_r$,
so the Ky Fan variational principle gives
\[
\mathcal K_r(A)
\geq
\operatorname{tr}(P_tA)
=
g(t).
\]
Since $g(0)=\operatorname{tr}(PA)$,
the conclusion follows.
\end{proof}

We now prove the polyhedral characterization
of the population dual-optimal set. Recall the active-group set
$\mathcal I_\lambda^\star$ and the cross-block
map $T_\lambda$ defined in
Section~\ref{subsec:fair-sharp-rate}.

\begin{lemma}[Polyhedral identification]
\label{lem:supp:fair-polyhedral-identification}
Suppose $\lambda>0$ and
Assumption~\ref{ass:fair-uniform-optimal-gap}
holds. Then
\[
\begin{aligned}
\mathcal A_\lambda^\star
=
\left\{
\alpha\in\Delta_G:
\alpha_a=0\ \text{ for all }a\notin\mathcal I_\lambda^\star,  T_\lambda(\alpha)=0
\right\}.
\end{aligned}
\]
\end{lemma}

\begin{proof}
Let $\mathcal C$ denote the set on the
right-hand side. Since
$T_\lambda(\alpha)$ is affine in $\alpha$,
$\mathcal C$ is a compact convex polyhedron.
First, let $\alpha\in\mathcal A_\lambda^\star$.
Since $\alpha$ is dual optimal,
$\phi_\lambda^\star-\phi_\lambda(\alpha)=0$.
By Lemma~\ref{lem:supp:fair-dual-gap-decomposition},
this gap is the sum of two nonnegative terms,
so both must vanish.
Because $\lambda>0$,
\[
\sum_{a=1}^G
\alpha_a\{B(P_\lambda^\star)-\mathcal E_a(P_\lambda^\star)\}=0
\]
implies
\[
\alpha_a=0
\text{ for all }a\notin\mathcal I_\lambda^\star.
\]
Moreover,
\[
\mathcal K_r\{M_\lambda(\alpha)\}
=
\operatorname{tr}\{P_\lambda^\star M_\lambda(\alpha)\}.
\]
Thus $P_\lambda^\star$ is a leading rank-$r$
spectral projector of $M_\lambda(\alpha)$.
In particular,
\[
(I-P_\lambda^\star)M_\lambda(\alpha)P_\lambda^\star=0,
\]
or equivalently,
$T_\lambda(\alpha)=0$.
Therefore,
\[
\mathcal A_\lambda^\star\subseteq\mathcal C.
\]
For the reverse inclusion, define
\[
\mathcal S
=
\left\{
\alpha\in\mathcal C:
\mathcal K_r\{M_\lambda(\alpha)\}
=
\operatorname{tr}\{P_\lambda^\star M_\lambda(\alpha)\}
\right\}.
\]
For every $\alpha\in\mathcal C$,
the definition of $\mathcal I_\lambda^\star$
implies
\[
\sum_{a=1}^G
\alpha_a
\{B(P_\lambda^\star)
-\mathcal E_a(P_\lambda^\star)\}
=0.
\]
Hence, by
Lemma~\ref{lem:supp:fair-dual-gap-decomposition},
\[
\phi_\lambda^\star-\phi_\lambda(\alpha)
=
\mathcal K_r\{M_\lambda(\alpha)\}
-
\operatorname{tr}\{P_\lambda^\star M_\lambda(\alpha)\},
\qquad \alpha\in\mathcal C.
\]
It follows that
$\mathcal S=\mathcal A_\lambda^\star$,
where we have used the previously established
inclusion
$\mathcal A_\lambda^\star\subseteq\mathcal C$.
In particular, $\mathcal S$ is nonempty.
Continuity of the Ky Fan functional
and the effective covariance map implies
that $\mathcal S$ is relatively closed
in $\mathcal C$.

We next show that $\mathcal S$ is
relatively open in $\mathcal C$.
Fix $\alpha\in\mathcal S$.
Then $\alpha\in\mathcal A_\lambda^\star$,
and Assumption
\ref{ass:fair-uniform-optimal-gap}
gives
\[
\operatorname{gap}_r
\{M_\lambda(\alpha)\}
\geq2\delta_\lambda>0.
\]
Furthermore, $P_\lambda^\star$ is the leading
rank-$r$ spectral projector of
$M_\lambda(\alpha)$.
For every $\beta\in\mathcal C$,
the identity
\[
T_\lambda(\beta)=0
\]
and the symmetry of
$M_\lambda(\beta)$ imply that
$M_\lambda(\beta)$ is block diagonal
relative to
\[
\mathbb R^p
=
\operatorname{ran}(P_\lambda^\star)
\oplus
\operatorname{ran}(I-P_\lambda^\star).
\]
At $\beta=\alpha$, Assumption
\ref{ass:fair-uniform-optimal-gap}
and the identity
$P_\lambda^\star=P_r\{M_\lambda(\alpha)\}$
imply that the smallest eigenvalue
of the first block exceeds the largest
eigenvalue of the second block by at least
$2\delta_\lambda>0$.
Since both blocks depend continuously
on $\beta$, Weyl's inequality ensures
that this strict separation persists
for all $\beta\in\mathcal C$
sufficiently close to $\alpha$.
As the eigenvalues of a block-diagonal
matrix are the union of the eigenvalues
of its diagonal blocks, the first block
continues to contain the largest $r$
eigenvalues of $M_\lambda(\beta)$.

For such $\beta$, $P_\lambda^\star$ remains
the leading rank-$r$ spectral
projector of $M_\lambda(\beta)$.
Consequently,
\[
\mathcal K_r\{M_\lambda(\beta)\}
=
\operatorname{tr}\{P_\lambda^\star M_\lambda(\beta)\},
\]
so $\beta\in\mathcal S$.

This proves that $\mathcal S$ is
relatively open in $\mathcal C$.
We have therefore shown that
$\mathcal S$ is nonempty and both
relatively open and relatively closed
in $\mathcal C$.
Because $\mathcal C$ is convex,
it is connected. Hence
\[
\mathcal S=\mathcal C.
\]
Since $\mathcal S=\mathcal A_\lambda^\star$,
the desired characterization follows.
\end{proof}

\begin{proof}[Proof of Theorem~
\ref{thm:fair-derived-dual-growth}]
The polyhedral characterization
\eqref{eq:fair-polyhedral-dual-set}
follows from
Lemma~\ref{lem:supp:fair-polyhedral-identification}.
It remains to establish the
quadratic-growth inequality.
Recall
\[
\sigma
=
\max_{1\leq a\leq G}
\|\Sigma_a\|_{\mathrm{op}}.
\]
Because the covariance matrices are
positive semidefinite,
\[
\begin{aligned}
\|M_\lambda(\alpha)\|_{\mathrm{op}}
\leq
\|\overline\Sigma\|_{\mathrm{op}}
+
\lambda\sum_{a=1}^G
\alpha_a\|\Sigma_a\|_{\mathrm{op}}
\leq(1+\lambda)\sigma.
\end{aligned}
\]
The assumed strict eigengap ensures
$\sigma>0$.
Applying
Lemma~\ref{lem:supp:fair-crossblock-mixing}
with
\[
A=M_\lambda(\alpha),
\qquad
P=P_\lambda^\star,
\qquad
L=(1+\lambda)\sigma
\]
yields
\[
\begin{aligned}
\mathcal K_r\{M_\lambda(\alpha)\}
-
\operatorname{tr}
\{P_\lambda^\star M_\lambda(\alpha)\}
\geq
\frac{
\|T_\lambda(\alpha)\|_F^2
}{
4(1+\lambda)\sigma
}.
\end{aligned}
\]
Combining this with
Lemma~\ref{lem:supp:fair-dual-gap-decomposition}
gives
\begin{equation}
\begin{aligned}
\phi_\lambda^\star-\phi_\lambda(\alpha)
\geq{}&
\lambda\sum_{a=1}^G
\alpha_a
\{B(P_\lambda^\star)-\mathcal E_a(P_\lambda^\star)\}
+
\frac{
\|T_\lambda(\alpha)\|_F^2
}{
4(1+\lambda)\sigma
}.
\end{aligned}
\label{eq:supp:fair-dual-gap-lower}
\end{equation}
Define
\[
t(\alpha)
=
\sum_{a\notin\mathcal I_\lambda^\star}
\alpha_a.
\]
If every group is active, we interpret
$t(\alpha)=0$.
Suppose first that inactive groups exist.
Recall that
\[
d_\lambda
=
\min_{a\notin\mathcal I_\lambda^\star}
\{B(P_\lambda^\star)
-\mathcal E_a(P_\lambda^\star)\}>0.
\]
Then
\[
\sum_{a=1}^G
\alpha_a\{B(P_\lambda^\star)-\mathcal E_a(P_\lambda^\star)\}
\geq
d_\lambda t(\alpha).
\]
Therefore,
\[
\begin{aligned}
\phi_\lambda^\star-\phi_\lambda(\alpha)
\geq{}&
\lambda d_\lambda t(\alpha)
+
\frac{
\|T_\lambda(\alpha)\|_F^2
}{
4(1+\lambda)\sigma
}.
\end{aligned}
\]
Since $\alpha\in\Delta_G$,
\[
0\leq t(\alpha)\leq1,
\]
and hence
\[
t(\alpha)\geq t(\alpha)^2.
\]
Let
\[
m_\lambda
=
\min\left\{
\lambda d_\lambda,
\frac{1}{4(1+\lambda)\sigma}
\right\}.
\]
It follows that
\begin{equation}
\begin{aligned}
\phi_\lambda^\star-\phi_\lambda(\alpha)
\geq
m_\lambda
\left\{
t(\alpha)^2+
\|T_\lambda(\alpha)\|_F^2
\right\}.
\end{aligned}
\label{eq:supp:fair-growth-components}
\end{equation}
If every group is active, the same
inequality holds with
\[
t(\alpha)=0,
\qquad
m_\lambda=\frac{1}{4(1+\lambda)\sigma}.
\]
By the polyhedral characterization
in Lemma~\ref{lem:supp:fair-polyhedral-identification},
$\mathcal A_\lambda^\star$ is the
feasible set of a finite system
of linear equalities and inequalities.
Hoffman's polyhedral error bound
\cite{Hoffman1952,LeventhalLewis2010}
therefore yields a finite constant
$\kappa_\lambda>0$ such that
\[
\begin{aligned}
\operatorname{dist}_2
(\alpha,\mathcal A_\lambda^\star)
\leq
\kappa_\lambda
\{t(\alpha)+\|T_\lambda(\alpha)\|_F\},
\end{aligned}
\]
for every $\alpha\in\Delta_G$.
Consequently,
\[
\begin{aligned}
\operatorname{dist}_2^2
(\alpha,\mathcal A_\lambda^\star)
\leq
\kappa_\lambda^2
\{t(\alpha)+\|T_\lambda(\alpha)\|_F\}^2
\leq
2\kappa_\lambda^2
\{t(\alpha)^2+\|T_\lambda(\alpha)\|_F^2\}.
\end{aligned}
\]
Combining this inequality with
\eqref{eq:supp:fair-growth-components}
gives
\[
\begin{aligned}
\phi_\lambda^\star-\phi_\lambda(\alpha)
&\geq
\frac{m_\lambda}{2\kappa_\lambda^2}
\operatorname{dist}_2^2
(\alpha,\mathcal A_\lambda^\star).
\end{aligned}
\]
Taking
\[
\underline\mu_\lambda
=
\frac{m_\lambda}{\kappa_\lambda^2}
\]
yields exactly
\eqref{eq:fair-derived-growth-modulus}
and proves
\eqref{eq:fair-derived-dual-growth}.
\end{proof}

\begin{remark}[Dependence of the growth modulus]
\label{remark:supp:fair-growth-modulus}
The preceding argument applies to every
finite-dimensional member of a sequence
of population models. Neither uniqueness
of the dual optimizer nor a negative-definite
dual Hessian is required.

The resulting modulus
$\underline\mu_\lambda$ may nevertheless
depend on $n$ through the covariance matrices,
the active-group regret gaps, the penalty,
and the Hoffman constant. In particular,
the theorem does not imply that
$\underline\mu_\lambda$ is bounded away
from zero uniformly over increasing-dimensional
models.

The sharp stochastic result below accounts
for this dependence explicitly.
\end{remark}

\subsection{Quantitative dual localization and
first-order projector recovery}
\label{supp:subsec:fair-sharp-rate}
We now develop the quantitative dual
localization and spectral perturbation
bounds required for
Theorem~\ref{thm:fair-sharp-projector-rate}.
The central issue is that convergence of
the empirical dual objective alone does
not automatically provide localization
uniformly along a sequence of changing
population problems. We therefore use
the quantitative quadratic-growth bound
derived above.
Recall
\[
\operatorname{dist}_2
(\alpha,\mathcal A_\lambda^\star)
=
\inf_{\beta\in\mathcal A_\lambda^\star}
\|\alpha-\beta\|_2.
\]
Since $\mathcal A_\lambda^\star$ is nonempty,
compact, and convex, its Euclidean projection
is uniquely defined.

\begin{lemma}[Quantitative dual localization]
\label{lem:supp:fair-dual-localization}
Suppose $\lambda>0$ and the quadratic-growth
inequality \eqref{eq:fair-derived-dual-growth}
holds with $\underline\mu_\lambda>0$.
Let
\[
E_n
=
\sup_{\alpha\in\Delta_G}
|\widehat\phi_\lambda(\alpha)
-\phi_\lambda(\alpha)|.
\]
Then every empirical dual optimizer
$\widehat\alpha_\lambda$ satisfies
\begin{equation}
\operatorname{dist}_2
(\widehat\alpha_\lambda,
\mathcal A_\lambda^\star)
\leq
2\sqrt{
\frac{E_n}{\underline\mu_\lambda}
}.
\label{eq:supp:fair-dual-localization}
\end{equation}
In particular, recalling the covariance
estimation error $e_n$ defined in
\eqref{eq:fair-covariance-error},
\[
\operatorname{dist}_2
(\widehat\alpha_\lambda,
\mathcal A_\lambda^\star)
\leq
2
\sqrt{
\frac{
r(1+2\lambda)e_n
}{
\underline\mu_\lambda
}
}.
\]
\end{lemma}

\begin{proof}
Let
$\widehat\alpha_\lambda\in
\widehat{\mathcal A}_\lambda$
and choose any
$\alpha^\star\in\mathcal A_\lambda^\star$.
By empirical optimality,
\[
\widehat\phi_\lambda
(\widehat\alpha_\lambda)
\geq
\widehat\phi_\lambda(\alpha^\star).
\]
Consequently,
\[
\begin{aligned}
\phi_\lambda^\star
-\phi_\lambda(\widehat\alpha_\lambda)
=
\phi_\lambda(\alpha^\star)
-\phi_\lambda(\widehat\alpha_\lambda)
\leq
|\phi_\lambda(\alpha^\star)
-\widehat\phi_\lambda(\alpha^\star)|
+
|\widehat\phi_\lambda(\widehat\alpha_\lambda)
-\phi_\lambda(\widehat\alpha_\lambda)|
\leq2E_n.
\end{aligned}
\]
The derived quadratic-growth inequality
\eqref{eq:fair-derived-dual-growth} gives
\[
\frac{\underline\mu_\lambda}{2}
\operatorname{dist}_2^2
(\widehat\alpha_\lambda,
\mathcal A_\lambda^\star)
\leq2E_n.
\]
Rearranging proves
\eqref{eq:supp:fair-dual-localization}.
The second assertion follows from the
dual bound in Proposition~
\ref{prop:fair-uniform-sample-criterion}.
\end{proof}

We next identify a neighborhood on which
both the population and empirical effective
covariances have uniformly separated leading
spectral clusters.
For $\lambda>0$ satisfying
Assumption~\ref{ass:fair-uniform-optimal-gap},
recall the constants $\delta_\lambda$
and $\sigma$ defined in Section~
\ref{subsec:fair-sharp-rate}.
Define
\begin{equation}
\rho_\lambda
=
\frac{
\delta_\lambda
}{
2\lambda\sigma\sqrt G
}
\label{eq:supp:fair-gap-neighborhood-radius}
\end{equation}
and
\begin{equation}
\mathcal N_\lambda
=
\left\{
\alpha\in\Delta_G:
\operatorname{dist}_2
(\alpha,\mathcal A_\lambda^\star)
\leq\rho_\lambda
\right\}.
\label{eq:supp:fair-gap-neighborhood}
\end{equation}

\begin{lemma}[Uniform local eigengap]
\label{lem:supp:fair-local-eigengap}
Suppose $\lambda>0$ and
Assumption~\ref{ass:fair-uniform-optimal-gap}
holds. Then
for every $\alpha\in\mathcal N_\lambda$,
\[
\operatorname{gap}_r\{M_\lambda(\alpha)\}
\geq\delta_\lambda.
\]
Furthermore, on the event
\[
(1+\lambda)e_n\leq\frac{\delta_\lambda}{4},
\]
we have
\[
\operatorname{gap}_r
\{\widehat M_\lambda(\alpha)\}
\geq\frac{\delta_\lambda}{2}
\]
uniformly over
$\alpha\in\mathcal N_\lambda$, where \(e_n\) is defined in \eqref{eq:fair-covariance-error}.
\end{lemma}

\begin{proof}
Fix $\alpha\in\mathcal N_\lambda$
and let
\[
\alpha^\dagger
=
\Pi_{\mathcal A_\lambda^\star}(\alpha).
\]
Since
\[
M_\lambda(\alpha)
-
M_\lambda(\alpha^\dagger)
=
\lambda
\sum_{a=1}^G
(\alpha_a-\alpha_a^\dagger)\Sigma_a,
\]
we obtain
\[
\begin{aligned}
\|M_\lambda(\alpha)
-M_\lambda(\alpha^\dagger)\|_{\mathrm{op}}
\leq
\lambda\sigma
\|\alpha-\alpha^\dagger\|_1
\leq
\lambda\sigma\sqrt G
\operatorname{dist}_2
(\alpha,\mathcal A_\lambda^\star)
\leq
\lambda\sigma\sqrt G\,\rho_\lambda
=
\frac{\delta_\lambda}{2}.
\end{aligned}
\]
By Assumption~\ref{ass:fair-uniform-optimal-gap},
\[
\operatorname{gap}_r
\{M_\lambda(\alpha^\dagger)\}
\geq2\delta_\lambda.
\]
Weyl's inequality therefore gives
\[
\begin{aligned}
\operatorname{gap}_r
\{M_\lambda(\alpha)\}
\geq
2\delta_\lambda
-
2\|M_\lambda(\alpha)
-M_\lambda(\alpha^\dagger)\|_{\mathrm{op}}
\geq\delta_\lambda.
\end{aligned}
\]
Similarly,
\[
\|\widehat M_\lambda(\alpha)
-M_\lambda(\alpha)\|_{\mathrm{op}}
\leq(1+\lambda)e_n.
\]
Another application of Weyl's inequality
gives
\[
\begin{aligned}
\operatorname{gap}_r
\{\widehat M_\lambda(\alpha)\}
\geq
\delta_\lambda-2(1+\lambda)e_n
\geq\frac{\delta_\lambda}{2}
\end{aligned}
\]
on the stated event.
Both bounds are uniform over
$\mathcal N_\lambda$.
\end{proof}

The following lemma controls the empirical
dual gradient in the eigengap neighborhood.

\begin{lemma}[Uniform dual-gradient perturbation]
\label{lem:supp:fair-dual-gradient}
Suppose $\lambda>0$ and
Assumption~\ref{ass:fair-uniform-optimal-gap}
holds.
On the event
\[
(1+\lambda)e_n\leq\frac{\delta_\lambda}{4},
\]
both $\phi_\lambda$ and
$\widehat\phi_\lambda$ are differentiable
throughout $\mathcal N_\lambda$.
Moreover, there exists a universal constant
$C_{\mathrm{sp}}>0$ such that
\begin{equation}
\begin{aligned}
\sup_{\alpha\in\mathcal N_\lambda}
\|\nabla\widehat\phi_\lambda(\alpha)
-\nabla\phi_\lambda(\alpha)\|_2
\leq
\lambda r\sqrt G
\left\{
2+
C_{\mathrm{sp}}
\frac{(1+\lambda)\sigma}{\delta_\lambda}
\right\}e_n.
\end{aligned}
\label{eq:supp:fair-dual-gradient-bound}
\end{equation}
\end{lemma}

\begin{proof}
For $\alpha\in\mathcal N_\lambda$,
write
\[
P(\alpha)=P_r\{M_\lambda(\alpha)\},
\qquad
\widehat P(\alpha)
=
P_r\{\widehat M_\lambda(\alpha)\}.
\]
Lemma~\ref{lem:supp:fair-local-eigengap}
ensures that both projectors are
uniquely defined on the stated event.
By differentiability of the Ky Fan
functional at matrices with a strict
rank-$r$ eigengap,
\[
\frac{\partial\phi_\lambda(\alpha)}
{\partial\alpha_a}
=
\lambda
\{c_a-\operatorname{tr}(P(\alpha)\Sigma_a)\}
\]
and
\[
\frac{\partial\widehat\phi_\lambda(\alpha)}
{\partial\alpha_a}
=
\lambda
\{\widehat c_a
-\operatorname{tr}(\widehat P(\alpha)
\widehat\Sigma_a)\}.
\]
For every group $a$,
\[
\begin{aligned}
\frac1\lambda
\left|
\frac{\partial\widehat\phi_\lambda(\alpha)}
{\partial\alpha_a}
-
\frac{\partial\phi_\lambda(\alpha)}
{\partial\alpha_a}
\right|
\leq
|\widehat c_a-c_a|
+
\left|
\operatorname{tr}
\{\widehat P(\alpha)
(\widehat\Sigma_a-\Sigma_a)\}
\right|
+
\left|
\operatorname{tr}
\{(\widehat P(\alpha)-P(\alpha))
\Sigma_a\}
\right|.
\end{aligned}
\]
By Weyl's inequality,
\[
|\widehat c_a-c_a|
\leq
r\|\widehat\Sigma_a-\Sigma_a\|_{\mathrm{op}}
\leq re_n.
\]
Similarly, because
$\widehat P(\alpha)$ is a
rank-$r$ orthogonal projector,
\[
\left|
\operatorname{tr}
\{\widehat P(\alpha)
(\widehat\Sigma_a-\Sigma_a)\}
\right|
\leq re_n.
\]
The population eigengap on
$\mathcal N_\lambda$ is at least
$\delta_\lambda$, and
\[
\|\widehat M_\lambda(\alpha)
-M_\lambda(\alpha)\|_{\mathrm{op}}
\leq(1+\lambda)e_n.
\]
The Davis--Kahan spectral projector
perturbation bound therefore gives
\[
\|\widehat P(\alpha)-P(\alpha)\|_F
\leq
C_{\mathrm{DK}}
\frac{
\sqrt r(1+\lambda)e_n
}{
\delta_\lambda
}
\]
for a universal constant
$C_{\mathrm{DK}}>0$.
Since the difference of two rank-$r$
projectors has rank at most $2r$,
\[
\begin{aligned}
\|\widehat P(\alpha)-P(\alpha)\|_*
\leq
\sqrt{2r}
\|\widehat P(\alpha)-P(\alpha)\|_F
\leq
C_{\mathrm{sp}}
\frac{
r(1+\lambda)e_n
}{
\delta_\lambda
}
\end{aligned}
\]
for another universal constant
$C_{\mathrm{sp}}>0$.
Consequently,
\[
\begin{aligned}
\left|
\operatorname{tr}
\{(\widehat P(\alpha)-P(\alpha))
\Sigma_a\}
\right|
\leq
\|\widehat P(\alpha)-P(\alpha)\|_*
\|\Sigma_a\|_{\mathrm{op}}
\leq
C_{\mathrm{sp}}
\frac{
r(1+\lambda)\sigma
}{
\delta_\lambda
}e_n.
\end{aligned}
\]
Combining the preceding bounds,
\[
\begin{aligned}
\left|
\frac{\partial\widehat\phi_\lambda(\alpha)}
{\partial\alpha_a}
-
\frac{\partial\phi_\lambda(\alpha)}
{\partial\alpha_a}
\right|
\leq
\lambda r
\left\{
2+
C_{\mathrm{sp}}
\frac{(1+\lambda)\sigma}{\delta_\lambda}
\right\}e_n.
\end{aligned}
\]
Taking the Euclidean norm over
$a=1,\ldots,G$ and then the supremum
over $\mathcal N_\lambda$ proves
\eqref{eq:supp:fair-dual-gradient-bound}.
\end{proof}

We next strengthen the preliminary localization
bound to a first-order bound.

\begin{lemma}[First-order dual-set perturbation]
\label{lem:supp:fair-first-order-dual}
Suppose $\lambda>0$ and
Assumption~\ref{ass:fair-uniform-optimal-gap}
holds. Let
$\widehat\alpha_\lambda
\in\widehat{\mathcal A}_\lambda$
belong to $\mathcal N_\lambda$,
and suppose
\[
(1+\lambda)e_n\leq\frac{\delta_\lambda}{4}.
\]
Then
\begin{equation}
\begin{aligned}
&
\operatorname{dist}_2
(\widehat\alpha_\lambda,
\mathcal A_\lambda^\star)
\leq
\frac{
2\lambda r\sqrt G
}{
\underline\mu_\lambda
}
\left\{
2+
C_{\mathrm{sp}}
\frac{(1+\lambda)\sigma}{\delta_\lambda}
\right\}e_n.
\end{aligned}
\label{eq:supp:fair-first-order-dual}
\end{equation}
\end{lemma}

\begin{proof}
For brevity, write
\[
\widehat\alpha=\widehat\alpha_\lambda,
\qquad
\alpha^\dagger
=
\Pi_{\mathcal A_\lambda^\star}
(\widehat\alpha),
\]
and
\[
d
=
\|\widehat\alpha-\alpha^\dagger\|_2.
\]
By the global quadratic-growth inequality
\eqref{eq:fair-derived-dual-growth},
\[
\frac{\underline\mu_\lambda}{2}d^2
\leq
\phi_\lambda(\alpha^\dagger)
-\phi_\lambda(\widehat\alpha).
\]
Since $\phi_\lambda$ is concave and
differentiable at $\widehat\alpha$,
\[
\begin{aligned}
\phi_\lambda(\alpha^\dagger)
-\phi_\lambda(\widehat\alpha)
\leq
\left\langle
\nabla\phi_\lambda(\widehat\alpha),
\alpha^\dagger-\widehat\alpha
\right\rangle.
\end{aligned}
\]
On the stated event,
$\widehat\phi_\lambda$ is also
differentiable at $\widehat\alpha$.
Because $\widehat\alpha$ maximizes
the concave function
$\widehat\phi_\lambda$ over $\Delta_G$,
the first-order variational inequality gives
\[
\left\langle
\nabla\widehat\phi_\lambda
(\widehat\alpha),
\alpha^\dagger-\widehat\alpha
\right\rangle
\leq0.
\]
Therefore,
\[
\begin{aligned}
\frac{\underline\mu_\lambda}{2}d^2
\leq
\left\langle
\nabla\phi_\lambda(\widehat\alpha)
-\nabla\widehat\phi_\lambda
(\widehat\alpha),
\alpha^\dagger-\widehat\alpha
\right\rangle
\leq
\|\nabla\phi_\lambda(\widehat\alpha)
-\nabla\widehat\phi_\lambda
(\widehat\alpha)\|_2\,d.
\end{aligned}
\]
If $d=0$, the desired inequality is
immediate. Otherwise, dividing by $d$
gives
\[
d
\leq
\frac{2}{\underline\mu_\lambda}
\|\nabla\phi_\lambda(\widehat\alpha)
-\nabla\widehat\phi_\lambda
(\widehat\alpha)\|_2.
\]
Applying
Lemma~\ref{lem:supp:fair-dual-gradient}
completes the proof.
\end{proof}

We are now ready to prove
Theorem~\ref{thm:fair-sharp-projector-rate}.
The first step is to show that the
deterministic small-perturbation condition
ensures localization of every empirical
dual optimizer within $\mathcal N_\lambda$
and preservation of the required
empirical eigengap.
\begin{proof}[Proof of Theorem~
\ref{thm:fair-sharp-projector-rate}]
Let
\[
\mathfrak c_\lambda
=
(1+\lambda)
+
\frac{\lambda^2rG\sigma}
{\underline\mu_\lambda}
\left\{
1+\frac{(1+\lambda)\sigma}
{\delta_\lambda}
\right\}.
\]
Recall the eigengap neighborhood
$\mathcal N_\lambda$ and its radius
\[
\rho_\lambda
=
\frac{\delta_\lambda}
{2\lambda\sigma\sqrt G}
\]
from
\eqref{eq:supp:fair-gap-neighborhood-radius}.

We first show that the small-perturbation
condition \eqref{eq:fair-deterministic-sharp-condition} ensures that every empirical
dual optimizer belongs to
$\mathcal N_\lambda$.
By Lemma~\ref{lem:supp:fair-dual-localization},
every
$\widehat\alpha_\lambda
\in\widehat{\mathcal A}_\lambda$
satisfies
\[
\operatorname{dist}_2^2
(\widehat\alpha_\lambda,
\mathcal A_\lambda^\star)
\leq
\frac{4r(1+2\lambda)e_n}
{\underline\mu_\lambda}.
\]
Consequently,
\[
\begin{aligned}
\frac{
\operatorname{dist}_2^2
(\widehat\alpha_\lambda,
\mathcal A_\lambda^\star)
}{
\rho_\lambda^2
}
\leq
\frac{
16\lambda^2rG\sigma^2
(1+2\lambda)e_n
}{
\underline\mu_\lambda
\delta_\lambda^2
}
\leq
32\frac{\mathfrak c_\lambda e_n}
{\delta_\lambda},
\end{aligned}
\]
where the last inequality follows from
$1+2\lambda\leq2(1+\lambda)$ and
the definition of $\mathfrak c_\lambda$.
Thus, whenever
\[
\mathfrak c_\lambda e_n
\leq\frac{\delta_\lambda}{64},
\]
all empirical dual optimizers belong
to $\mathcal N_\lambda$.
Moreover, since
$\mathfrak c_\lambda\geq1+\lambda$,
the same condition implies
\[
(1+\lambda)e_n
\leq\frac{\delta_\lambda}{64}
<
\frac{\delta_\lambda}{4}.
\]
Lemmas~\ref{lem:supp:fair-local-eigengap}
and~\ref{lem:supp:fair-dual-gradient}
therefore apply throughout the
relevant neighborhood. Since every empirical dual optimizer
belongs to $\mathcal N_\lambda$,
the conditions of
Lemma~\ref{lem:supp:fair-first-order-dual}
are also satisfied.
By Lemma~\ref{lem:supp:fair-first-order-dual},
\[
\begin{aligned}
\operatorname{dist}_2
(\widehat\alpha_\lambda,
\mathcal A_\lambda^\star)
\leq
\frac{2\lambda r\sqrt G}
{\underline\mu_\lambda}
\left\{
2+
C_{\mathrm{sp}}
\frac{(1+\lambda)\sigma}
{\delta_\lambda}
\right\}e_n.
\end{aligned}
\]
Put $C_1=4+2C_{\mathrm{sp}}$.
Then
\begin{equation}
\begin{aligned}
\operatorname{dist}_2
(\widehat\alpha_\lambda,
\mathcal A_\lambda^\star)
\leq
C_1\frac{\lambda r\sqrt G}
{\underline\mu_\lambda}
\left\{
1+
\frac{(1+\lambda)\sigma}
{\delta_\lambda}
\right\}e_n.
\end{aligned}
\label{eq:supp:fair-deterministic-dual}
\end{equation}
Let
\[
\alpha^\dagger
=
\Pi_{\mathcal A_\lambda^\star}
(\widehat\alpha_\lambda).
\]
By the definitions of the population
and empirical effective covariances,
\[
\begin{aligned}
\|\widehat M_\lambda
(\widehat\alpha_\lambda)
-M_\lambda(\alpha^\dagger)\|_{\mathrm{op}}
\leq
(1+\lambda)e_n
+
\lambda\sigma\sqrt G\,
\operatorname{dist}_2
(\widehat\alpha_\lambda,
\mathcal A_\lambda^\star)
\leq
C_1\mathfrak c_\lambda e_n.
\end{aligned}
\]
Choose the universal constant
\[
c
=
\min\left\{
\frac1{64},
\frac1{2C_1}
\right\}.
\]
Under
\eqref{eq:fair-deterministic-sharp-condition},
the preceding arguments apply and
\[
\|\widehat M_\lambda
(\widehat\alpha_\lambda)
-M_\lambda(\alpha^\dagger)\|_{\mathrm{op}}
\leq
\frac{\delta_\lambda}{2}.
\]
Since
$\alpha^\dagger\in\mathcal A_\lambda^\star$,
Assumption~
\ref{ass:fair-uniform-optimal-gap}
gives
\[
\operatorname{gap}_r
\{M_\lambda(\alpha^\dagger)\}
\geq2\delta_\lambda.
\]
Weyl's inequality therefore yields
\[
\begin{aligned}
\operatorname{gap}_r
\{\widehat M_\lambda
(\widehat\alpha_\lambda)\}
\geq
2\delta_\lambda
-
2\|\widehat M_\lambda
(\widehat\alpha_\lambda)
-M_\lambda(\alpha^\dagger)\|_{\mathrm{op}}
\geq\delta_\lambda>0.
\end{aligned}
\]
By empirical primal--dual optimality,
the strict eigengap implies that the
empirical Fair PCA optimizer is unique
and satisfies
\[
\widehat H_\lambda
=
P_r\{\widehat M_\lambda
(\widehat\alpha_\lambda)\}.
\]
Furthermore, the population spectral
characterization gives
\[
P_\lambda^\star
=
P_r\{M_\lambda(\alpha^\dagger)\}.
\]
The Davis--Kahan theorem implies
\[
\begin{aligned}
\|\widehat H_\lambda-P_\lambda^\star\|_F
\leq
C_{\mathrm{DK}}
\frac{\sqrt r}{\delta_\lambda}
\|\widehat M_\lambda
(\widehat\alpha_\lambda)
-M_\lambda(\alpha^\dagger)\|_{\mathrm{op}}
\leq
C_{\mathrm{DK}}C_1
\frac{\sqrt r}{\delta_\lambda}
\mathfrak c_\lambda e_n.
\end{aligned}
\]
Taking
\[
C=C_1\max\{1,C_{\mathrm{DK}}\}
\]
establishes
\eqref{eq:fair-deterministic-dual-bound}
and
\eqref{eq:fair-deterministic-projector-bound}.
\end{proof}

\begin{proof}[Proof of Corollary~
\ref{cor:fair-sharp-projector-rate}]
Under Assumption~\ref{ass:fair-subgaussian},
Lemma~\ref{lem:supp:fair-simultaneous-concentration}
gives
\[
e_n=O_p(\eta_n).
\]
Since
\[
b_{\lambda,n}
=
\mathfrak c_\lambda\eta_n,
\]
the condition
\[
\frac{b_{\lambda,n}}{\delta_\lambda}
\longrightarrow0
\]
implies
\[
\frac{\mathfrak c_\lambda e_n}
{\delta_\lambda}
\xrightarrow{p}0.
\]
Consequently,
\[
\mathbb P\left\{
\mathfrak c_\lambda e_n
\leq c\delta_\lambda
\right\}
\longrightarrow1,
\]
where $c>0$ is the universal constant
in Theorem~
\ref{thm:fair-sharp-projector-rate}.
On this event, the deterministic theorem
gives
\[
\begin{aligned}
\operatorname{dist}_2
(\widehat\alpha_\lambda,
\mathcal A_\lambda^\star)
\leq
C\frac{\lambda r\sqrt G}
{\underline\mu_\lambda}
\left\{
1+\frac{(1+\lambda)\sigma}
{\delta_\lambda}
\right\}e_n.
\end{aligned}
\]
Using $e_n=O_p(\eta_n)$
and the definition of $s_{\lambda,n}$
proves
\eqref{eq:fair-sharp-dual-rate}.

The same theorem guarantees that
the empirical primal optimizer is
unique and is a rank-$r$ projector
with probability tending to one.
Moreover,
\[
\|\widehat H_\lambda-P_\lambda^\star\|_F
\leq
C\frac{\sqrt r}{\delta_\lambda}
\mathfrak c_\lambda e_n
\]
on the same event. It follows that
\[
\begin{aligned}
\|\widehat H_\lambda-P_\lambda^\star\|_F
=
O_p\left(
\frac{\sqrt r}{\delta_\lambda}
\mathfrak c_\lambda\eta_n
\right)
=
O_p\left(
\frac{\sqrt r}{\delta_\lambda}
b_{\lambda,n}
\right).
\end{aligned}
\]
This proves
\eqref{eq:fair-sharp-projector-rate}.
Finally, under
\eqref{eq:fair-sharp-consistency-condition},
the right-hand side converges to zero,
yielding Frobenius consistency.
\end{proof}
\subsection{An example separating the eigengap conditions}
\label{supp:subsec:fair-eigengap-separation}
We illustrate that
Assumption~\ref{ass:fair-optimal-eigengap}
does not imply
Assumption~\ref{ass:fair-uniform-optimal-gap}.
Let $p=2$, $r=1$, $G=3$, $\lambda=1$,
$q=(1/10,1/10,4/5)$, and
\[
\Sigma_1=
\begin{pmatrix}
5&2\\2&2
\end{pmatrix},
\qquad
\Sigma_2=
\begin{pmatrix}
5&-2\\-2&2
\end{pmatrix},
\qquad
\Sigma_3=
\begin{pmatrix}
5&0\\0&6
\end{pmatrix}.
\]
All three covariance matrices are positive
definite, with $c_1=c_2=c_3=6$, and
\[
\overline\Sigma
=
\operatorname{diag}(5,26/5).
\]
For $\alpha\in\Delta_3$, the effective
covariance is
\[
M_1(\alpha)
=
\begin{pmatrix}
10 & 2(\alpha_1-\alpha_2)\\
2(\alpha_1-\alpha_2)
&56/5-4(\alpha_1+\alpha_2)
\end{pmatrix}.
\]
Since its first diagonal entry is always
$10$, we have
$\mathcal K_1\{M_1(\alpha)\} =
\lambda_{\max}\{M_1(\alpha)\}\geq10$.
Moreover, this lower bound is attained
precisely when
\[
\alpha_1=\alpha_2=t,
\qquad
\frac{3}{20}\leq t\leq\frac12.
\]
Indeed, the off-diagonal entry must vanish,
and the second diagonal entry must not
exceed $10$.
Since
$\phi_1(\alpha)=6-\mathcal K_1\{M_1(\alpha)\}$,
the complete dual-optimal set is
\[
\mathcal A_1^\star
=
\left\{
(t,t,1-2t):
\frac{3}{20}\leq t\leq\frac12
\right\}.
\]
For every $\alpha(t)\in\mathcal A_1^\star$,
\[
M_1\{\alpha(t)\}
=
\operatorname{diag}
\left(10,\frac{56}{5}-8t\right),
\]
and therefore
\[
\operatorname{gap}_1
\{M_1(\alpha(t))\}
=
8t-\frac65.
\]
At $t=1/2$, the eigengap equals $14/5>0$,
so Assumption~\ref{ass:fair-optimal-eigengap}
holds and the population Fair PCA projector
is uniquely identified as
$P_1^\star=e_1e_1^\top$.
At $t=3/20$, however,
$M_1\{\alpha(t)\}=10I_2$,
whose eigengap is zero.
Consequently,
\[
\inf_{\alpha\in\mathcal A_1^\star}
\operatorname{gap}_1\{M_1(\alpha)\}
=0,
\]
and
Assumption~\ref{ass:fair-uniform-optimal-gap}
fails.

The example also illustrates why spectral ordering
matters in the characterization of the dual-optimal
set. At $P_1^\star=e_1e_1^\top$, all three groups
attain the same reconstruction regret, so the
active-group support condition imposes no restriction
on $\alpha\in\Delta_3$. Moreover, the cross-subspace
condition
\[
(I-P_1^\star)M_1(\alpha)P_1^\star=0
\]
is equivalent to $\alpha_1=\alpha_2=t$, where
$0\leq t\leq1/2$. For such weights,
\[
M_1\{\alpha(t)\}
=
\operatorname{diag}
\left(10,\frac{56}{5}-8t\right).
\]

Thus, the cross-subspace condition ensures that $\operatorname{ran}(P_1^\star)$ is an invariant subspace of the effective covariance, but does not guarantee that it is a leading eigenspace. Specifically, when $t<3/20$, the second eigenvalue exceeds the first, so $e_1$ is not a leading eigenvector and the corresponding weights are not dual optimal. At $t=3/20$, the two eigenvalues coincide, whereas for $t>3/20$, $e_1$ becomes the unique leading eigenvector. Consequently, without the stronger eigengap condition, the characterization in \eqref{eq:fair-polyhedral-dual-set} of Section~\ref{sec:fair-sample-theory} would incorrectly include weights with $t<3/20$ that are not dual optimal. This example highlights the distinction between subspace invariance and spectral optimality, illustrating why the cross-subspace condition alone is insufficient to characterize the dual-optimal set and clarifying the role of Assumption~\ref{ass:fair-uniform-optimal-gap}.

\section{Proofs for Section 5:
Uniform Path Recovery and Penalty Selection}
\label{supp:sec:fair-selection}
This section establishes uniform recovery
of the reconstruction--fairness path,
additional guarantees for adaptive
penalty selection, and the statistical
properties of fairness-tolerance
calibration.

As in Section~\ref{sec:fair-sample-theory},
the dimension, rank, number of groups,
group sample sizes, and covariance
matrices may depend on $n$.
The dependence is suppressed where
no confusion arises.

\subsection{Uniform reconstruction--fairness
path recovery}
\label{supp:subsec:fair-selection-path}

We first establish the deterministic
performance bounds used in
Theorem~\ref{thm:fair-selection-uniform-path}.
Recall that
\[
e_n
=
\max_{1\leq a\leq G}
\|\widehat\Sigma_a-\Sigma_a\|_{\mathrm{op}}
\]
is defined in
\eqref{eq:fair-covariance-error}.
For every $H\in\mathcal F_r$,
the uniform criterion arguments
of Section~\ref{sec:fair-sample-theory}
give
\begin{equation}
|\widehat B(H)-B(H)|
\leq2re_n.
\label{eq:supp:fair-selection-regret-error}
\end{equation}
We also require a corresponding
bound for relative reconstruction risk.

\begin{lemma}[Relative reconstruction
criterion bound]
\label{lem:supp:fair-selection-relative-risk}
For every $H\in\mathcal F_r$,
\begin{equation}
|\widehat C(H)-C(H)|
\leq2re_n.
\label{eq:supp:fair-selection-relative-error}
\end{equation}
\end{lemma}

\begin{proof}
Recall that
\[
C(H)
=
R(H)-R_0,
\qquad
\widehat C(H)
=
\widehat R(H)-\widehat R_0.
\]
Since
\[
R(H)
=
\operatorname{tr}(\overline\Sigma)
-\operatorname{tr}(H\overline\Sigma),
\]
we have
\[
C(H)
=
\mathcal K_r(\overline\Sigma)
-
\operatorname{tr}(H\overline\Sigma).
\]
Similarly,
\[
\widehat C(H)
=
\mathcal K_r(\widehat{\overline\Sigma})
-
\operatorname{tr}
(H\widehat{\overline\Sigma}).
\]
The Ky Fan Lipschitz inequality
\eqref{eq:supp:fair-ky-fan-lipschitz}
gives
\[
\left|
\mathcal K_r(\widehat{\overline\Sigma})
-
\mathcal K_r(\overline\Sigma)
\right|
\leq re_n.
\]
Moreover, since $H\in\mathcal F_r$,
\[
\begin{aligned}
\left|
\operatorname{tr}
\{H(\widehat{\overline\Sigma}
-\overline\Sigma)\}
\right|\leq
\|H\|_*
\|\widehat{\overline\Sigma}
-\overline\Sigma\|_{\mathrm{op}}
\leq re_n.
\end{aligned}
\]
Combining the last two inequalities
proves the claim.
\end{proof}

\begin{proof}[Proof of Theorem~
\ref{thm:fair-selection-uniform-path}]
Let $c,C>0$ be universal constants
for which Theorem~
\ref{thm:fair-sharp-projector-rate}
holds, increasing $C$ if necessary.
Suppose
\[
\mathfrak c_\Lambda e_n
\leq c\delta_\Lambda.
\]
For every $\lambda\in\Lambda$,
\[
\mathfrak c_\lambda
\leq\mathfrak c_\Lambda,
\qquad
\delta_\lambda
\geq\delta_\Lambda.
\]
Consequently,
\[
\mathfrak c_\lambda e_n
\leq c\delta_\lambda.
\]
The deterministic condition of
Theorem~\ref{thm:fair-sharp-projector-rate}
therefore holds simultaneously
for every $\lambda\in\Lambda$.
The theorem implies that each empirical
optimizer is unique and is a rank-$r$
orthogonal projector, and that
\[
\|\widehat H_\lambda-P_\lambda^\star\|_F
\leq
C\frac{\sqrt r\,\mathfrak c_\lambda e_n}
{\delta_\lambda}.
\]
Taking the supremum over $\Lambda$
gives
\eqref{eq:fair-selection-uniform-projector}.
We now control the two performance
coordinates.
For symmetric matrices $H,H'$,
\[
|R(H)-R(H')|
\leq
\|\overline\Sigma\|_{\mathrm{op}}
\|H-H'\|_*.
\]
Since
$\|\overline\Sigma\|_{\mathrm{op}}
\leq\sigma$,
\begin{equation}
|C(H)-C(H')|
\leq
\sigma\|H-H'\|_*.
\label{eq:supp:fair-selection-risk-lipschitz}
\end{equation}
Likewise, for each group,
\[
\begin{aligned}
|\mathcal E_a(H)-\mathcal E_a(H')|
=
|\operatorname{tr}
\{(H-H')\Sigma_a\}|\leq
\sigma\|H-H'\|_*.
\end{aligned}
\]
Using
\[
|\max_a x_a-\max_a y_a|
\leq\max_a|x_a-y_a|,
\]
we obtain
\begin{equation}
|B(H)-B(H')|
\leq
\sigma\|H-H'\|_*.
\label{eq:supp:fair-selection-regret-lipschitz}
\end{equation}
On the stated perturbation event,
both $\widehat H_\lambda$ and
$P_\lambda^\star$ are rank-$r$
orthogonal projectors.
Their difference has rank at most
$2r$, and hence
\[
\begin{aligned}
\|\widehat H_\lambda-P_\lambda^\star\|_*
\leq
\sqrt{2r}
\|\widehat H_\lambda-P_\lambda^\star\|_F\leq
C\frac{
r\mathfrak c_\Lambda e_n
}{
\delta_\Lambda
},
\end{aligned}
\]
where the universal constant $C$
has been enlarged if necessary.
By Lemma~
\ref{lem:supp:fair-selection-relative-risk},
\[
\begin{aligned}
|\widehat C(\widehat H_\lambda)
-C(P_\lambda^\star)|\leq
|\widehat C(\widehat H_\lambda)
-C(\widehat H_\lambda)|
+
|C(\widehat H_\lambda)
-C(P_\lambda^\star)|\leq
2re_n+
C\frac{
r\sigma\mathfrak c_\Lambda e_n
}{
\delta_\Lambda
}.
\end{aligned}
\]
Similarly, by
\eqref{eq:supp:fair-selection-regret-error}
and
\eqref{eq:supp:fair-selection-regret-lipschitz},
\[
\begin{aligned}
|\widehat B(\widehat H_\lambda)
-B(P_\lambda^\star)|\leq
2re_n+
C\frac{
r\sigma\mathfrak c_\Lambda e_n
}{
\delta_\Lambda
}.
\end{aligned}
\]
Both inequalities hold uniformly
over $\lambda\in\Lambda$.
This establishes
\eqref{eq:fair-selection-uniform-performance}.
Applying
\eqref{eq:supp:fair-selection-risk-lipschitz}
and
\eqref{eq:supp:fair-selection-regret-lipschitz}
without the empirical criterion errors
also gives
\eqref{eq:fair-selection-fitted-performance}.
Finally, Assumption~
\ref{ass:fair-subgaussian}
and the simultaneous covariance
concentration bound imply
\[
e_n=O_p(\eta_n).
\]
Condition
\eqref{eq:fair-selection-uniform-localization-rate}
ensures that
\[
\mathbb P\{
\mathfrak c_\Lambda e_n
\leq c\delta_\Lambda
\}
\longrightarrow1.
\]
Substitution into the deterministic
performance bound proves
\eqref{eq:fair-selection-path-rate}.
The final consistency assertion
follows from
\eqref{eq:fair-selection-path-consistency}.
\end{proof}

We next establish the regularity of
the population projector path, which
will be used in the analysis of
adaptive penalty selection.

\begin{lemma}[Lipschitz continuity
of the population projector path]
\label{lem:supp:fair-selection-path-lipschitz}
Suppose Assumption
\ref{ass:fair-uniform-optimal-gap}
holds for every $\lambda\in\Lambda$,
with $\delta_\Lambda>0$ as defined
in \eqref{eq:fair-selection-uniform-constants}.
Then, for all $\lambda,\mu\in\Lambda$,
\begin{equation}
\|P_\lambda^\star-P_\mu^\star\|_F
\leq
\frac{\sqrt{2r}\,\sigma}
{\delta_\Lambda}
|\lambda-\mu|.
\label{eq:supp:fair-selection-path-lipschitz}
\end{equation}
\end{lemma}

\begin{proof}
For every $\lambda\in\Lambda$,
the uniform optimal-set eigengap
assumption provides an optimal
effective covariance with rank-$r$
eigengap at least
$2\delta_\lambda\geq2\delta_\Lambda$.
By the population quadratic margin
\eqref{eq:supp:fair-population-margin},
\[
Q_\lambda(H)-Q_\lambda(P_\lambda^\star)
\geq
\delta_\Lambda
\|H-P_\lambda^\star\|_F^2,
\qquad H\in\mathcal F_r.
\]
Taking $H=P_\mu^\star$ gives
\begin{equation}
\begin{aligned}
\delta_\Lambda
\|P_\mu^\star-P_\lambda^\star\|_F^2\leq
Q_\lambda(P_\mu^\star)
-Q_\lambda(P_\lambda^\star).
\end{aligned}
\label{eq:supp:fair-selection-path-margin}
\end{equation}
Since $P_\mu^\star$ minimizes
$Q_\mu$,
\[
Q_\mu(P_\mu^\star)
\leq Q_\mu(P_\lambda^\star).
\]
Using
\[
Q_\lambda(H)
=
-\operatorname{tr}(H\overline\Sigma)
+\lambda B(H),
\]
we therefore obtain
\[
\begin{aligned}
Q_\lambda(P_\mu^\star)
-Q_\lambda(P_\lambda^\star)&=
Q_\mu(P_\mu^\star)
-Q_\mu(P_\lambda^\star)+
(\lambda-\mu)
\{B(P_\mu^\star)
-B(P_\lambda^\star)\}\\
&\leq
(\lambda-\mu)
\{B(P_\mu^\star)
-B(P_\lambda^\star)\}\\
&\leq
|\lambda-\mu|
|B(P_\mu^\star)
-B(P_\lambda^\star)|.
\end{aligned}
\]
By
\eqref{eq:supp:fair-selection-regret-lipschitz},
\[
\begin{aligned}
|B(P_\mu^\star)-B(P_\lambda^\star)|\leq
\sigma\|P_\mu^\star-P_\lambda^\star\|_*\leq
\sigma\sqrt{2r}
\|P_\mu^\star-P_\lambda^\star\|_F,
\end{aligned}
\]
where the last inequality follows
because both matrices are rank-$r$
projectors.
Combining the preceding inequalities
gives
\[
\begin{aligned}
\delta_\Lambda
\|P_\mu^\star-P_\lambda^\star\|_F^2\leq
\sigma\sqrt{2r}
|\lambda-\mu|
\|P_\mu^\star-P_\lambda^\star\|_F.
\end{aligned}
\]
If the projector difference vanishes,
the conclusion is immediate.
Otherwise, division by its
Frobenius norm proves
\eqref{eq:supp:fair-selection-path-lipschitz}.
\end{proof}
\subsection{Adaptive penalty selection
and subspace consistency}
\label{supp:subsec:fair-selection-preference}
\begin{proposition}[Additional oracle inequalities
and penalty consistency]
\label{prop:supp:fair-selection-penalty}
Suppose Assumption
\ref{ass:fair-uniform-optimal-gap}
holds for every $\lambda\in\Lambda$,
with $\delta_\Lambda>0$ and
$\mathfrak c_\Lambda<\infty$.
Let $\Phi$ be Lipschitz continuous
with constant $L_\Phi$ under the
$\ell_\infty$ norm.
On the event
\eqref{eq:fair-selection-uniform-condition},
the selected penalty satisfies
\begin{equation}
\begin{aligned}
0
\leq
D(\widehat\lambda)
-
\min_{\lambda\in\Lambda_n}D(\lambda)\leq
2L_\Phi\Delta_n,
\end{aligned}
\label{eq:fair-selection-oracle-inequality}
\end{equation}
and
\begin{equation}
\begin{aligned}
0
\leq
D(\widehat\lambda)
-
\min_{\lambda\in\Lambda}D(\lambda)\leq
2L_\Phi\Delta_n
+
\frac{2L_\Phi r\sigma^2}
{\delta_\Lambda}h_n.
\end{aligned}
\label{eq:fair-selection-continuum-oracle}
\end{equation}
Under Assumption~\ref{ass:fair-subgaussian}
and conditions
\eqref{eq:fair-selection-uniform-localization-rate}
and
\eqref{eq:fair-selection-path-consistency},
\[
D(\widehat\lambda)
-
\min_{\lambda\in\Lambda_n}D(\lambda)
\xrightarrow{p}0.
\]
If, in addition,
\eqref{eq:fair-selection-grid-rate}
holds, then
\begin{equation}
D(\widehat\lambda)
-
\min_{\lambda\in\Lambda}D(\lambda)
\xrightarrow{p}0.
\label{eq:fair-selection-continuum-consistency}
\end{equation}
Furthermore, under
Assumption~\ref{ass:fair-subgaussian}
and conditions
\eqref{eq:fair-selection-uniform-localization-rate}
and
\eqref{eq:fair-selection-path-consistency},
let
\[
\lambda_n^\star
\in
\arg\min_{\lambda\in\Lambda_n}D(\lambda).
\]
Suppose that, for every
$\varepsilon>0$, there exists
$c_\varepsilon>0$ such that,
for all sufficiently large $n$,
\begin{equation}
\inf_{\substack{\lambda\in\Lambda_n\\
|\lambda-\lambda_n^\star|
\geq\varepsilon}}
\{D(\lambda)-D(\lambda_n^\star)\}
\geq c_\varepsilon,
\label{eq:fair-selection-identification}
\end{equation}
whenever the set in the infimum
is nonempty. Then
\begin{equation}
|\widehat\lambda-\lambda_n^\star|
\xrightarrow{p}0.
\label{eq:fair-selection-penalty-consistency}
\end{equation}
\end{proposition}
\begin{proof}[Proof of Proposition
\ref{prop:supp:fair-selection-penalty}]
For every $\lambda\in\Lambda_n$,
the Lipschitz property of $\Phi$
gives
\[
\begin{aligned}
|\widehat D(\lambda)-D(\lambda)|
\leq
L_\Phi
\|\widehat{\mathbf f}(\lambda)
-\mathbf f(\lambda)\|_\infty\leq L_\Phi\Delta_n.
\end{aligned}
\]
Consequently,
\begin{equation}
\sup_{\lambda\in\Lambda_n}
|\widehat D(\lambda)-D(\lambda)|
\leq L_\Phi\Delta_n.
\label{eq:supp:fair-selection-preference-error}
\end{equation}
Let
\[
\lambda_n^\star
\in
\arg\min_{\lambda\in\Lambda_n}
D(\lambda).
\]
By the empirical optimality of
$\widehat\lambda$,
\[
\widehat D(\widehat\lambda)
\leq
\widehat D(\lambda_n^\star).
\]
Thus
\[
\begin{aligned}
&D(\widehat\lambda)
-D(\lambda_n^\star)=
D(\widehat\lambda)
-\widehat D(\widehat\lambda)+
\widehat D(\widehat\lambda)
-\widehat D(\lambda_n^\star)+
\widehat D(\lambda_n^\star)
-D(\lambda_n^\star)\leq2L_\Phi\Delta_n.
\end{aligned}
\]
The left-hand side is nonnegative
by population optimality, proving
\eqref{eq:fair-selection-oracle-inequality}.

We next compare the selected penalty
with the optimum over the entire
interval $\Lambda$.
By Lemma~
\ref{lem:supp:fair-selection-path-lipschitz}
and the Lipschitz properties of
$C$ and $B$,
\[
\begin{aligned}
\|\mathbf f(\lambda)
-\mathbf f(\mu)\|_\infty\leq
\sigma
\|P_\lambda^\star-P_\mu^\star\|_*\leq
\sigma\sqrt{2r}
\|P_\lambda^\star-P_\mu^\star\|_F\leq
\frac{2r\sigma^2}{\delta_\Lambda}
|\lambda-\mu|
\end{aligned}
\]
for all $\lambda,\mu\in\Lambda$.
Consequently,
\begin{equation}
|D(\lambda)-D(\mu)|
\leq
\frac{2L_\Phi r\sigma^2}
{\delta_\Lambda}
|\lambda-\mu|.
\label{eq:supp:fair-selection-preference-lipschitz}
\end{equation}
Since $\Lambda$ is compact and
$D$ is continuous, there exists
\[
\lambda^\circ
\in
\arg\min_{\lambda\in\Lambda}D(\lambda).
\]
By the definition of $h_n$,
there exists $\mu_n\in\Lambda_n$
such that
\[
|\mu_n-\lambda^\circ|
\leq h_n.
\]
It follows that
\[
\begin{aligned}
D(\widehat\lambda)-D(\lambda^\circ)=
D(\widehat\lambda)
-\min_{\lambda\in\Lambda_n}D(\lambda)+
\min_{\lambda\in\Lambda_n}D(\lambda)
-D(\lambda^\circ)&\leq
2L_\Phi\Delta_n
+D(\mu_n)-D(\lambda^\circ)\\
&\leq
2L_\Phi\Delta_n
+
\frac{2L_\Phi r\sigma^2}
{\delta_\Lambda}h_n.
\end{aligned}
\]
This proves
\eqref{eq:fair-selection-continuum-oracle}.

Under the stated concentration
and rate conditions,
Theorem~\ref{thm:fair-selection-uniform-path}
gives $\Delta_n=o_p(1)$.
The finite-grid oracle excess
therefore converges to zero
in probability. If, in addition,
\eqref{eq:fair-selection-grid-rate}
holds, the continuum oracle excess
also converges to zero, proving
\eqref{eq:fair-selection-continuum-consistency}.

It remains to prove penalty consistency.
Fix $\varepsilon>0$.
By
\eqref{eq:fair-selection-identification},
for all sufficiently large $n$,
\[
|\widehat\lambda-\lambda_n^\star|
\geq\varepsilon
\]
implies
\[
D(\widehat\lambda)
-D(\lambda_n^\star)
\geq c_\varepsilon.
\]
Combining this with the oracle inequality
on the event
\eqref{eq:fair-selection-uniform-condition},
we obtain
\[
\begin{aligned}
\mathbb P\{
|\widehat\lambda-\lambda_n^\star|
\geq\varepsilon
\}\leq
\mathbb P\{
\mathfrak c_\Lambda e_n
>c\delta_\Lambda
\}
+
\mathbb P\{
2L_\Phi\Delta_n
\geq c_\varepsilon
\}\longrightarrow0.
\end{aligned}
\]
This proves
\eqref{eq:fair-selection-penalty-consistency}.
\end{proof}
\begin{proof}[Proof of Corollary
\ref{cor:fair-selection-oracle}]
The uniform criterion bounds give
\[
\sup_{H\in\mathcal F_r}
|\widehat C(H)-C(H)|
\leq2re_n,
\qquad
\sup_{H\in\mathcal F_r}
|\widehat B(H)-B(H)|
\leq2re_n.
\]
By Lipschitz continuity of $\Phi$,
\[
\begin{aligned}
\left|
\Phi\{C(\widehat H_{\widehat\lambda}),
B(\widehat H_{\widehat\lambda})\}
-\widehat D(\widehat\lambda)
\right|\leq2L_\Phi re_n.
\end{aligned}
\]
By empirical optimality and
\eqref{eq:supp:fair-selection-preference-error},
\[
\begin{aligned}
\left|
\widehat D(\widehat\lambda)
-\min_{\lambda\in\Lambda_n}D(\lambda)
\right|=
\left|
\min_{\lambda\in\Lambda_n}\widehat D(\lambda)
-\min_{\lambda\in\Lambda_n}D(\lambda)
\right|\leq L_\Phi\Delta_n.
\end{aligned}
\]
Moreover,
\eqref{eq:supp:fair-selection-preference-lipschitz}
and the definition of $h_n$ imply
\[
\begin{aligned}
0
\leq
\min_{\lambda\in\Lambda_n}D(\lambda)
-\min_{\lambda\in\Lambda}D(\lambda)\leq
\frac{2L_\Phi r\sigma^2}
{\delta_\Lambda}h_n.
\end{aligned}
\]
Combining these three inequalities
establishes
\eqref{eq:fair-selection-oracle-performance}.

Under the stated stochastic conditions,
Theorem~\ref{thm:fair-selection-uniform-path}
gives $\Delta_n=o_p(1)$.
Furthermore, $e_n=O_p(\eta_n)$
and
\eqref{eq:fair-selection-path-consistency}
imply $re_n=o_p(1)$.
The small-perturbation event holds
with probability tending to one.
Together with
\eqref{eq:fair-selection-grid-rate},
these facts prove
\eqref{eq:fair-selection-oracle-consistency}.
\end{proof}

\begin{corollary}[Consistency of adaptively
selected subspaces]
\label{cor:fair-selection-subspace}
Suppose the population assumptions of
Theorem~\ref{thm:fair-selection-uniform-path}
hold, together with
Assumption~\ref{ass:fair-subgaussian}
and conditions
\eqref{eq:fair-selection-uniform-localization-rate}
and
\eqref{eq:fair-selection-path-consistency}.
Suppose further that
\begin{equation}
\frac{\sqrt r\,\mathfrak c_\Lambda\eta_n}
{\delta_\Lambda}
\longrightarrow0,
\qquad
\frac{\sqrt r\,\sigma}{\delta_\Lambda}
=O(1).
\label{eq:supp:fair-selection-subspace-conditions}
\end{equation}
If the grid identification condition
\eqref{eq:fair-selection-identification}
holds, then
\begin{equation}
\|\widehat H_{\widehat\lambda}
-P_{\lambda_n^\star}^\star\|_F
\xrightarrow{p}0.
\label{eq:fair-selection-subspace-consistency}
\end{equation}
Alternatively, let
\[
\lambda_\Lambda^\star
\in
\arg\min_{\lambda\in\Lambda}D(\lambda).
\]
Suppose that, for every
$\varepsilon>0$, there exists
$c_\varepsilon>0$ such that,
for all sufficiently large $n$,
\begin{equation}
\inf_{\substack{\lambda\in\Lambda\\
|\lambda-\lambda_\Lambda^\star|
\geq\varepsilon}}
\{D(\lambda)-D(\lambda_\Lambda^\star)\}
\geq c_\varepsilon,
\label{eq:fair-selection-continuum-identification}
\end{equation}
whenever the set in the infimum
is nonempty.
If, in addition,
\eqref{eq:fair-selection-grid-rate}
holds, then
\begin{equation}
\|\widehat H_{\widehat\lambda}
-P_{\lambda_\Lambda^\star}^\star\|_F
\xrightarrow{p}0.
\label{eq:fair-selection-continuum-subspace}
\end{equation}
\end{corollary}
\begin{proof}[Proof of Corollary
\ref{cor:fair-selection-subspace}]
By the triangle inequality,
\[
\begin{aligned}
\|\widehat H_{\widehat\lambda}
-P_{\lambda_n^\star}^\star\|_F\leq
\|\widehat H_{\widehat\lambda}
-P_{\widehat\lambda}^\star\|_F+
\|P_{\widehat\lambda}^\star
-P_{\lambda_n^\star}^\star\|_F.
\end{aligned}
\]
On the uniform perturbation event,
Theorem~\ref{thm:fair-selection-uniform-path}
gives
\[
\|\widehat H_{\widehat\lambda}
-P_{\widehat\lambda}^\star\|_F
\leq
C\frac{
\sqrt r\,\mathfrak c_\Lambda e_n
}{
\delta_\Lambda
}.
\]
Under the stated concentration
and rate assumptions, this term
is $o_p(1)$.
For the second term,
Lemma~\ref{lem:supp:fair-selection-path-lipschitz}
gives
\[
\begin{aligned}
\|P_{\widehat\lambda}^\star
-P_{\lambda_n^\star}^\star\|_F\leq
\frac{\sqrt{2r}\,\sigma}
{\delta_\Lambda}
|\widehat\lambda-\lambda_n^\star|.
\end{aligned}
\]
The multiplier
$\sqrt r\,\sigma/\delta_\Lambda$
is bounded by assumption, while
\eqref{eq:fair-selection-penalty-consistency}
gives
\[
|\widehat\lambda-\lambda_n^\star|
=o_p(1).
\]
Thus the second term is also
$o_p(1)$, completing the proof.
For the continuum assertion,
the oracle inequality
\eqref{eq:fair-selection-continuum-oracle}
and the stated rate conditions imply
\[
D(\widehat\lambda)
-D(\lambda_\Lambda^\star)
=o_p(1).
\]
By the identification condition
\eqref{eq:fair-selection-continuum-identification},
the same argument used to establish
\eqref{eq:fair-selection-penalty-consistency}
gives
\[
|\widehat\lambda-\lambda_\Lambda^\star|
\xrightarrow{p}0.
\]
By the triangle inequality and
Lemma~\ref{lem:supp:fair-selection-path-lipschitz},
\[
\begin{aligned}
\|\widehat H_{\widehat\lambda}
-P_{\lambda_\Lambda^\star}^\star\|_F\leq
\|\widehat H_{\widehat\lambda}
-P_{\widehat\lambda}^\star\|_F+
\frac{\sqrt{2r}\sigma}{\delta_\Lambda}
|\widehat\lambda-\lambda_\Lambda^\star|.
\end{aligned}
\]
The first term is $o_p(1)$
by the uniform projector bound, while
the second is $o_p(1)$ by the
identification condition and
$\sqrt r\,\sigma/\delta_\Lambda=O(1)$.
This proves
\eqref{eq:fair-selection-continuum-subspace}.
\end{proof}
\subsection{Penalty selection from a
prescribed fairness tolerance}
\label{supp:subsec:fair-selection-tolerance}

\begin{proof}[Proof of Proposition
\ref{prop:fair-selection-tolerance-calibration}]
Recall that
\[
\sup_{H\in\mathcal F_r}
|\widehat B(H)-B(H)|
\leq2re_n
\]
by the uniform regret bound
established in the proof of
Proposition
\ref{prop:fair-uniform-sample-criterion}.
Let
\[
\widehat B_{\min}
=
\min_{H\in\mathcal F_r}\widehat B(H).
\]
The uniform bound implies
\[
\widehat B_{\min}
\leq B_{\min}+2re_n<\tau.
\]
Thus the empirical fairness
constraint is strictly feasible.
Since $(r/p)I_p$ belongs to the
relative interior of $\mathcal F_r$,
a sufficiently small convex
combination of $(r/p)I_p$ and a
strictly feasible point remains
strictly feasible and belongs to
the relative interior of $\mathcal F_r$.
As $\widehat R$ is affine and
$\widehat B$ is convex and continuous,
Slater's condition holds.
Consequently, strong Lagrange
duality applies to
\[
\min_{\substack{H\in\mathcal F_r\\
\widehat B(H)\leq\tau}}
\widehat R(H).
\]
Consequently, its dual problem is
\[
\max_{\lambda\geq0}
\{\widehat V(\lambda)-\lambda\tau\}.
\]
To verify attainment, choose
$H_{\min}\in\mathcal F_r$
such that
$\widehat B(H_{\min})
=\widehat B_{\min}$.
Then
\[
\begin{aligned}
\widehat V(\lambda)-\lambda\tau
\leq
\widehat R(H_{\min})
+\lambda(\widehat B_{\min}-\tau)\longrightarrow-\infty
\end{aligned}
\]
as $\lambda\to\infty$.
Continuity therefore ensures
the existence of a finite
optimal multiplier
$\widehat\lambda(\tau)$.

Strong duality and complementary
slackness imply that, for any
optimal multiplier, there exists
an empirical penalized optimizer
$\widehat H_{\widehat\lambda(\tau)}$
satisfying
\[
\widehat B
(\widehat H_{\widehat\lambda(\tau)})
\leq\tau
\]
and
\[
\widehat\lambda(\tau)
\{\widehat B
(\widehat H_{\widehat\lambda(\tau)})
-\tau\}=0.
\]
This proves
\eqref{eq:fair-selection-tolerance-slackness}.
By the uniform regret bound,
\[
\begin{aligned}
B(\widehat H_{\widehat\lambda(\tau)})
&\leq
\widehat B
(\widehat H_{\widehat\lambda(\tau)})
+2re_n\leq\tau+2re_n.
\end{aligned}
\]
Taking the positive part proves
\eqref{eq:fair-selection-tolerance-violation}.
If $\widehat\lambda(\tau)>0$,
complementary slackness gives
\[
\widehat B
(\widehat H_{\widehat\lambda(\tau)})
=\tau.
\]
Hence
\[
\begin{aligned}
|B(\widehat H_{\widehat\lambda(\tau)})
-\tau|=
|B(\widehat H_{\widehat\lambda(\tau)})
-\widehat B
(\widehat H_{\widehat\lambda(\tau)})|
\leq2re_n,
\end{aligned}
\]
which establishes
\eqref{eq:fair-selection-tolerance-accuracy}.
\end{proof}

\section{Numerical illustrations of the theoretical results}
\label{supp:sec:fair-theory-numerics}
We provide three numerical investigations
of the theoretical results developed in
the main paper.
The first examines the local
reconstruction--fairness expansions around
ordinary PCA, including the distinction
between positive and vanishing local
fairness-response coefficients.
The second examines the statistical recovery
of the empirical Fair PCA projector as
group sample sizes increase.
The third investigates projector recovery
and dual-weight estimation when the
population dual-optimal set is nonunique.

\subsection{Local reconstruction--fairness expansions}
\label{supp:subsec:fair-theory-local}

We first illustrate
Theorem~\ref{thm:fair-local-expansion}
and Proposition~\ref{prop:fair-local-nondegeneracy}
using three fixed population covariance
configurations.
All calculations in this subsection use the
population covariance matrices; no sampling
is involved.
In Scenario A, let $p=2$, $r=1$, $G=2$,
$q=(1/2,1/2)$, and
\[
\Sigma_1^{(A)}
=
\begin{pmatrix}
2&1\\
1&2
\end{pmatrix},
\qquad
\Sigma_2^{(A)}
=
\begin{pmatrix}
4&-1\\
-1&1
\end{pmatrix}.
\]
The aggregate covariance is
$\overline\Sigma^{(A)}
=\operatorname{diag}(3,3/2)$.
Group 1 is the unique worst-regret group
under ordinary PCA, and
\[
\mathcal I_0^{(A)}=\{1\},
\qquad
\zeta_0^{(A)}=\frac23.
\]
In Scenario B, let $p=3$, $r=1$, $G=3$,
$q=(1/3,1/3,1/3)$, and
\[
\begin{aligned}
\Sigma_1^{(B)}=
\begin{pmatrix}
2&1&0\\
1&2&0\\
0&0&1
\end{pmatrix},
\Sigma_2^{(B)}
&=
\begin{pmatrix}
2&0&1\\
0&1&0\\
1&0&2
\end{pmatrix},
\Sigma_3^{(B)}
=
\begin{pmatrix}
12&-1&-1\\
-1&1&0\\
-1&0&1
\end{pmatrix}.
\end{aligned}
\]
Here
\[
\overline\Sigma^{(B)}
=
\operatorname{diag}(16/3,4/3,4/3).
\]
Groups 1 and 2 attain the same largest
regret under ordinary PCA, whereas group 3
has strictly smaller regret.
Thus $\mathcal I_0^{(B)}=\{1,2\}$.
Writing $\alpha=(t,1-t,0)$ gives
\[
\Gamma_0(\alpha)
=
\frac{t^2+(1-t)^2}{4},
\]
so that
\[
\zeta_0^{(B)}
=
\min_{0\leq t\leq1}
\Gamma_0(t,1-t,0)
=
\frac18.
\]
This configuration illustrates a positive
local coefficient with two distinct active
groups whose cross-subspace interactions
occur in different directions.
Finally, Scenario C has $p=2$, $r=1$,
$G=2$, $q=(1/2,1/2)$, and
\[
\Sigma_1^{(C)}
=
\begin{pmatrix}
3&1\\
1&1
\end{pmatrix},
\qquad
\Sigma_2^{(C)}
=
\begin{pmatrix}
3&-1\\
-1&1
\end{pmatrix}.
\]
The aggregate covariance is
$\overline\Sigma^{(C)}
=\operatorname{diag}(3,1)$.
Both groups are active under ordinary PCA,
but their cross-subspace interactions cancel
when $\alpha=(1/2,1/2)$.
Consequently,
\[
\mathcal I_0^{(C)}=\{1,2\},
\qquad
\zeta_0^{(C)}=0.
\]
This configuration illustrates the
locally constant population
regularization path.
For each configuration, we evaluate the
population Fair PCA projector
$P_\lambda^\star$ over
\[
\lambda\in
\{0.001,0.002,0.005,0.01,0.02,
0.05,0.075,0.1,0.15,0.2\}.
\]
The solutions are computed using the
effective-covariance spectral formulation
and checked by the primal--dual gap and
rank-one eigengap.
All reconstruction risks and group regrets
are then evaluated directly from the
population covariance matrices.

For the nondegenerate Scenarios A and B,
define
\[
\Delta B_\lambda
=
B(P_0^\star)-B(P_\lambda^\star),
\qquad
\Delta R_\lambda
=
R(P_\lambda^\star)-R(P_0^\star).
\]
Theorem~\ref{thm:fair-local-expansion}
and \eqref{eq:fair-local-quadratic-frontier}
predict that the three ratios
\[
\begin{aligned}
T_B(\lambda)
=
\frac{\Delta B_\lambda}
{2\zeta_0\lambda},
T_R(\lambda)
=
\frac{\Delta R_\lambda}
{\zeta_0\lambda^2},
T_F(\lambda)
=
\frac{4\zeta_0\Delta R_\lambda}
{(\Delta B_\lambda)^2}
\end{aligned}
\]
all converge to one as
$\lambda\downarrow0$.

\begin{figure}[t]
\centering
\includegraphics[width=0.95\textwidth]
{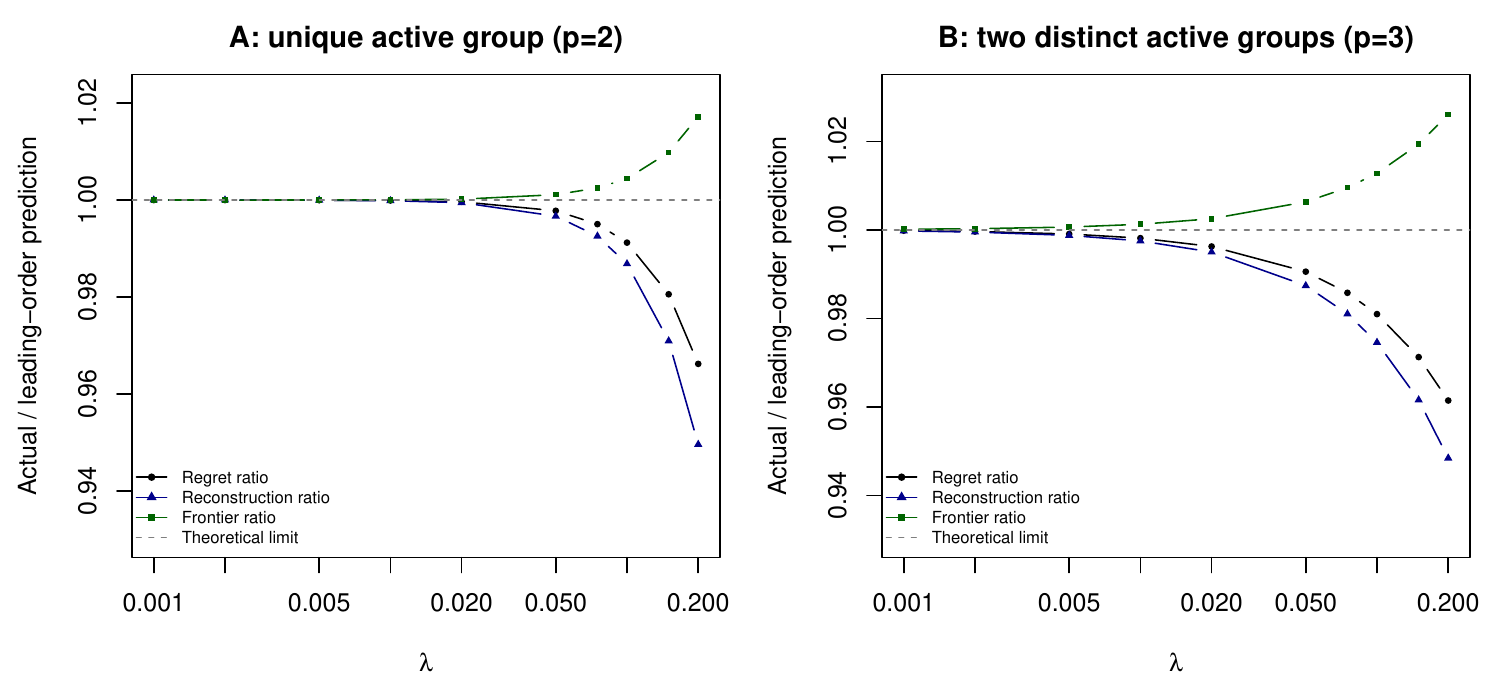}
\caption{
Numerical illustration of the local
reconstruction--fairness expansions.
The left panel corresponds to Scenario A
with one active group and $\zeta_0=2/3$;
the right panel corresponds to Scenario B
with two distinct active groups and
$\zeta_0=1/8$.
The curves show the ratios $T_B(\lambda)$,
$T_R(\lambda)$, and $T_F(\lambda)$.
The horizontal dashed line indicates
their common theoretical limit of one
as $\lambda\downarrow0$.
All quantities are evaluated using
population covariance matrices.
}
\label{fig:supp-fair-local-ratios}
\end{figure}

Figure~\ref{fig:supp-fair-local-ratios}
displays these ratios for Scenarios A and B.
The numerical results agree closely with
their theoretical limits.
At $\lambda=0.001$, the regret,
reconstruction, and frontier ratios are
approximately $(0.999999,0.999999,1.000000)$
in Scenario A and
$(0.999812,0.999750,1.000125)$
in Scenario B.
Even at $\lambda=0.1$, the corresponding
ratios are approximately
$(0.9912,0.9869,1.0044)$
and $(0.9810,0.9746,1.0128)$,
respectively.
Thus the leading-order approximations
remain accurate over a nontrivial range
of small penalties in these configurations.

\begin{figure}[t]
\centering
\includegraphics[width=0.78\textwidth]
{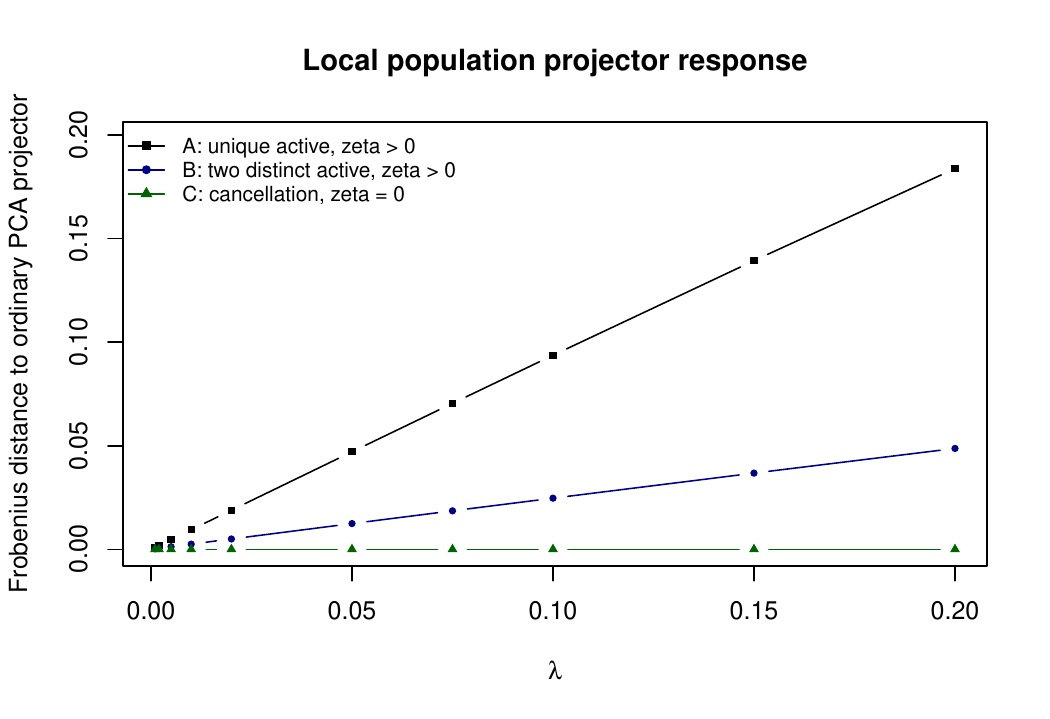}
\caption{
Distance between the population Fair PCA
projector and ordinary PCA as the fairness
penalty varies.
Scenarios A and B have $\zeta_0>0$
and exhibit an immediate change in the
population-optimal subspace.
Scenario C has $\zeta_0=0$ and its
population-optimal projector remains
unchanged over the displayed interval.
}
\label{fig:supp-fair-local-motion}
\end{figure}

We also examine the population projector
distance
$\|P_\lambda^\star-P_0^\star\|_F$.
As shown in
Figure~\ref{fig:supp-fair-local-motion},
this distance is positive for
small positive penalties in the
two nondegenerate scenarios.
In Scenario C, the projector distance,
regret improvement, and reconstruction
cost are all zero, to numerical precision,
throughout the evaluated penalty interval.
This illustrates the local-constancy
conclusion of
Proposition~\ref{prop:fair-local-nondegeneracy}.

All computed population solutions satisfy
the numerical primal--dual and eigengap
certification criteria.
Across the three scenarios and all
tested penalties, the maximum absolute
primal--dual gap is below $2.7\times10^{-15}$,
and the smallest effective-covariance
rank-one eigengap exceeds $1.5$.
These deterministic calculations illustrate
both the quantitative local expansions
and the qualitative distinction between
$\zeta_0>0$ and $\zeta_0=0$.

\subsection{Statistical recovery of the Fair PCA subspace}
\label{supp:subsec:fair-theory-projector}

We next investigate the convergence of the
empirical Fair PCA projector, focusing on
the first-order statistical recovery result
in Corollary~\ref{cor:fair-sharp-projector-rate}.

We consider $G=2$ groups with $p=20$,
$r=2$, aggregation weights
$q=(0.65,0.35)$, and a fixed fairness
penalty $\lambda=0.75$.
Let $e_j$ denote the $j$th canonical
basis vector and set
\[
u
=
\cos(0.55)e_2+\sin(0.55)e_3.
\]
The group covariance matrices are
\[
\begin{aligned}
\Sigma_1
&=
I_p+25e_1e_1^\top+8e_2e_2^\top,
\\
\Sigma_2
&=
I_p+25e_1e_1^\top+7uu^\top.
\end{aligned}
\]
Thus the groups share a strong first
principal direction but differ in the
orientation and strength of their
second signal component.

The population Fair PCA projector
$P_\lambda^\star$ is obtained from the
population dual spectral problem.
The numerical population dual optimizer
has first component approximately
$0.257972$, and the associated
unnormalized effective covariance
has rank-two eigengap approximately
$11.1645$.
The population primal--dual gap is
negligible.
The population dual optimizer is unique
in this configuration, so this experiment
examines first-order recovery in a
regular setting rather than the
nonunique-dual case.

For each group, we generate independent
observations
$X_{a,1},\ldots,X_{a,n}$
from $N_p(0,\Sigma_a)$, with
\[
n\in\{100,200,400,800,1600,3200\}.
\]
The two groups have the same sample size
$n$, and 300 independent Monte Carlo
replications are performed for each $n$.
Since the population means are known to
be zero, we use
\[
\widehat\Sigma_a
=
\frac1n\sum_{i=1}^n
X_{a,i}X_{a,i}^\top.
\]
For each replication, we compute the
empirical Fair PCA solution
$\widehat H_\lambda$ using the dual
spectral formulation.
The main performance measure is
the projector Frobenius error
\[
E_n
=
\|\widehat H_\lambda-P_\lambda^\star\|_F.
\]
We also calculate
\[
e_n
=
\max_{a=1,2}
\|\widehat\Sigma_a-\Sigma_a\|_{\mathrm{op}}
\]
to examine the relationship between
covariance estimation and subspace recovery.

Table~\ref{tab:supp-fair-projector-rate}
reports the mean and median projector
errors over the 300 replications.
The final column gives the median
root-$n$-scaled error.

\begin{table}[t]
\centering
\caption{
Statistical recovery of the Fair PCA projector.
For each sample size $n$ per group, entries
are based on 300 Monte Carlo replications.
Projector errors are measured by
$\|\widehat H_\lambda-P_\lambda^\star\|_F$.
}
\label{tab:supp-fair-projector-rate}
\small
\begin{tabular}{rccc}
\hline
$n$
&
Mean error
&
Median error
&
Median $\sqrt n E_n$
\\
\hline
100  & 0.1943 & 0.1950 & 1.9496 \\
200  & 0.1378 & 0.1365 & 1.9302 \\
400  & 0.0978 & 0.0969 & 1.9379 \\
800  & 0.0676 & 0.0670 & 1.8940 \\
1600 & 0.0481 & 0.0476 & 1.9046 \\
3200 & 0.0343 & 0.0341 & 1.9297 \\
\hline
\end{tabular}
\end{table}

The projector error decreases steadily
as the sample size increases, while the
median root-$n$-scaled errors remain
between approximately $1.89$ and $1.95$.
This is consistent with an $n^{-1/2}$
convergence rate in the present
fixed-dimensional Gaussian setting.

\begin{figure}[t]
\centering
\includegraphics[width=0.82\textwidth]
{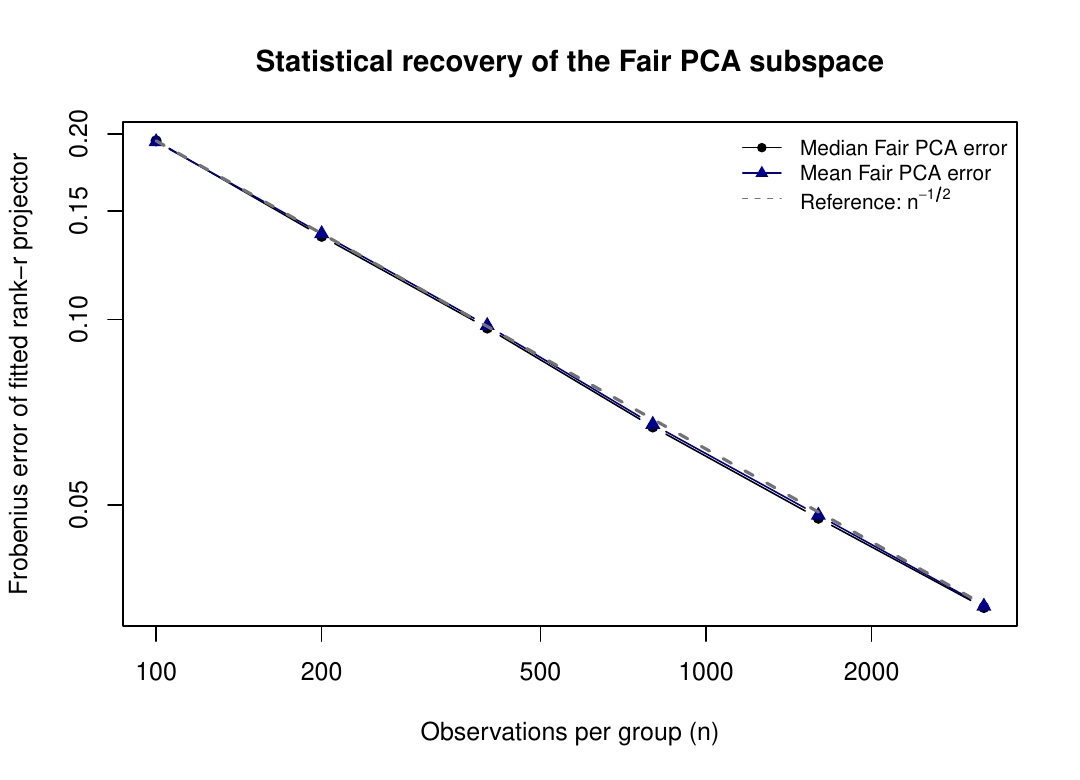}
\caption{
Monte Carlo convergence of the empirical
Fair PCA projector for $p=20$, $r=2$,
$G=2$, and $\lambda=0.75$.
The curves show the mean and median
Frobenius projector errors over
300 replications at each group sample size.
The dashed curve is proportional to
$n^{-1/2}$ and provides a reference
for the first-order statistical rate.
}
\label{fig:supp-fair-projector-rate}
\end{figure}

Figure~\ref{fig:supp-fair-projector-rate}
shows the mean and median projector
errors on logarithmic axes, together
with a reference curve proportional
to $n^{-1/2}$.
A least-squares fit of log mean error
against log sample size gives an
empirical slope of approximately $-0.503$.
The corresponding slope based on median
error is approximately $-0.505$.
For comparison, the mean covariance
estimation error has an empirical
log--log slope of approximately $-0.509$.

All 1,800 empirical fits passed the
numerical certification checks.
The largest absolute recorded
primal--dual gap was below
$2.4\times10^{-7}$.
The observed errors and convergence
slopes agree closely with the
first-order covariance-perturbation
behavior described by
Theorem~\ref{thm:fair-sharp-projector-rate}
and Corollary~\ref{cor:fair-sharp-projector-rate}.

\subsection{Statistical recovery with nonunique dual weights}
\label{supp:subsec:fair-theory-nonunique-dual}

We next investigate statistical recovery when the
population dual-optimal weights are nonunique.
This setting illustrates a distinctive feature of
Theorem~\ref{thm:fair-sharp-projector-rate}:
the empirical dual weights need only approach
the population dual-optimal set, rather than
any particular optimal weight vector, for
consistent recovery of the principal subspace.

Consider $G=3$ groups with $p=2$, $r=1$,
aggregation weights $q=(0.1,0.5,0.4)$, and
fairness penalty $\lambda=0.75$.
The group covariance matrices are
\[
\begin{aligned}
\Sigma_1
&=
\begin{pmatrix}
4&1\\
1&4
\end{pmatrix},
&
\Sigma_2
&=
\begin{pmatrix}
4&-1\\
-1&4
\end{pmatrix},
&
\Sigma_3
&=
\begin{pmatrix}
4&3/2\\
3/2&11/4
\end{pmatrix}.
\end{aligned}
\]
Each group has largest covariance eigenvalue
$c_a=\mathcal K_1(\Sigma_a)=5$.
The aggregate covariance is
\[
\overline\Sigma
=
\begin{pmatrix}
4&1/5\\
1/5&7/2
\end{pmatrix}.
\]
Let $e_1=(1,0)^\top$ and
$P_\lambda^\star=e_1e_1^\top$.
At this projector, all three groups have
reconstruction regret equal to one.
Direct evaluation of the population dual
problem gives the optimal set
\[
\mathcal A_\lambda^\star
=
\left\{
\alpha(t)=
\begin{pmatrix}
11/30-5t/4\\
19/30+t/4\\
t
\end{pmatrix}
:
0\leq t\leq22/75
\right\}.
\]
Thus the population dual-optimal set is a
nondegenerate line segment, rather than a
single weight vector.
For every $\alpha(t)\in\mathcal A_\lambda^\star$,
the effective covariance is
\[
M_\lambda\{\alpha(t)\}
=
\begin{pmatrix}
7&0\\
0&13/2-15t/16
\end{pmatrix}.
\]
Its rank-one eigengap is
\[
\frac12+\frac{15t}{16}
\geq\frac12.
\]
Consequently, all population-optimal dual
weights identify the same unique leading
rank-one projector $P_\lambda^\star$,
and the eigengap is uniformly positive
over the optimal set.
The ordinary PCA projector differs from
$P_\lambda^\star$, so this configuration
also represents a nontrivial fairness
adjustment.

For each group, we generate independent
observations from $N_2(0,\Sigma_a)$.
The group sample sizes are equal, with
\[
n\in\{200,400,800,1600,3200,6400\}.
\]
We perform 300 independent Monte Carlo
replications at each sample size.
Since the population means are known to
be zero, the group covariance estimators are
\[
\widehat\Sigma_a
=
\frac1n\sum_{i=1}^n
X_{a,i}X_{a,i}^\top.
\]
For each replication, we solve the empirical
dual problem over $\Delta_3$ and obtain
a candidate empirical Fair PCA projector
$\widehat H_\lambda$ and dual weight vector
$\widehat\alpha_\lambda$.
Numerical primal--dual and eigengap checks
are used to certify the rank-one solutions.

We examine the projector error
\[
E_n
=
\|\widehat H_\lambda-P_\lambda^\star\|_F
\]
and the distance from the empirical dual
weights to the population-optimal set,
\[
D_n
=
\operatorname{dist}_2
(\widehat\alpha_\lambda,
\mathcal A_\lambda^\star).
\]
To distinguish recovery of the optimal set
from identification of a particular vector,
we also consider
\[
F_n
=
\|\widehat\alpha_\lambda
-\alpha_{\mathrm{ref}}^\star\|_2,
\]
where $\alpha_{\mathrm{ref}}^\star$ is the
midpoint of $\mathcal A_\lambda^\star$,
corresponding to $t=11/75$.
The reference vector is chosen only for
comparison and is not singled out by the
population optimization problem.

Table~\ref{tab:supp-fair-nonunique-dual}
reports the mean errors and the number
of numerically certified solutions.
All 300 replications at each of
$n=1600,3200,6400$ pass the certification
checks.
At the three smaller sample sizes,
11, 19, and 10 replications, respectively,
do not yield certified rank-one solutions
under the implemented checks.
The reported error summaries therefore
use only certified replications, and the
small-sample results should be interpreted
with this qualification.

\begin{table}[t]
\centering
\caption{
Statistical recovery with nonunique
population dual weights.
Each sample size has 300 Monte Carlo
replications.
The columns report the number of
certified rank-one empirical solutions,
mean Frobenius projector error,
mean distance to the dual-optimal set,
and mean distance to its fixed midpoint.
Means are computed over certified
replications only.
}
\label{tab:supp-fair-nonunique-dual}
\small
\begin{tabular}{rcccc}
\hline
$n$ & Certified & Mean $E_n$
& Mean $D_n$ & Mean $F_n$\\
\hline
200  & 289/300 & 0.3059 & 0.1741 & 0.3214\\
400  & 281/300 & 0.1740 & 0.1128 & 0.2755\\
800  & 290/300 & 0.0816 & 0.0850 & 0.2601\\
1600 & 300/300 & 0.0445 & 0.0615 & 0.2463\\
3200 & 300/300 & 0.0320 & 0.0431 & 0.2449\\
6400 & 300/300 & 0.0233 & 0.0297 & 0.2399\\
\hline
\end{tabular}
\end{table}

The results show that both the projector
error and the distance to the population
dual-optimal set decrease as the sample
size increases.
Over the fully certified sample sizes
$n=1600,3200,6400$, the empirical
log--log slopes of their mean errors
are approximately $-0.466$ and $-0.525$,
respectively.
These values are consistent with
first-order statistical recovery in
the present fixed-dimensional setting.

\begin{figure}[t]
\centering
\includegraphics[width=0.82\textwidth]
{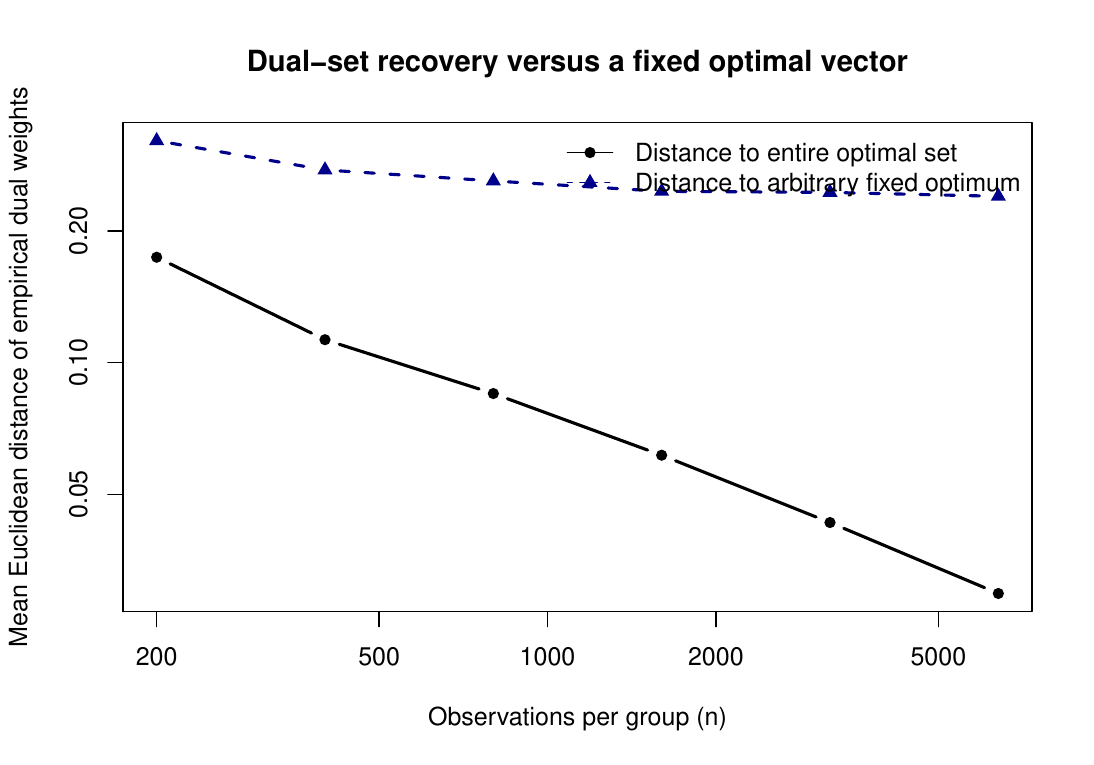}
\caption{
Recovery of the population dual-optimal
set versus a fixed optimal weight vector.
The solid curve reports the mean Euclidean
distance from the empirical dual weights
to $\mathcal A_\lambda^\star$.
The other curve reports the mean distance
to its midpoint
$\alpha_{\mathrm{ref}}^\star$.
Means are based on numerically certified
replications.
The distance to the optimal set decreases,
whereas the distance to the arbitrarily
chosen reference vector remains appreciable.
}
\label{fig:supp-fair-nonunique-distances}
\end{figure}

By contrast, the mean distance to the
fixed reference weight remains close
to $0.24$ at the larger sample sizes.
Figure~\ref{fig:supp-fair-nonunique-distances}
displays the different behavior of
distance to the optimal set and distance
to the chosen reference point.

\begin{figure}[t]
\centering
\includegraphics[width=0.82\textwidth]
{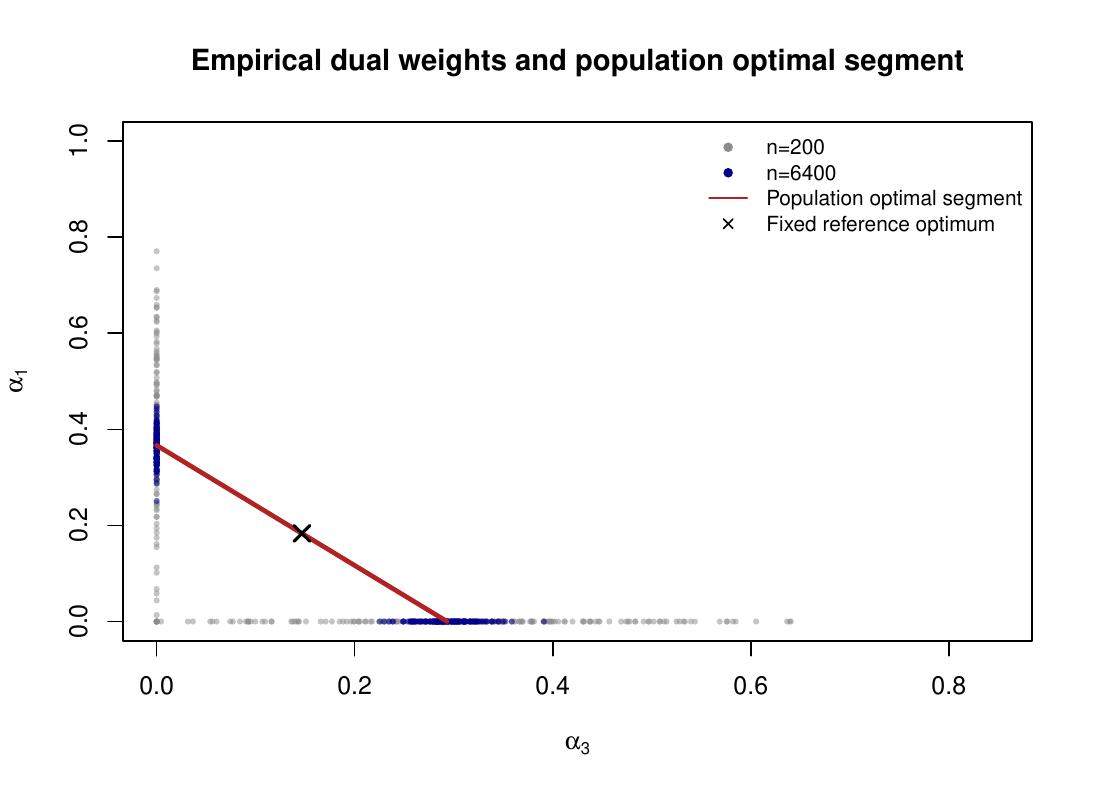}
\caption{
Empirical dual weights and the
population-optimal dual set.
The horizontal and vertical coordinates
are $\alpha_3$ and $\alpha_1$,
respectively.
Points correspond to numerically
certified empirical dual solutions
at $n=200$ and $n=6400$.
The line segment represents
$\mathcal A_\lambda^\star$,
and the cross marks its midpoint.
The empirical weights approach the
optimal segment without concentrating
around its midpoint.
}
\label{fig:supp-fair-nonunique-weights}
\end{figure}

Figure~\ref{fig:supp-fair-nonunique-weights}
further compares the empirical dual weights
with the population-optimal line segment.
As the sample size increases, the estimated
weights concentrate closer to the segment,
but may lie near different points along it.
At $n=6400$, approximately $48\%$ and
$49\%$ of the projected empirical weights
lie within the first and last $10\%$ of
the segment, respectively.
This pattern reflects the selection of
empirical dual optimizers by the numerical
procedure; it does not imply that every
optimizer-selection rule must favor the
endpoints.

The results illustrate why the statistical
analysis is formulated in terms of distance
to the dual-optimal set.
Although the empirical dual weights need
not approach a prespecified population-optimal
vector, the projector and dual-set errors
decrease as the covariance estimates improve.
This behavior is consistent with
Theorem~\ref{thm:fair-sharp-projector-rate}
and Corollary~\ref{cor:fair-sharp-projector-rate},
which permit nonunique population dual
optimizers while providing first-order
recovery guarantees under the stated
eigengap and regularity conditions.

\section{Finite-penalty and minimax-regret Fair PCA}
\label{supp:sec:numerical}

The simulation in the main paper
examines the empirical reconstruction--fairness
trade-off. We now investigate a complementary
finite-sample question: how moderate fairness
penalties compare with the two empirical
endpoints when reconstruction risk and
worst-group regret are evaluated under the
population covariances.
Since different aggregation weights are
considered, we write
\[
R_q(H)=\sum_{a=1}^2q_aR_a(H)
\]
and denote the corresponding empirical
Fair PCA estimator by
$\widehat H_{\lambda,q}$.
The empirical minimax-regret endpoint is
\[
\widehat H_\infty
\in\arg\min_{H\in\mathcal F_r}\widehat B(H),
\]
which does not depend on $q$.
We consider $G=2$ groups with $p=50$, $r=3$,
and covariance matrices
\[
\begin{aligned}
\Sigma_1
&=I_p+75e_1e_1^\top+35e_2e_2^\top
  +1.5e_3e_3^\top,\\
\Sigma_2
&=I_p+75e_1e_1^\top+35e_2e_2^\top
  +10uu^\top,
\end{aligned}
\]
where
\[
u=\cos(0.35)e_3+\sin(0.35)e_4
\]
and $e_j$ denotes the $j$th canonical basis
vector. Thus the groups share two strong
principal directions but differ in the
strength and orientation of their third
signal component.
Independent observations are generated from
$N_p(0,\Sigma_a)$, with
\[
(n_1,n_2)\in\{(135,15),(120,30)\}.
\]
For each configuration, we consider equal
weights $q=(0.5,0.5)$ and sample-proportion
weights $q_a=n_a/(n_1+n_2)$.
For each configuration, we consider equal
weights $q=(0.5,0.5)$ and sample-proportion
weights $q_a=n_a/(n_1+n_2)$.
The empirical estimators are computed using
their corresponding dual formulations, with
the minimax-regret endpoint solved directly
rather than approximated by a large finite
penalty.
All performance measures are evaluated
using the population covariance matrices.

The empirical covariance matrices are
computed using the known zero means.
We consider
\[
\lambda\in
\{0,0.025,0.05,0.1,0.2,0.5,0.8,1.5,\infty\}.
\]
Table~\ref{tab:supp-fair-finite-minimax}
reports the Monte Carlo mean reconstruction
risk and worst-group regret over 2,000
replications per configuration.
Within each Monte Carlo replication,
both weighting schemes are applied to
the same simulated dataset.
All reported solutions satisfy the numerical
optimality and eigengap certification criteria.
The empirical minimax-regret estimator
$\widehat H_\infty$ is independent of the
weights $q$, as is its population
worst-group regret $B(\widehat H_\infty)$.
Its aggregate reconstruction risk
$R_q(\widehat H_\infty)$ nevertheless varies
with $q$, since the latter determines the
relative contributions of the two groups
to the reconstruction criterion.
\begin{table}[t]
\centering
\caption{
Population reconstruction risk and worst-group
regret of empirical Fair PCA under different
penalties and weights.
Entries are Monte Carlo means over 2,000
replications. The cases $\lambda=0$ and
$\lambda=\infty$ correspond to ordinary
sample PCA and empirical minimax-regret
PCA, respectively.
}
\label{tab:supp-fair-finite-minimax}
\footnotesize
\begin{tabular}{cccrr}
\hline
$(n_1,n_2)$
&
$q_1$
&
$\lambda$
&
$R_q(\widehat H)$
&
$B(\widehat H)$
\\
\hline
$(135,15)$
& 0.9 & $0$      & 48.149 & 2.847 \\
&     & $0.025$  & 48.139 & 2.721 \\
&     & $0.05$   & 48.145 & 2.661 \\
&     & $0.1$    & 48.191 & 2.660 \\
&     & $0.2$    & 48.344 & 2.853 \\
&     & $0.5$    & 48.909 & 3.644 \\
&     & $0.8$    & 49.252 & 4.078 \\
&     & $1.5$    & 49.292 & 4.125 \\
&     & $\infty$ & 49.292 & 4.126 \\
\cline{2-5}
& 0.5 & $0$      & 50.303 & 4.348 \\
&     & $0.025$  & 50.257 & 4.296 \\
&     & $0.05$   & 50.216 & 4.250 \\
&     & $0.1$    & 50.158 & 4.185 \\
&     & $0.2$    & 50.113 & 4.133 \\
&     & $0.5$    & 50.107 & 4.125 \\
&     & $0.8$    & 50.107 & 4.126 \\
&     & $1.5$    & 50.107 & 4.126 \\
&     & $\infty$ & 50.107 & 4.126 \\
\hline
$(120,30)$
& 0.8 & $0$      & 48.143 & 2.078 \\
&     & $0.025$  & 48.136 & 2.028 \\
&     & $0.05$   & 48.135 & 1.998 \\
&     & $0.1$    & 48.148 & 1.974 \\
&     & $0.2$    & 48.207 & 2.007 \\
&     & $0.5$    & 48.382 & 2.194 \\
&     & $0.8$    & 48.389 & 2.201 \\
&     & $1.5$    & 48.389 & 2.201 \\
&     & $\infty$ & 48.389 & 2.201 \\
\cline{2-5}
& 0.5 & $0$      & 48.856 & 2.376 \\
&     & $0.025$  & 48.823 & 2.340 \\
&     & $0.05$   & 48.794 & 2.309 \\
&     & $0.1$    & 48.748 & 2.259 \\
&     & $0.2$    & 48.704 & 2.212 \\
&     & $0.5$    & 48.694 & 2.201 \\
&     & $0.8$    & 48.694 & 2.201 \\
&     & $1.5$    & 48.694 & 2.201 \\
&     & $\infty$ & 48.694 & 2.201 \\
\hline
\end{tabular}
\end{table}

Under sample-proportion weighting,
moderate fairness penalties achieve smaller
mean population worst-group regret than
both ordinary sample PCA and the empirical
minimax-regret estimator, sometimes together
with improved reconstruction risk.
For $(n_1,n_2)=(135,15)$, mean regret
decreases from $2.847$ at $\lambda=0$
to $2.721$ at $\lambda=0.025$, compared
with $4.126$ at the empirical minimax
endpoint, while mean reconstruction risk
decreases slightly from $48.149$ to
$48.139$. Moreover,
\[
B(\widehat H_{\lambda,q})
<
\min\{
B(\widehat H_{0,q}),
B(\widehat H_\infty)
\}
\]
holds in $95.2\%$ of replications at
$\lambda=0.025$, with
$R_q(\widehat H_{\lambda,q})
<R_q(\widehat H_{0,q})$
holding simultaneously in $69.8\%$.
For $(n_1,n_2)=(120,30)$, the
corresponding frequencies at
$\lambda=0.05$ are $82.2\%$ and
$51.2\%$, respectively.

The behavior differs under equal weighting.
Here, mean population worst-group regret
generally decreases towards that of the
empirical minimax endpoint as the penalty
increases. Nevertheless, even small
penalties improve both mean population
reconstruction risk and worst-group regret
relative to ordinary sample PCA.
Thus, intermediate penalties can provide
favorable population performance without
necessarily outperforming the empirical
minimax endpoint. Under sample-proportion
weighting, by contrast, mean population
regret is nonmonotone along the empirical
regularization path, attaining its smallest
values at moderate penalties.

The sensitivity to weights
reflects both population and estimation
effects. The weights determine the aggregate
covariance underlying the reconstruction
criterion, while also controlling the
contribution of each group's covariance
estimation error, since
\[
\widehat{\overline\Sigma}_q
-
\overline\Sigma_q
=
\sum_{a=1}^2
q_a(\widehat\Sigma_a-\Sigma_a).
\]
Under unequal sample sizes, equal weighting
assigns greater influence to the less
precisely estimated covariance of the
smaller group than sample-proportion weighting.
Consequently, the weighting schemes can
produce substantially different empirical
subspaces and regularization paths.
In contrast, the empirical minimax-regret
estimator minimizes $\widehat B(H)$ alone
and is independent of the prescribed
weights, although it remains
sensitive to covariance-estimation error.

These results illustrate the potential
finite-sample benefits of balancing
reconstruction and fairness rather than
optimizing either empirical criterion
alone. Ordinary sample PCA minimizes
empirical reconstruction risk, while the
empirical minimax-regret estimator minimizes
empirical worst-group regret; neither
property guarantees the best corresponding
population performance in finite samples.
A moderate penalty may improve population
performance relative to one or both
empirical endpoints, without being uniformly
preferable to either. This does not
contradict the optimality of population
ordinary PCA for reconstruction or of
population minimax-regret PCA for fairness,
but illustrates a possible statistical
benefit of moderate regularization when
the group covariances must be estimated.

\section{Additional analyses for NHANES accelerometry}
\label{supp:sec:fair-nhanes}

We provide additional analyses of the NHANES accelerometry
data considered in Section~\ref{sec:fair-nhanes}, examining
sensitivity to the projection rank, reconstruction weights,
and temporal resolution. Unless otherwise stated, we use
the same 6,055 participants, four age groups, and original
MIMS measurements as in the main analysis.

The projection-rank and temporal-resolution analyses
use five repeated 80:20 training--test splits,
corresponding to the first five partitions of the
primary experiment. The reconstruction-weight and
minimax-endpoint analyses use all 30 splits.
Partitions are stratified jointly by age group
and sex. All estimators are fitted on training
observations and evaluated on the corresponding
test samples. We report medians across splits of the
percentage reduction in worst-group regret
and the percentage change in aggregate
reconstruction risk, each calculated relative
to ordinary PCA within the same split.
Positive regret reductions indicate improved
fairness, whereas positive reconstruction costs
indicate reduced aggregate reconstruction accuracy.
The repeated splits assess sensitivity to the
random partition and are not treated as
independent replications for inference.

\subsection{Sensitivity to the projection rank}
\label{supp:subsec:fair-nhanes-rank}

We first examine ranks $r=1,2,3,5$, retaining $p=288$
and equal reconstruction weights $q_a=1/4$.
Table~\ref{tab:supp:fair-nhanes-rank} reports the
ordinary-PCA baselines and the test performance of
Fair PCA at two representative penalty levels.
\begin{table}[t]
\centering
\caption{
Sensitivity to the projection rank at $p=288$.
Explained variance refers to ordinary PCA on the
training sample. The column $B_0$ gives the median
ordinary-PCA test worst-group regret. Fairness gains
and reconstruction costs are median percentages
relative to ordinary PCA over five splits.
}
\label{tab:supp:fair-nhanes-rank}
\small
\begin{tabular}{crr rrrr}
\toprule
&
&
&
\multicolumn{2}{c}{$\lambda=0.1$}
&
\multicolumn{2}{c}{$\lambda=0.2$}
\\
\cmidrule(lr){4-5}
\cmidrule(lr){6-7}
$r$
& Var. (\%)
& $B_0$
& Gain (\%)
& Cost (\%)
& Gain (\%)
& Cost (\%)
\\
\midrule
1 & 30.6 & 243.09 & 19.00 & 0.013 & 22.61 & 0.041 \\
2 & 44.9 & 291.25 & 17.51 & 0.061 & 14.22 & 0.166 \\
3 & 50.6 & 414.73 & 11.21 & 0.158 &  9.99 & 0.307 \\
5 & 58.0 & 380.56 & 13.07 & 0.130 & 21.51 & 0.344 \\
\bottomrule
\end{tabular}
\end{table}
Positive median improvements in worst-group regret
are observed at every rank. At $\lambda=0.1$, test
regret decreases in four of the five splits for
$r=1$ and $r=3$, and in all five splits for
$r=2$ and $r=5$. The reconstruction costs remain
small, although their magnitudes vary with the rank.

The effect of the fairness penalty also depends on the rank.
For example, increasing $\lambda$ from $0.1$ to
$0.2$ improves the median test regret reduction
from 13.07\% to 21.51\% at $r=5$, accompanied
by an increase in reconstruction cost from
0.130\% to 0.344\%. In contrast, the corresponding
median regret reduction decreases at $r=2$.
Moreover, ordinary-PCA worst-group regret need
not decrease with $r$, since the group-specific
oracle benchmarks entering $B$ also depend
on the rank.

These results indicate that the improvement
observed in the main analysis is not specific
to $r=2$. They do not, however, identify a
uniformly preferable rank, as changing $r$
also changes the reconstruction objective and
the group-specific regret benchmarks.

\subsection{Sensitivity to reconstruction weights}
\label{supp:subsec:fair-nhanes-weights}

We examine the role of the reconstruction
weights $q_a$, fixing $p=288$ and $r=2$.
In addition to equal weights $q_a=1/4$,
we consider sample-proportion weights
\[
q_a=\frac{n_{a,\mathrm{tr}}}{n_{\mathrm{tr}}},
\qquad a=1,\ldots,4,
\]
where $n_{a,\mathrm{tr}}$ denotes the
training sample size of group $a$.
These weights are held fixed when evaluating
the fitted projectors on the test sample.
The proportions are approximately
$(0.287,0.254,0.295,0.164)$.

\begin{table}[t]
\centering
\caption{
Sensitivity to reconstruction weights at
$p=288$ and $r=2$. Entries are medians
over 30 train/test splits. Fairness gains
and reconstruction costs are percentage
changes relative to ordinary PCA fitted
under the same weighting scheme.
The final column gives the number of splits
with improved test worst-group regret.
}
\label{tab:supp:fair-nhanes-weights}
\small
\begin{tabular}{lcrrc}
\toprule
Weights
& $\lambda$
& Gain (\%)
& Cost (\%)
& Improved \\
\midrule
Equal
& 0.05 & 10.25 & 0.015 & 29/30 \\
& 0.10 & 18.54 & 0.053 & 29/30 \\
& 0.20 & 24.75 & 0.153 & 27/30 \\
& 0.40 & 24.89 & 0.180 & 28/30 \\
\addlinespace
Sample proportion
& 0.05 &  9.59 & 0.011 & 27/30 \\
& 0.10 & 16.62 & 0.041 & 25/30 \\
& 0.20 & 16.45 & 0.090 & 25/30 \\
& 0.40 & 18.34 & 0.109 & 24/30 \\
\bottomrule
\end{tabular}
\end{table}

Table~\ref{tab:supp:fair-nhanes-weights}
shows positive median reductions in test
worst-group regret under both weighting schemes,
with generally larger percentage improvements
under equal weighting. At $\lambda=0.1$,
the median reductions are 18.54\% and 16.62\%
under equal and sample-proportion weighting,
respectively. However, the ordinary-PCA
benchmarks differ: their median test
worst-group regrets are 291.30 and 265.72,
respectively. The corresponding absolute
regrets under Fair PCA are 238.97 and 225.52.
Thus, a larger percentage improvement does not
necessarily imply a smaller final worst-group regret.

This behavior is consistent with the
effective-covariance characterization.
Relative to equal weighting,
sample-proportion weighting assigns greater
weight to the 6--17 and 40--64 age groups
and less weight to the oldest group,
affecting both the ordinary-PCA benchmark
and the subsequent fairness adjustment.
Although the fairness penalty reweights
the group covariance matrices, the finite-penalty
solution generally remains dependent on the
initial reconstruction weights. Neither weighting
scheme is universally preferable, as their
relative performance depends on the group
covariance matrices.

\subsection{Sensitivity to temporal resolution}
\label{supp:subsec:fair-nhanes-resolution}

We also investigate the temporal resolution 
used to construct the daily activity profiles.
We divide the 24-hour period into intervals of
60, 30, 15, 10, or 5 minutes, producing dimensions
$p=24,48,96,144,288$, respectively.
All representations use the same 6,055 participants
and identical training--test partitions.
Each coordinate remains the mean valid MIMS
activity per minute within its clock-time
interval, averaged over eligible monitoring days.
We fix $r=2$, $\lambda=0.1$, and equal group weights.
Table~\ref{tab:supp:fair-nhanes-resolution}
summarizes the results.

\begin{table}[t]
\centering
\caption{
Sensitivity to temporal resolution at $r=2$ and
$\lambda=0.1$. The empirical effective-rank column
reports the maximum effective rank across the four
group covariance matrices, estimated from the
full matched cohort. Performance summaries are
medians over five identical train/test partitions.
}
\label{tab:supp:fair-nhanes-resolution}
\small
\begin{tabular}{rrrrrc}
\toprule
$p$
& Minutes
& Effective rank
& Gain (\%)
& Cost (\%)
& Improved \\
\midrule
 24 & 60 & 2.60 & 18.85 & 0.103 & 4/5 \\
 48 & 30 & 2.83 & 18.56 & 0.086 & 4/5 \\
 96 & 15 & 3.05 & 18.16 & 0.073 & 4/5 \\
144 & 10 & 3.17 & 17.97 & 0.068 & 4/5 \\
288 &  5 & 3.39 & 17.51 & 0.061 & 5/5 \\
\bottomrule
\end{tabular}
\end{table}

The median reduction in test worst-group regret
varies only modestly across resolutions, from
18.85\% at $p=24$ to 17.51\% at $p=288$.
When normalized by ordinary-PCA reconstruction
risk, the median reconstruction cost decreases
from 0.103\% to 0.061\%. However, as shown in
Section~\ref{supp:subsec:fair-nhanes-captured},
this pattern depends on the choice of
normalization. The favorable empirical
reconstruction--fairness trade-off nevertheless
persists across the examined resolutions.

One possible explanation is the concentrated
covariance structure of these activity profiles.
To quantify this, we compute the empirical
effective rank of each group covariance,
\[
\frac{
\operatorname{tr}(\widehat\Sigma_a)
}{
\|\widehat\Sigma_a\|_{\mathrm{op}}
}.
\]

The maximum empirical effective rank across
groups increases from approximately 2.60 to
3.39 as the ambient dimension increases from
24 to 288, indicating a concentration of
variance in relatively few spectral directions.
This provides a possible explanation for the
similar empirical performance across resolutions.

\subsection{Reconstruction cost relative to captured variance}
\label{supp:subsec:fair-nhanes-captured}
The preceding analyses normalize the increase
in test reconstruction risk by the test
reconstruction risk of training-fitted ordinary
PCA. We also consider normalization by the
maximum variance captured by a rank-$r$ projector.
The distinction is relevant when comparing
different projection ranks or dimensions,
for which reconstruction risk and captured
variance may have different relative magnitudes.
Let $\widehat H_0$ and $\widehat H_\lambda$ denote
the estimators fitted on the training sample, and define
\[
\Delta_{\mathrm{te}}(\lambda)
=
\widehat R_{\mathrm{te}}(\widehat H_\lambda)
-
\widehat R_{\mathrm{te}}(\widehat H_0).
\]
For each train/test split, we consider the two
percentage reconstruction costs
\[
100\,
\frac{\Delta_{\mathrm{te}}(\lambda)}
{\widehat R_{\mathrm{te}}(\widehat H_0)}
\qquad\text{and}\qquad
100\,
\frac{\Delta_{\mathrm{te}}(\lambda)}
{\mathcal K_r(\widehat{\bar\Sigma}_{\mathrm{te}})},
\]
where
\[
\widehat{\bar\Sigma}_{\mathrm{te}}
=
\sum_{a=1}^G q_a\widehat\Sigma_{\mathrm{te},a}.
\]
The first measure expresses the reconstruction cost
relative to the test reconstruction risk of
training-fitted ordinary PCA, whereas the second
expresses it relative to the maximum rank-$r$
variance captured on the test sample. The group
weights $q_a$ are the same as those used in fitting,
and test observations are used only for evaluation.

\begin{table}[t]
\centering
\caption{
Test reconstruction costs in the primary NHANES
experiment ($p=288$, $r=2$, equal group weights).
The two columns use the reconstruction risk of
training-fitted ordinary PCA and the maximum
rank-$r$ captured variance on the test sample,
respectively, as denominators. Entries are median
percentages over 30 train/test splits.
}
\label{tab:supp:fair-nhanes-captured-main}
\begin{tabular}{rcc}
\toprule
$\lambda$
& Cost / Reconstruction risk (\%)
& Cost / Captured variance (\%) \\
\midrule
0.025 & 0.004 & 0.005 \\
0.050 & 0.015 & 0.019 \\
0.100 & 0.053 & 0.065 \\
0.200 & 0.153 & 0.187 \\
0.400 & 0.180 & 0.217 \\
\bottomrule
\end{tabular}
\end{table}

Table~\ref{tab:supp:fair-nhanes-captured-main}
shows that the reconstruction costs remain small
under either normalization. At $\lambda=0.1$,
the median cost is 0.053\% relative to ordinary-PCA
reconstruction risk and 0.065\% relative to
captured variance. Thus, the favorable
reconstruction--fairness trade-off observed
in the primary experiment is preserved under
the alternative normalization.

The distinction becomes more pronounced when
comparing temporal resolutions.
Table~\ref{tab:supp:fair-nhanes-captured-resolution}
reports both measures for the matched-participant
experiments, with $r=2$, $\lambda=0.1$, and
equal group weights.

\begin{table}[t]
\centering
\caption{
Test reconstruction costs across temporal
resolutions under the two normalizations.
Entries are median percentages over five
matched train/test splits.
}
\label{tab:supp:fair-nhanes-captured-resolution}
\begin{tabular}{rcc}
\toprule
$p$
& Cost / Reconstruction risk (\%)
& Cost / Captured variance (\%) \\
\midrule
 24 & 0.103 & 0.067 \\
 48 & 0.086 & 0.071 \\
 96 & 0.073 & 0.073 \\
144 & 0.068 & 0.074 \\
288 & 0.061 & 0.076 \\
\bottomrule
\end{tabular}
\end{table}

When normalized by ordinary-PCA reconstruction risk,
the median cost decreases from 0.103\% at $p=24$
to 0.061\% at $p=288$. In contrast, the cost
normalized by captured variance remains relatively
stable, ranging from 0.067\% to 0.076\%.
This difference reflects changes in the relative
magnitudes of reconstruction risk and captured
variance across temporal resolutions. Consequently,
the decrease observed under reconstruction-risk
normalization does not imply a corresponding decrease
in the fraction of captured variance sacrificed
by Fair PCA.

We also examine the alternative normalization
in the rank and weighting analyses. At
$\lambda=0.1$, the median reconstruction costs
relative to captured variance are
0.030\%, 0.076\%, 0.155\%, and 0.094\%
for ranks $r=1,2,3,5$, respectively.
For equal and sample-proportion reconstruction
weights at $r=2$, the corresponding median
costs over 30 splits are 0.065\% and
0.051\%, respectively. The costs remain
small across these configurations, although
their relative magnitudes depend on the
projection rank and weighting scheme.

Overall, the favorable empirical
reconstruction--fairness trade-off persists
under both normalizations. Nevertheless,
comparisons of relative reconstruction costs
across ranks or dimensions require care,
as the choice of denominator can affect
their apparent magnitude and ordering.

\subsection{The minimax-regret endpoint}
\label{supp:subsec:fair-nhanes-minimax}

We further examine the limiting case
$\lambda=\infty$, corresponding to the
empirical minimax-regret problem
\[
\widehat H_\infty
\in
\arg\min_{H\in\mathcal F_r}
\widehat B_{\mathrm{tr}}(H).
\]
Unlike the finite-penalty objective, this
problem does not depend on the reconstruction
weights $q_a$. Consequently, for a given
training sample, the minimax solution set
is identical under equal and sample-proportion
weighting. In the present experiments, all
30 minimax fits satisfy the numerical
optimality and eigengap checks, yielding
rank-$r$ projectors.

Table~\ref{tab:supp:fair-nhanes-minimax}
compares the endpoint with representative
finite penalties for $p=288$ and $r=2$,
using the same 30 train/test splits.
Reconstruction costs are reported relative
to both ordinary-PCA reconstruction risk
and rank-$r$ captured variance, as defined
in Section~\ref{supp:subsec:fair-nhanes-captured}.

\begin{table}[t]
\centering
\caption{
Comparison of finite fairness penalties with
the minimax-regret endpoint for NHANES
accelerometry ($p=288$, $r=2$).
Entries are medians over 30 train/test splits.
Training and test gains are percentage reductions
in worst-group regret. Risk and variance costs
are percentage reconstruction costs normalized
by $\widehat R_{\mathrm{te}}(\widehat H_0)$
and $\mathcal K_r(\widehat{\bar\Sigma}_{\mathrm{te}})$,
respectively. The final column reports the number
of splits with improved test worst-group regret.
}
\label{tab:supp:fair-nhanes-minimax}
\small
\begin{tabular}{llrrrrc}
\toprule
Weights
& $\lambda$
& Train gain
& Test gain
& Risk cost (\%)
& Variance cost (\%)
& Improved \\
\midrule
Equal
& 0.1      & 22.72 & 18.54 & 0.053 & 0.065 & 29/30 \\
& 0.2      & 38.15 & 24.75 & 0.153 & 0.187 & 27/30 \\
& 0.4      & 39.82 & 24.89 & 0.180 & 0.217 & 28/30 \\
& $\infty$ & 41.52 & 23.25 & 0.243 & 0.298 & 26/30 \\
\addlinespace
Sample proportion
& 0.1      & 22.11 & 16.62 & 0.041 & 0.051 & 25/30 \\
& 0.2      & 31.51 & 16.45 & 0.090 & 0.113 & 25/30 \\
& 0.4      & 32.77 & 18.34 & 0.109 & 0.133 & 24/30 \\
& $\infty$ & 34.25 & 15.20 & 0.146 & 0.176 & 23/30 \\
\bottomrule
\end{tabular}
\end{table}

Under both weighting schemes, the minimax
endpoint achieves the largest median reduction
in training worst-group regret. However,
this does not translate into the largest
median test improvement. Under equal weighting,
the median test regret reduction is 23.25\%
at $\lambda=\infty$, compared with 24.89\%
at $\lambda=0.4$, while the corresponding
reconstruction cost increases. A similar
pattern is observed under sample-proportion
weighting.

Although the fitted minimax projector is the
same under both weighting schemes, the
percentage fairness improvements differ
because the ordinary-PCA benchmarks depend
on $q_a$. In particular, the absolute test
worst-group regret of the minimax projector
is identical, whereas the reconstruction
cost is evaluated using different aggregate
reconstruction risks.

The minimax endpoint also redistributes
group-specific regrets more extensively
than moderate penalties. Its fitted dual
weights generally involve several age groups,
and the identity of the worst-represented
test group varies across splits.

In this application, $\lambda=0.4$ achieves a more favorable median out-of-sample reconstruction--fairness trade-off than the minimax endpoint under both weighting schemes, yielding a larger median reduction in test worst-group regret at a smaller median increase in aggregate reconstruction risk. This occurs despite the minimax endpoint achieving the greatest reduction in training worst-group regret. These findings suggest that, in finite samples, retaining aggregate reconstruction risk in the objective can improve not only out-of-sample reconstruction accuracy but also fairness on held-out observations, relative to minimizing training worst-group regret alone.

\end{document}